\documentclass[oneside, a4paper,11pt,reqno]{amsart}
\usepackage{amsthm}
\usepackage{graphics,dsfont}
\usepackage[T1]{fontenc}

\usepackage{euscript}

\usepackage{amsmath,amssymb}

\usepackage{dsfont}
\usepackage{graphicx,psfrag}

\usepackage{multicol}

\newcommand{\mysection}{\setcounter{equation}{0} \section}

\newtheorem{theo}{Theorem}[section]

\newtheorem{lemma}[theo]{Lemma}

\newtheorem{cor}[theo]{Corollary}

\newtheorem{fact}[theo]{Fact}

\def\un{ \mathds{1}}

\def\un{ {\bf 1}}

\def\a{\alpha}

\def\b{\beta}

\def\bti{\widetilde{\b}}
\def\bhat{\widehat{\b}}
\def\Bhat{\widehat{B}}
\def\Bt{\widetilde{\beta}}

\def\D{\Delta}
\def\d{\delta}
\def\e{\varepsilon}

\def\E{\mathbb{E}}
\def\o{\omega}

\def\g{\gamma}

\def\l{\lambda}
\def\L{\mathcal{L}}
\def\Lb{\widehat{L}}
\def\Ib{\widehat{I}}
\def\N{\mathbb{N}}
\def\P{\mathbb{P}}
\def\p{\psi}

\def\R{\mathbb{R}}

\def\Rt{\widetilde{R}}
\def\r{\rho}
\def\s{\sigma}

\def\z{\zeta}

\def\el{\overset{\L}{=}}

\def\EA{E_5}
\def\EB{E_6}
\def\EG{E_{10}}
\def\EGb{E_{15}}
\def\EH{E_{16}}
\def\EK{E_{12}}
\def\EL{E_2}
\def\ELb{E_3}
\def\EP{E_{17}}
\def\EQ{E_{18}}
\def\ER{E_{11}}
\def\ET{E_{13}}
\def\EU{E_7}
\def\EV{E_8}
\def\EW{E_1}
\def\fn{E_4}
\def\EY{E_{14}}

\def\v1{v_{1}}

\def\bc{Borel--Cantelli }

\def\j1{J_{\Lambda}^<}

\def\k{\kappa}
\def\wk{W_{\k}}

\def\po{P_{\wk}}

\def\lx{L_X}
\def\lB{L_B}
\def\lm{\lx^*}
\def\lE{L^*}
\def\lp{\lx^{+}}
\def\LA{\widetilde{L}_0}

\def\lneg{L_X^{*-}(+\infty)}

\def\lbar{\overline{L}}
\def\hbar{\overline{H}}

\def\rt{\widetilde{R}_{2+2\k}}
\def\rb{\widehat{R}_{2+2\k}}
\def\ty{\Lambda_Y}

\title[Maximum local time of a diffusion process in
a drifted Brownian potential]
{The maximum of the local time of a diffusion process\\ in
a drifted Brownian potential}

\date{\today}

\author{Alexis Devulder}
\address{
Laboratoire de Mathématiques de Versailles, UVSQ, CNRS, Université Paris-Saclay, 78035 Versailles, France
}
\email{devulder@math.uvsq.fr }

\keywords{Random environment,
diffusion in a random potential, maximum local time, L\'evy class, law of the iterated logarithm.\\
This research was partially supported by the french ANR project MEMEMO2 2010 BLAN 0125.}

\begin{document}


\begin{abstract}
We consider a one-dimensional diffusion process $X$ in a $(-\kappa/2)$-drifted
Brownian potential for $\kappa\neq 0$. We are interested in the maximum of its
local time, and study its almost sure asymptotic behaviour,
which is proved to be different from the behaviour of the
maximum local time of the transient random walk in random
environment. We also obtain the convergence in law
of the maximum local time of $X$ under the annealed law after suitable renormalization when $\kappa \geq 1$.
Moreover, we characterize all the upper and lower
classes for the hitting times of $X$, in the sense of Paul L\'evy,
and provide laws of the iterated logarithm for the diffusion $X$ itself. To this aim, we use annealed technics.
\end{abstract}

\maketitle





\mysection{Introduction}
\subsection{Presentation of the model}

We consider a diffusion process in random environment, defined as
follows. For $\k\in\R$, we introduce the random potential
\begin{equation}\label{eqdefinitionWk}
\wk(x):=W(x)-\frac{\k}{2}x,\qquad x\in\R,
\end{equation}
where $(W(x),\ x\in\R)$ is a standard two-sided
Brownian motion.
Informally, a diffusion process ($X(t),\ t\geq
0)$ in the random potential $W_{\k}$ is defined by
$$
\left\{
\begin{array}{l}
\textnormal{d}X(t)=\textnormal{d}\b(t)-\frac{1}{2}\wk'(X(t))
\textnormal{d}t,\\
X(0)=0,
\end{array}\right.
$$
where $(\b(t),\ t\geq 0)$ is a Brownian motion
independent of $W$. More rigorously, $(X(t),\ t\geq 0)$ is a diffusion process
such that $X(0)=0$, and whose conditional generator given
$W_{\k}$ is
$$
\frac{1}{2}e^{W_{\k}(x)}\frac{\text{d}}{\text{d} x}
\left(e^{-W_{\k}(x)}\frac{\text{d}}{\text{d} x}\right).
$$
Let $P$ be the probability measure associated to $\wk$.
We denote by $\po^{}$ the law of $X$ conditionally on the
environment $\wk$, and call it the {\it quenched law}. We also define
the {\it annealed law} $\P$ as follows:
$$
    \P(\cdot )
:=
    \int \po^{}(\cdot )P(\wk\in\text{d}\o).
$$
Notice in particular that $X$ is a Markov process under $\po^{}$, but not under $\P$.
Such a diffusion can also be constructed from a Brownian motion through (random) changes of time and scale
(see \eqref{eqChangeOfTimeAndScale} below).
This diffusion $X$, introduced by Schumacher \cite{Papier_3_Sch1} and
Brox \cite{Papier_3_B3}, is generally considered as the
continuous time analogue of random walks in random environment
(RWRE), which have many applications in physics and biology (see
e.g.~Le Doussal et al.~\cite{Papier_3_LMF1}); for an account of
general properties of RWRE, we refer to R\'ev\'esz
\cite{Papier_3_R1} and Zeitouni \cite{Papier_3_Z1}. This
diffusion has been studied for example by Kawazu
and Tanaka \cite{Papier_3_KT1}, see Theorem \ref{theoKawazuTanaka} below, later improved by Hu, Shi and Yor~\cite{Papier_3_HSY1}.
Large deviations results are proved in Taleb \cite{Papier_3_Ta} and Talet \cite{Papier_3_Ta2}
(see also Devulder \cite{Devulder2} for some properties of the rate function),
and moderate deviations are given by Hu and Shi \cite{Papier_3_HS2} in the recurrent case,
and by Faraud \cite{Faraud} in the transient case.
A localization result and an aging theorem are provided by Andreoletti and Devulder \cite{Andreoletti_Devulder}
in the case $0<\kappa<1$.
For a relation between RWRE and the
diffusion $X$, see e.g. Shi \cite{Papier_3_S2}.
See also Carmona \cite{Papier_3_Carmona}, Cheliotis \cite{cheliotisAsymmetric},
Mathieu \cite{Papier_3_Mathieu}, Singh \cite{SinghIHP}, \cite{Singh} and Tanaka \cite{tanaka}
for diffusions in other potentials.

In this paper, we are interested in the transient case, that is,
we suppose $\k\neq 0$. If $X$ is a diffusion in the random
potential $\wk$, then $-X$ is a diffusion in the random potential
$(\wk(-x),\ x\in\R)$ which has the same law as $(W_{-\k}(x),\
x\in\R)$. Hence we may assume without loss of generality that
$\k>0$. In this case, $X(t)\to_{t\to+\infty}+\infty$ $\P$--almost
surely.

Our goal is to study the
asymptotics of the
maximum
of the local time of $X$. Corresponding problems for RWRE have attracted much attention, and have
been studied, for example, in R\'ev\'esz (\cite{Papier_3_R1},
Chapter 29), Shi \cite{Papier_3_S4},
Gantert et al. \cite{GPS}, \cite{Papier_3_GS1},
Hu et al. \cite{Papier_3_HS1},
Dembo et al. \cite{DGPS}
and Andreoletti (\cite{Papier_3_A}, see also \cite{Papier_3_An2}).
Moreover the local time of such processes in random environment plays an important role in estimation problems
(see e.g. Comets et al. \cite{CometsEstimation}), in persistence (see Devulder \cite{Devulder_Persistence})
and in the study of processes in random scenery (see  Zindy \cite{Zindy}).

\subsection{Maximum local time}

We denote by $(\lx(t,x),\ t\ge 0,\ x\in\R)$ the local time of $X$,
which is the jointly continuous process satisfying, for any
positive measurable function $f$,
\begin{equation}
    \label{eqDefinitiondeWk}
    \int_0^t f(X(s))\text{d}s=\int_{-\infty}^{+\infty}
    f(x)\lx(t,x)\text{d}x,\qquad t\geq 0.
\end{equation}
The existence of such a process was proved by Hu and Shi (\cite{Papier_3_HS1}, eq. (2.6)); see \eqref{tempsLocal} below for an expression of $L_X$.
We are interested in the {\it maximum local time} of $X$ at
time $t$, defined as
$$
\lm(t):=\sup\limits _{x\in\R}\lx(t,x),\qquad t\geq 0.
$$
In the recurrent case $\k=0$, Hu and Shi \cite{Papier_3_HS1} first proved that for any $x\in\R$,
$$
\frac{\log L_X(t,x)}{\log t}\overset{\L}{\longrightarrow} U \wedge \hat{U},
$$
where $U$ and $\hat{U}$ are two independent random variables uniformly distributed in $[0,1]$,
and $``\overset{\L}{\longrightarrow}"$ denotes convergence in law under the annealed law $\P$.
Moreover, thoughout the paper, $\log$ denotes the natural logarithm.
The limit law of $\lm(t)$, suitably renormalized,
is determined by Andreoletti and Diel \cite{Papier_3_AnDi} when $\k=0$:
\begin{equation}\label{eqMaxlocRec}
    \frac{\lm(t)}{t}
\overset{\L}{\longrightarrow}
    \Big(
    \int_{-\infty}^\infty e^{-\widetilde{W}(x)}\text{d} x
    \Big)^{-1},
\end{equation}
where $\big(\widetilde{W}(x),\ x\in \R\big)$ is a two-sided Brownian motion conditioned to stay positive.
Furthermore, Shi \cite{Papier_3_S4} proved the following surprising result: $\P$-almost surely when $\k=0$,
\begin{equation}\label{eqShiLimsupRecurrent}
    \limsup_{t\to+\infty} \lm(t)/(t\log\log\log t)
\geq
    1/32.
\end{equation}
The question whether this is the good renormalization remained open during 13 years,
until Diel \cite{Diel} gave a positive answer to this question. He proved indeed that  in this recurrent case $\kappa=0$,
$$
    \limsup_{t\to+\infty} \lm(t)/(t\log\log\log t)
    \leq e^2/2,
\qquad
    j_0^2/64
    \leq
    \liminf_{t\to+\infty} \lm(t)/[ t/(\log\log\log t)]
    \leq
    e^2\pi^2/4
$$
$\P$-almost surely,
where $j_0$ is the smallest strictly positive root of the Bessel function $J_0$.
Moreover, the convergence in law \eqref{eqMaxlocRec} is extended to the case of stable L\'évy environment by Diel and Voisin \cite{DielVoisin}.
Finally, related questions about favorite sites, that is, locations in which the local time is maximum at time $t$, are
considered by
Hu and Shi \cite{HuShiProblemMostVisitedSite},
Cheliotis \cite{CheliotisFavorite},
and Andreoletti et al. \cite{Andreoletti_Devulder_Vechambre}.



\subsection{Results}
We define the first hitting time of $r$ by $X$ as follows:
\begin{equation}
    \label{H}
    H(r):=\inf\{t\geq 0,\quad X(t)>r\},\qquad r\geq 0.
\end{equation}
We recall that there are three different regimes for $H$ in the transient case $\k>0$:

\begin{theo}(Kawazu and Tanaka, \cite{Papier_3_KT1})
 \label{theoKawazuTanaka}
 When $r$ tends to infinity,
 \begin{eqnarray}
     \frac{H(r)}{r^{1/\k}}
& \overset{\L}{\longrightarrow}&
    c_0 S_{\k}^{ca} ,
\qquad
     0<\k<1,
\label{KawazuTanakakappaInferieur1}
\\
     \frac{H(r)}{r\log r}
& \overset{P.}{\longrightarrow}&
    4,
\qquad
    \k=1,
\label{KawazuTanakakappa1}\\
     \frac{H(r)}{r}
& \overset{a.s.}{\longrightarrow}&
    \frac{4}{\k-1},
\qquad
     \k>1,
\label{KawazuTanakakappaSuperieur1}
\end{eqnarray}
where $c_0=c_0(\k)>0$ is a finite constant, the symbols
$``\overset{\L}{\longrightarrow}"$,
$``\overset{P.}{\longrightarrow}"$ and
$``\overset{a.s.}{\longrightarrow}"$ denote respectively
convergence in law, in probability and almost sure convergence,
with respect to the annealed probability $\P$. Moreover, for $0<\k<1$,
$S_\k^{ca}$ is a completely asymmetric stable variable of index
$\k$, and is a positive variable (see
\eqref{TransfoLaplace1} for its characteristic function).
\end{theo}




The asymptotics of the maximum local time $\lm(t)$ heavily depend on the value of
$\k$. We start with the upper asymptotics of $L_X^*(t)$:

\begin{theo}
 \label{theoLimsupDifferentRWRE}
 If $0<\k<1$, then
 $$
  \limsup_{t\to+\infty} \frac{\lm(t)}{t} = +\infty \qquad \P
  \text{--a.s.}
 $$
\end{theo}

Theorem \ref{theoLimsupDifferentRWRE} tells us that in the case
$0<\k<1$, the maximum local time of $X$ has a completely different
behaviour from the maximum local time of RWRE (the latter is
trivially bounded by $t/2$ for any positive integer $t$, for
example). Such a peculiar phenomenon has already been observed (see \eqref{eqShiLimsupRecurrent})
by Shi \cite{Papier_3_S4} in the recurrent case, and is even more
surprising here since $X$ is transient.


Theorem \ref{theoLimsupKsup1} gives, in the case $\k>1$, an
integral test which completely characterizes the upper functions
of $\lm(t)$, in the sense of Paul L\'evy.

\begin{theo}
 \label{theoLimsupKsup1}
 Let $a(\cdot)$ be a positive nondecreasing function. If $\k>1$, then
 $$
  \sum_{n=1}^{\infty}\frac{1}{n a(n)} \
  \left\{
  \begin{array}{l}
      <+\infty
      \\
      =+\infty
  \end{array}
  \right.
  \Longleftrightarrow \limsup_{t\to\infty} \frac{\lm(t)}{[t
  a(t)]^{1/\k}} =
  \left\{
  \begin{array}{l}
      0
      \\
      +\infty
  \end{array}
  \right.
  \qquad\P\text{--a.s.}
 $$
\end{theo}

This is in agreement with a result of Gantert and Shi \cite{Papier_3_GS1} for RWRE.
We notice in particular that
$\limsup_{t\to+\infty} \lm(t)/t$ is almost surely $+\infty$ when $0<\k<1$ by Theorem \ref{theoLimsupDifferentRWRE},
whereas it is $0$ when $\k>1$ by Theorem \ref{theoLimsupKsup1}. We have not been able to prove whether  $\limsup_{t\to+\infty} \lm(t)/t$ is infinite in
the very delicate case $\k=1$, since a  proof similar to that of Theorem \ref{theoLimsupDifferentRWRE}
just shows that it is greater than a positive deterministic constant (see Remark page \pageref{InegLimsupCasK1} for more details).


We now turn to the lower asymptotics of $\lm(t)$.

\begin{theo}
 \label{theoLiminfKsup1}
 We have
 \begin{eqnarray*}
     \liminf_{t\to\infty}\frac{\lm(t)}{t/\log\log t}
  & \leq & \k^2c_1(\k)\qquad  \P\text{--a.s.}\quad \textnormal{if
     } 0<\k<1,
     \\
     \liminf_{t\to\infty}\frac{\lm(t)}{t/[(\log t)\log\log t]}
  & \leq & \frac{1}{2}\qquad  \P\text{--a.s.}\quad
     \textnormal{if } \k=1,
     \\
     \liminf_{t\to\infty}\frac{\lm(t)}{(t/\log\log t)^{1/\k}}
  & = & 4\left(\frac{(\k-1)\k^2}{8}\right)^{1/\k}\qquad
     \P\text{--a.s.}\quad \textnormal{if }  \k>1,
 \end{eqnarray*}
where $c_1(\kappa)$ is defined in \eqref{SchilderKinf1}.
\end{theo}

%

\bigskip

\begin{theo}
 \label{theoLiminfKinf1}
 We have, for any $\e>0$,
 \begin{eqnarray*}
     \liminf_{t\to\infty} \frac{\lm(t)}{t/[(\log
     t)^{1/\k} (\log\log t)^{(2/\k)+\e} ]}
  &=& +\infty \qquad  \P\text{--a.s.}\quad \textnormal{if }
     0<\k\leq 1.
 \end{eqnarray*}
\end{theo}

\bigskip
 In the case $0<\k\leq 1$, Theorems \ref{theoLiminfKsup1}
and \ref{theoLiminfKinf1} give different bounds, for technical
reasons.

We also get the convergence in law under the annealed law $\P$ of $\lm(t)$, suitably renormalized, when $\k\geq 1$:

\begin{theo}\label{theoLawLm}
We have as $t\to+\infty$, under the annealed law $\P$,
 \begin{eqnarray*}
     \frac{\lm(t)}{t/\log t}
& \overset{\L}{\longrightarrow} &
    \frac{1}{2\mathcal{E}}
    \qquad \textnormal{if }     \k= 1,
\\
     \frac{\lm(t)}{t^{1/\k}}
& \overset{\L}{\longrightarrow} &
    4\big[\k^2(\k-1)/8\big]^{1/\k}\mathcal{E}^{-1/\k}
    \qquad \textnormal{if }     \k> 1,
 \end{eqnarray*}
where
$\mathcal{E}$ denotes an exponential variable with mean $1$.
\end{theo}

We notice that in the previous theorem, the case $0<\kappa<1$ is lacking. Indeed, we did not succeed in obtaining it with the
annealed technics of the present paper, because
due to \eqref{KawazuTanakakappaInferieur1},
$H(r)$ suitably renormalized converges in law but does not converge in probability to a positive constant in this case.
This is why we used quenched technics in Andreoletti et al. \cite{Andreoletti_Devulder_Vechambre} to prove that $\lm(t)/t$
converges in law under $\P$ as $t\to+\infty$ when $0<\kappa<1$. To this aim, we used and extended to local time the quenched tools developed
in Andreoletti et al. \cite{Andreoletti_Devulder} to get the localization of $X$ in this case $0<\k<1$,
combined with some additional  tools such as  two dimensional L\'evy processes and convergence in Skorokhod topology.

So, Theorem \ref{theoLawLm} completes the results of
\cite{Andreoletti_Devulder_Vechambre} and  \cite{Papier_3_AnDi}  (see our \eqref{eqMaxlocRec}), that is,
  these 3 results give the convergence in law of $\lm(t)$ suitably renormalized
for any value of $\k\in\R$.

In the proof of Theorems \ref{theoLimsupDifferentRWRE}, \ref{theoLiminfKsup1}
and \ref{theoLiminfKinf1}, we will
frequently need to use the almost sure asymptotics of the
first hitting times $H(\cdot)$. In view of the last part \eqref{KawazuTanakakappaSuperieur1} of
Theorem \ref{theoKawazuTanaka}, we only need to study the case
$\k \in (0,1]$.

\begin{theo}
 \label{ThUpperLevyKinf1}
 Let $a(\cdot)$ be a positive nondecreasing function.
 If $0<\k<1$, then
 $$
  \sum_{n=1}^{\infty}\frac{1}{n a(n)} \
  \left\{
  \begin{array}{l}
      <+\infty
      \\
      =+\infty
  \end{array}
  \right.
  \Longleftrightarrow \limsup_{r\to\infty} \frac{H(r)}{[r
  a(r)]^{1/\k}}=
  \left\{
  \begin{array}{l}
      0
      \\
      +\infty
  \end{array}
  \right.
  \qquad\P\text{--a.s.}
 $$
 If $\k=1$, the statement holds under the additional
 assumption that $\limsup_{r\to+\infty} (\log r)/a(r)<\infty$.
\end{theo}

\begin{theo}
 \label{thLowerLevyKinf1}
We have ($\Gamma$ denotes the usual gamma function)
\begin{eqnarray}
    \liminf_{r\to+\infty} \frac{H(r)}{r^{1/\k} / (\log\log r)^{(1/\k)-1} }
& = &
    \frac{8 \k [\pi \k]^{1/\k}  (1-\k)^{\frac{1-\k}{\k}}}{\big[2 \Gamma^2(\k)\sin(\pi\k)\big]^{1/\k}}
    =:
    c_2(\k)
\quad
    \P \text{--a.s. if }\,
    0<\kappa<1,~~\hphantom{aa}
\label{Maine3}
\\
  \liminf_{r\to+\infty} \frac{H(r)}{r\log r} & = & 4 \qquad
  \P \text{--a.s. if }  \k=1.
\label{Maine4}
\end{eqnarray}
\end{theo}

\medskip

The following corollary follows immediately from Theorem \ref{ThUpperLevyKinf1}
and gives a negative answer to a question raised in  Hu, Shi and Yor (\cite{Papier_3_HSY1}, Remark 1.3 p. 3917):

\begin{cor}
The convergence in probability $H(r)/(r\log r) \to 4$ in Theorem
\ref{theoKawazuTanaka} in the case $\k=1$ cannot be strengthened
into an almost sure convergence.
\end{cor}


We observe that in the case $0<\k<1$, the process $H(\cdot)$ has
the same almost sure asymptotics as $\k$--stable subordinators (see
Bertoin \cite{Papier_3_Ber1} p.~92).

Finally, define $\log_1:=\log$ and
$\log_k:=\log_{k-1}\circ\log$ for $k>1$.
Theorems \ref{ThUpperLevyKinf1} and \ref{thLowerLevyKinf1}, and the fact
that $X(t)$ is not very far from $\sup_{0\leq s \leq t} X(s)$
(see Lemma \ref{LemmaNotTooFarOnTheLeft} below) lead to

\begin{cor}\label{CorrolaryLILofX}
Recall that $c_2(\k)$ is defined in \eqref{Maine3}.
We have for $k\in\N^*$,
\begin{eqnarray}
    \limsup_{t\to\infty}\frac{X(t)}{t^\k(\log\log  t)^{1-\k}}
& = &
    \frac{2\Gamma^2(\k)\sin(\pi\k)}{\pi 8^\k  \k^{\k+1}(1-\k)^{1-\k}}
=
    \frac{1}{[c_2(\k)]^\k}
\qquad
    \P\text{--a.s. if  } 0<\k<1,\hphantom{aa}
\label{eqLimsupXcaskentre0et1}
\\
    \limsup_{t\to\infty}\frac{X(t)}{t/\log t}
& = &
    \frac{1}{4}
\qquad
    \P\text{--a.s. if  } \k=1,
\label{eqLimsupXcaskegal1}
\end{eqnarray}
\begin{equation}
   \left\{
  \begin{array}{l}
      \alpha
      \leq 1
      \\
      \alpha
      >1
  \end{array}
  \right.
    \Longleftrightarrow
\liminf_{t\to+\infty}\frac{X(t)}{t^\k/[(\log t)\dots(\log_{k-1} t)(\log_k t)^\alpha]}
   \left\{
  \begin{array}{l}
       = 0
      \\
      =+\infty
  \end{array}
  \right.
\qquad
    \P\text{--a.s. if  } 0<\k\leq 1,
\label{eqLiminfXcaskentre0et1}
\end{equation}
where for $k=1$, $(\log t)\dots(\log_{k-1} t)=1$ by convention.
These results remain true if we replace $X(t)$ by $\sup_{0\leq s \leq t} X(s)$.
\end{cor}

Corresponding results in the recurrent case $\k=0$ are proved by Hu et al. \cite{HuShiLimits},
extended later by Singh \cite{SinghIHP} to some asymptotically stable potentials
and following results of Deheuvels et al. \cite{{deheuvelsRevesz}} for Sinai's walk.


\smallskip

Our proof hinges upon stochastic calculus. In particular, one key ingredient of the proofs of
Theorems \ref{theoLimsupDifferentRWRE}--\ref{thLowerLevyKinf1} is
an approximation of the joint law of the hitting time $H[F(r)]$ of $F(r)\approx r$ by $X$
and the maximum local time $\lm[H(F(r))]$ of $X$ at this time, stated in Lemma \ref{lemmaApproxLxLi}, and proved in Section \ref{sectannexe}.
Another important tool is a modification of the Borel-Cantelli lemma, stated in Lemma \ref{lemmaBorelCantelli},
which, loosely speaking, says that one can chop the real half line $[0,\infty)$ into regions in which the diffusion $X$ behaves in an "independent" way.

\smallskip

The rest of the paper is organized as follows. In Section
\ref{sectSimpl}, we give some preliminaries on local time and Bessel
processes. We present in Section~\ref{sectMainEstimates} some
estimates which will be needed later on; the proof of
one
key estimate (Lemma \ref{lemmaApproxLxLi}) is
postponed until Section \ref{sectannexe}. Section
\ref{SectEtudeLmdeHdeR} is devoted to the study of the almost sure asymptotics of $\lm[H(r)]$,
stated in Theorems \ref{lemmaLimsupL0} and \ref{lemmaLiminfL0}.
In Section \ref{SectEtudeHdeR}, we study
the L\'evy classes for the hitting times $H(r)$ and prove
Theorems \ref{ThUpperLevyKinf1} and \ref{thLowerLevyKinf1} and Corollary \ref{CorrolaryLILofX}. In
Section \ref{SectJointStudy}, we study $\lm[H(r)]/H(r)$ and
prove Theorems
\ref{theoLimsupDifferentRWRE}--\ref{theoLawLm}.
Section \ref{sectannexe} is devoted to the proof of Lemma
\ref{lemmaApproxLxLi}.
Finally, we prove in Section \ref{SectProofLemmasFin}
some lemmas dealing with Bessel processes, Jacobi processes and
Brownian motion.

Throughout the paper, the letter $c$ with a subscript denotes
constants that are finite and positive.



\mysection{Some preliminaries}

\subsection{Preliminaries on local time and Bessel processes}
 \label{sectSimpl}

We first define, for any Brownian motion $(B(t),\ t\geq 0)$ and $r>0$,
the hitting time
$$
    \s_B(r)
:=
    \inf\{t>0,\ B(t)=r\}.
$$
Moreover, we denote by $(\lB(t,x),\ t\geq 0,\ x\in\R)$
the local time of $B$, i.e., the jointly continuous process
satisfying
$
    \int_0^t f(B(s))\text{d}s
=
    \int_{-\infty}^{+\infty} f(x)\lB(t,x)\text{d}x
$ for any positive measurable function
$f$. We define the inverse local time of $B$ at $0$ as
$$
    \tau_ {B}(a)
:=
    \inf\{t\geq 0,\ L_{B}(t,0)> a\},\qquad a>0.
$$
Furthermore, for any $\d\in[0,\infty)$ and $x\in[0,\infty)$, the
unique strong solution of the stochastic differential equation
$$
    Z(t)
=
    x+2\int_0^t\sqrt{Z(s)}\text{d}\b(s)+\d t,
$$
where $(\b(s),s\geq 0)$ is a (one dimensional) Brownian motion,    is
called a $\d$--dimensional squared Bessel process starting from
$x$.
A Bessel process with dimension $\d$ (or equivalently with order $\delta/2-1$) starting from $x\geq 0$ is
defined as the (nonnegative) square root of a $\d$--dimensional
squared Bessel process starting from $x^2$
\big(see e.g. Borodin et al. \cite{BorodinSalminem}, 39 p. 73
for a more general definition
as a linear diffusion with generator
$\frac{1}{2} \frac{d^2}{d x^2}{}+\frac{\delta-1}{2x}\frac{d}{dx}$ for every $\delta\in\R$;
 see also G{\"o}ing-Jaeschke et al. \cite{Going-Yor} definition 3 p. 329\big). We recall some
important results.

\begin{fact}(first Ray--Knight theorem)
 \label{FactRN1}
 Consider $r>0$ and a Brownian motion $(B(t),\ t\geq 0)$. The
 process $(L_B(\s_B(r),r-x),\ x\geq 0)$ is a continuous
 inhomogeneous Markov process, starting from~$0$. It is a
 $2$--dimensional squared Bessel process for $x\in[0,r]$ and a
 $0$--dimensional squared Bessel process for $x\geq r$.
\end{fact}

\begin{fact}(second Ray--Knight theorem)
 \label{FactRN2}
 Fix $r>0$, and let $(B(t),\ t\geq 0)$ be a Brownian motion. The
 process $(L_B(\tau_B(r),x),\ x\geq 0)$ is a\, $0$--dimensional
 squared Bessel process starting from $r$.
\end{fact}

See e.g.\ Revuz and Yor (\cite{Papier_3_RY3}, chap.~XI) for more
details about Ray--Knight theorems and Bessel processes.
Following the method used by Hu et al. (\cite{Papier_3_HSY1}, see eq. (3.8)), we also
need the following well known result:

\begin{fact}(Lamperti representation theorem,
 see Yor \cite{Yor2} eq. (2.e))
\label{FactLamperti}
Consider $\wk(x)=W(x)-\k x/2$ as in \eqref{eqdefinitionWk} with $\kappa>0$,
where $(W(x),\ x\geq 0)$ is a Brownian motion. There exists a
$(2-2\k)$--dimensional Bessel process $(\rho(t),\ t\geq 0)$,
starting from $\rho(0)=2$, such that
$\exp[\wk(t)/2] = \rho(A(t))/2$ for all $t\geq 0$,
where $A(r):=\int_0^r e^{W_\k(s)}\textnormal{d}s$, $r\geq 0$.
\end{fact}

We also recall the following extension to Bessel processes of
Williams' time reversal theorem (see Yor \cite{Papier_3_Y11},
p.~80; see also G{\"o}ing-Jaeschke et al. \cite{Going-Yor} eq. (34)).

\begin{fact}
 \label{factWilliam}
 One has, for $\d<2$,
 $$
  (R_\d(T_0-s),\ s\leq T_0)\el(R_{4-\d}(s),\ s\leq\g_a),
 $$
 where $\el$ denotes equality in law,
 $(R_\d(s),\ s\geq 0)$ denotes a $\d$--dimensional Bessel
 process starting from $a>0$, $T_0:=\inf\{s\geq 0,\ R_\d(s)=0\}$,
 $(R_{4-\d}(s),\ s\geq 0)$ is a $(4-\d)$--dimensional Bessel
 process starting from $0$, and $\g_a:=\sup\{s\geq 0,\
 R_{4-\d}(s)=a\}$.
\end{fact}

Let $S_\k^{ca}$ be a (positive) completely asymmetric stable
variable of index $\k$ for $0<\kappa<1$, and $C_8^{ca}$ a (positive) completely asymmetric
Cauchy variable of parameter $8$. Their characteristic
functions are given by:
\begin{equation}
    \E e^{itS_\k^{ca}}
  =  \exp\left[-|t|^\k\left(1-i\, \textnormal{sgn}(t)
    \tan\big(\frac{\pi\k}{2}\big)\right) \right],
    \label{TransfoLaplace1}
\quad\
    \E e^{itC_8^{ca}}
  =  \exp\Big[- 8 \Big(|t|+i t \frac{2}{\pi} \log|t| \Big)
    \Big].
\end{equation}
Throughout the paper, we set $\l:= 4(1+\k)$.
If $(B(t),\ t\geq 0)$ denotes, as before,  a Brownian motion,
we introduce
\begin{eqnarray}
    K_\b({\k})
 & := & \int_0^{+\infty}x^{1/\k-2}L_{\b}(\tau_\b(\l),x)\text{d}x,
    \qquad 0<\k<1,
    \label{e3p18}
    \\
    C_\b
 & := & \int_0^1\frac{L_\b(\tau_\b(8),x)-8}{x}\text{d}x+
    \int_1^{+\infty}\frac{L_\b(\tau_\b(8),x)}{x}\text{d}x.
    \label{e2p18}
\end{eqnarray}
We have the following equalities in law:

\begin{fact}(Biane and Yor \cite{Papier_3_BY2})\label{FactBianeYor}
For
$0<\k<1$,
\begin{equation*}
    C_\beta
 \el
    8c_3 +(\pi/2)C_8^{ca},
\qquad
    K_\beta(\k)
 \el
    \big(\k^{2-1/\k} c_4(\k)/4\big)S_\k^{ca},
\end{equation*}
where $c_3>0$ denotes an unimportant constant, and
\begin{equation}
    \label{DefPsi}
    \psi(\k) := \left(\frac{\pi \k}{4\Gamma^2(\k)\sin(\pi
    \k/2)}\right)^{1/\k},
\qquad
    c_4(\k)
:= 8
    \psi(\k)\l^{1/\k}\k^{-1/\k}.
\end{equation}
\end{fact}
This fact is proved in (Biane and Yor \cite{Papier_3_BY2});
the identity in law related to $C_\beta$
is given in its paragraph (4.3.2) pp 64-66 and
the one related to $K_\beta(\k)$
follows from its (1.a) p. 24.

Finally, the first Ray-Knight theorem leads to the following formula. For
$v>0$ and $y>0$,
\begin{equation}
    \label{borodin}
    \P\bigg(\sup_{0\leq s\leq \tau_{\beta}(v)} \beta(s)<y\bigg)
=
    \P\big[L_\b(\sigma_{\beta}(y),0)>v\big]
=
    \P\big(R_2^2(y)>v\big)
=
    \exp\left(-\frac{v}{2y}\right),
\end{equation}
where $(R_2(s),\ s\geq 0)$ is a  $2$--dimensional Bessel process
starting from~$0$.



\subsection{Some preliminaries on the diffusion}\label{sectMainEstimates}

We assume in the rest of the paper that $\k>0$, and so $X$ is a.s. transient to the right.
We start by introducing
$$
    A(x)
:=
    \int_0^x e^{\wk(y)}\text{d}y,
\quad
    x\in\R,
\qquad
    A_\infty
:=
    \int_0^\infty e^{\wk(y)}\text{d}y
<
    \infty
\text{ a.s.}
$$
We recall that $A$ is a scale function of $X$ under the quenched law $\po^{}$ (see e.g. Shi \cite{Papier_3_S2} eq. (2.2)).
That is, if $\po^y$ denotes the law of the diffusion $X$ in the potential $\wk$, starting from $y$ instead of $0$, we have conditionally on the potential $\wk$,
\begin{equation}\label{eqScaleFunctionA}
    \po^y\big[H(z)<H(x)\big]
=
    \big[A(y)-A(x)\big]/\big[A(z)-A(x)\big],
\qquad
    x<y<z.
\end{equation}
We observe that,
since $\k>0$,
$
    A(x)
\to
    A_\infty
<\infty
$
a.s. when $x\to +\infty$.

For technical reasons, we have to introduce the random function
$F$ as follows. Fix $r>0$. Since the function $x\mapsto
A_{\infty}-A(x)=:D(x)$ is almost surely continuous and (strictly)
decreasing and has limits $+\infty$ and $0$ respectively on $-\infty$ and $+\infty$, there exists a unique
$F(r)\in\R$, depending only on the process $\wk$, such that
\begin{equation}
    A_{\infty}-A(F(r))=\exp(-\k r/2)=:\d(r).
    \label{e1p5}
\end{equation}
Our first estimate describes how close
$F(r)$ is to $r$, for large $r$.

\begin{lemma}
 \label{lemmaEncadrementFdeR}
 Let $\k>0$ and  $0<\d_0<1/2$. Define for $r>0$,
 \begin{equation}
    \label{eqEncadrementFdeR}
    \EW(r)
 :=
    \big\{\big(1-5r^{-\d_0}/\k\big)r
            \leq
            F(r)
            \leq
            \big(1+5 r^{-\d_0}/\k\big)r\big\}.
 \end{equation}
 Then for all large $r$,
 \begin{equation}
    \label{eqProbaEW}
    \P\big[\EW(r)^c\big]
\leq
    \exp\big(-r^{1-2\d_0}\big).
 \end{equation}
 As a consequence, for any $\varepsilon>0$, we have, almost
 surely, for all large $r$,
 \begin{equation}
     \label{F}
     (1-\varepsilon) r \le F(r) \le (1+\varepsilon)r.
 \end{equation}
\end{lemma}

\bigskip

\noindent{\bf Proof of Lemma \ref{lemmaEncadrementFdeR}.}
Let $0<\d_0<1/2$, and fix $r>0$. We have
\begin{equation}\label{anest}
    \P\big[\EW(r)^c\big]
\leq
    \P\big[F(r)< (1-5r^{-\d_0}/\k)r\big]
    +\P\big[F(r)> (1+5r^{-\d_0}/\k)r\big].
\end{equation}
Define $s_\pm:=(1 \pm 5r^{-\d_0}/\k)r$, and 
$A_{\infty}^{(s)}:=\int_s^{\infty}\exp(\wk(u)-\wk(s))\text{d}u$ for $s\geq 0$.
Observe that $D$ is stricly decreasing, $D(F(r))=e^{-\k r/2}$ and that
$D(s_\pm)=A_{\infty}^{(s_\pm)}\exp(\wk(s_\pm))$.
Consequently,
\begin{equation*}
    \P\big[F(r)<(1-5r^{-\d_0}/\k)r\big]
\leq
    \P\big[D(F(r))>D(s_-)\big]
=
    \P\big[-\k r/2>\log\big(A_{\infty}^{(s_-)}\big)+\wk(s_-)\big].
\end{equation*}
Moreover,
$A_{\infty}^{(s_\pm)}\el A_{\infty}  \el 2/\g_\k$,
where $\g_\k$ is a gamma variable of parameter $(\k,1)$ (see Dufresne \cite{Papier_3_D1} or Borodin et al. \cite{BorodinSalminem} IV.48 p. 78),
i.e., $\g_\k$ has density $\frac{1}{\Gamma(\k)} e^{-x} x^{\k -1}\un_{\R_+}(x)$. Hence
\begin{eqnarray*}
    \P\big[F(r)< (1-5r^{-\d_0}/\k)r\big]
& \leq &
    \P\big[\log(2/\gamma_k)<-r^{1-\d_0}\big]
    +\P\big[W(s_-)<-3r^{1-\d_0}/2\big]
\\
& \leq &
    2\exp\big(-9r^{1-2\d_0}/8\big),
\end{eqnarray*}
for large $r$,
since $\P[W(1)<-x]\leq e^{-x^2/2}$ for $x\geq 1$.
Similarly, we have for large $r$,
\begin{equation*}
    \P\big[F(r)>s_+\big]
\leq
    \P\big[\log(2/\gamma_k)>r^{1-\d_0}/2\big]
    +\P\big[W(s_+)>2r^{1-\d_0}\big]
\leq
    \exp\big(-9r^{1-2\d_0}/8\big).
\end{equation*}
This yields \eqref{eqProbaEW} in view
of \eqref{anest}.


Then $
    \sum_{n\geq 1} \P\big[\EW(n)^c\big]
<
    \infty
$,
so  \eqref{F} follows from the Borel--Cantelli lemma and the monotonicity of
$F(\cdot)$.
\hfill$\Box$


\medskip

In the rest of the paper, we define, for $\d_1>0$ and any $r>0$,
\begin{equation}
    c_{5}:= 2(\lambda/\k)^{\d_1},
\qquad
    \psi_\pm (r):= 1\pm {\frac{c_{5}}{r^{\d_1}}},
\qquad
    t_\pm (r):= \frac{\k\psi_\pm(r)r}{\l} .
    \label{e1p22}
\end{equation}
Taking $\psi_\pm(r)$ as defined above instead of simply $1\pm\e$ is necessary e.g. in Lemma \ref{lemmaEncadrementQ} below.
Moreover, if $(\b(s),\ s\geq 0 )$ is a Brownian
motion and $v>0$, we define the Brownian motion $(\b_v(s),\ s\geq
0)$ by $\beta_v(s):=(1/v) \beta(v^2 s)$, $s\geq 0$.

We prove in Section \ref{sectannexe} the following
approximation of the joint law of
$
\big(
    \lm[H(F(r))], H(F(r))
\big)
$.

\begin{lemma}
 \label{lemmaApproxLxLi}

Let $\k>0$ and $\e\in (0,1)$. For $\d_1>0$ small enough, there exists $c_{6}>0$ and $\a>0$ such that for $r$ large enough,
there exist a Brownian motion $(\b(t),\ t\ge 0)$ such that the following holds:

\noindent{\bf (i)} Whenever $\k>0$,
we have
$$
\P [ \EL(r)] \ge 1- r^{-\alpha},
$$
where
\begin{eqnarray}
     \label{e3p25}
     \EL(r)
&:=& \left\{(1-\e)\Lb_-(r)\leq \lm[H(F(r))]\leq
     (1+\e)\Lb_+(r)\right\} ,
\\
    \label{e2p22}
    \Lb_\pm(r)
&:=&
    4[\k t_\pm(r) ]^{1/\k}
    \bigg[\sup\limits_{0\leq u \leq \tau_{\beta_{t_\pm(r)}}(\l)} \beta_{t_\pm(r)}(u)\bigg]^{1/\k}
    =4\bigg[\sup\limits_{0\leq u \leq \tau_{\beta}(\l t_\pm(r))} \k\beta(u)\bigg]^{1/\k} .
\end{eqnarray}
{\bf(ii)} If $0<\k\leq 1$,
we have
$$
    \P [ \ELb(r) ]
\ge
    1- r^{-\alpha},
$$
 where, using the notation introduced in \eqref{e3p18} and \eqref{e2p18},
 \begin{eqnarray}
     \label{e4p25}
     \ELb(r)
  &:=& \left\{ (1-\e) \Ib_-(r)\le H(F(r))\le (1+\e)
     \Ib_+(r)\right\} ,
     \\
     \label{eqI''CasKinf1}
        \Ib_\pm(r)
  &:=&
     \left\{
     \begin{array}{lr} 4 \k^{1/\k-2}t_\pm(r)^{1/\k}
       \big[ K_{\b_{t_\pm (r)}}(\k) \pm c_{6}
        t_\pm (r)^{1-1/\k}\big],
&
    \quad 0<\k<1,
\\[2mm]
       4t_\pm(r) \big[C_{\b_{t_\pm (r)}} + 8 \log t_\pm (r) \big],
&
    \quad \k=1 .
\end{array}\right.
\end{eqnarray}
\end{lemma}

\smallskip

Notice in particular that the Brownian motion $\beta$ is the same in {\bf (i)} and {\bf (ii)};
this allows to approximate the law of quantities depending on both $\lm[H(F(r))]$ and $H(F(r))$,
such as $\lm[H(F(r))]/H(F(r))$, which is useful in Section \ref{SectJointStudy}.
This is possible because we kept the random function $F(r)$ in the expressions $\lm[H(F(r))]$ and $H(F(r))$,
in order to have the same Brownian motion $\beta$ in the left hand side and the right hand side of
the inequalities defining $\EL(r)$ and $\ELb(r)$.

The proof of Lemma \ref{lemmaApproxLxLi} is postponed to Section
\ref{sectannexe}.

\medskip

With an abuse of notation, for $z\geq 0$, we denote by
$X\circ{\Theta_{H(z)}}$ the process $(X(H(z)+t)-z,\ t\geq 0)$.
Notice that due to the strong Markov property applied at stopping time $H(z)$ under the quenched law $\po$ , $X\circ{\Theta_{H(z)}}$ is,
 conditionally on $\wk$,  a diffusion in the $(-\k/2)$-drifted Brownian potential
$\wk\circ \Theta_z:=(\wk(x+z)-\wk(z),\ x\in\R)$, starting from $0$. Define
$H_{X\circ\Theta_{H(z)}}(s)=H(z+s)-H(z)$, $s\geq 0$, which is the hitting time of $s$ by $X\circ{\Theta_{H(z)}}$.
In view of \eqref{e1p5}, we also define $F_{\wk\circ \Theta_z}$ by
$
    \int_{F_{\wk\circ \Theta_z}(r)}^\infty
    e^{\wk\circ \Theta_z(u)} \text{d} u
=
    \delta(r)
$,
$r>0$.
That is, $F_{\wk\circ \Theta_z}$ plays the same role for $\wk\circ \Theta_z$ (resp. for $X\circ{\Theta_{H(z)}}$)
as $F$ does for $\wk$ (resp. for $X$).
Similarly, $\lE_{X\circ\Theta_{H(z)}}$ and $(L^*\circ H)_{X\circ\Theta_{H(z)}}$
denote respectively the processes $\lE$ and $L^*\circ H$ for the diffusion
$X\circ\Theta_{H(z)}$, with $(\lE)_X:=\lm$.
The following lemma is a modification of the \bc lemma.


\bigskip

\begin{lemma}
 \label{lemmaBorelCantelli}
 Let $\k>0$,
 $\a>0$, $r_n:=\exp(n^{\a})$ and $Z_n:=\sum_{k=1}^n r_k$ for $n\ge 1$.
 Assume $f$ is a continuous function $(0,+\infty)^2\to\R$ and
 $(\D_n)_{n\ge 1}$ is a sequence of open sets in $\R$ such that
 \begin{equation}
     \label{e1p9}
     \sum_{n\geq 1} \P \big\{ f[(H\circ F)(r_{2n}), (\lm\circ
     H\circ F)(r_{2n}) ] \in \D_n \big\} = +\infty.
 \end{equation}
Then for any $0<\varepsilon<1/2$, $\P$ almost surely, there exist
infinitely many $n$ such that for some
 $t_n \in [(1-\e)r_{2n}, (1+\e)r_{2n}]$,
$$
    f\big[H_{X \circ \Theta_{H(Z_{2n-1})}} (t_n), (\lE\circ
        H)_{X \circ \Theta_{H(Z_{2n-1})}} (t_n) \big]
\in
    \D_n.
$$
The results remain true if $r_n=n^n$ for every $n\geq 1$.
\end{lemma}

\bigskip

\noindent {\bf Proof of Lemma \ref{lemmaBorelCantelli}.} We
divide $\R_+$ into some regions in which the diffusion $X$ will
behave ``independently'', in order to apply the \bc lemma.

To this aim, let $n\ge 1$ and
$$
    \fn(n)
:=
    \bigg\{ \inf_{\{t: \; H(Z_{2n-1})\leq t\leq
        H(Z_{2n}+r_{2n+1}/2)\}} X(t)>Z_{2n-2}+ \frac{1}{2} r_{2n-1}
    \bigg\}.
$$
Define $x_n:=r_{2n-1}/2$.
For any environment, i.e., for any realization of $\wk$, $X$ is a
Markov process under $\po^{}$, and $H(Z_{2n-1})$ is a stopping time. Hence,
$\po^{}(\fn(n)^c)$ is the probability that the diffusion in the potential $\wk$ started at
$Z_{2n-1}$ hits level $Z_{2n-2}+x_n$ before $Z_{2n}+x_{n+1}$, that is
\begin{equation}\label{eqProbaE5}
    \po^{}\big[\fn(n)^c\big]
=
    \Bigg( 1+ \frac{ \int_{Z_{2n-2} +x_n}^{Z_{2n-1}}
        e^{\wk(u)} \text{d}u} { \int_{Z_{2n-1}}^{Z_{2n}+x_{n+1}}
        e^{\wk(u)} \text{d}u}
    \Bigg)^{-1}
\leq \frac{\int_{Z_{2n-1}}^{Z_{2n}+x_{n+1}}
    e^{\wk(u)} \text{d}u}{\int_{Z_{2n-2} +x_n}^{Z_{2n-1}}
    e^{\wk(u)} \text{d}u},
\end{equation}
where we used \eqref{eqScaleFunctionA}.
Observe that $r_{2n-1}-x_n=x_n$ and define for some $0<\e_0<\k/4$,
\begin{equation*}
    \EA(n)
:=
    \bigg\{ \sup_{0\leq u \leq r_{2n-1}-x_n } \left| \wk(u+ Z_{2n-2}+x_n) - \wk(Z_{2n-2} +x_n) +\frac{\k}{2} u \right|
            \leq \e_0(r_{2n-1}-x_n)
            \vphantom{\sup_{0\leq u \le r_{2n-1}-x_n }
              \left| \wk(u+Z_{2n}+x_n) - \wk(Z_{2n}+x_n) +
              \frac{\k}{2} u \right|}
    \bigg\}
\end{equation*}
and
$
    \EB(n)
 := \left\{ \sup_{u\ge 0} [ \wk(u+Z_{2n-1}) - \wk(Z_{2n-1})] \leq v_n \right\}
$, where $v_n:=2(\log n)/\k$.
Since $\sup_{0\leq u\leq x_n}W(u)\el |W(x_n)|$
and $\sup_{x\geq 0} \wk(x)$ has an exponential law of parameter $\k$ (see e.g. Borodin et al. \cite{BorodinSalminem} 1.1.4 (1) p. 251),
we have for large $n$,
\begin{equation}
\label{eqprobaEAbis}
    \P\big[\EA(n)^c\big]
=
    \P\Big(\sup_{0\leq u \leq x_n}|W(u)|>\e_0 x_n\Big)
\leq
    4\exp\Big[- \frac{\e_0^2 x_n}{2}\Big]
\ \  \text{ and } \
    \P\big[\EB(n)^c\big]
=
    \exp(-\k v_n)=n^{-2}.
\end{equation}
Moreover by \eqref{eqProbaE5}, we have for $n$ large enough,
on $\EA(n)\cap\EB(n)$,
%
\begin{eqnarray}
    \po^{}\big[\fn(n)^c\big]
& \le &
    \k\frac{(r_{2n} +x_{n+1} ) \exp[ v_n+ \wk(Z_{2n-1})]}
    {\exp[ \wk( Z_{2n-2} +x_n) -\e_0(r_{2n-1}
    -x_n)]}
\nonumber
\\
& \le &
    \k( r_{2n} +x_{n+1} ) \exp[ v_n + (2\e_0-\k/2)
    (r_{2n-1} -x_n)].
    \label{eqpoFncBis}
\end{eqnarray}
Now, integrate \eqref{eqpoFncBis} over $\EA(n)\cap \EB(n)$.
Since $\P[\EA(n)^c]$ and $\P[\EB(n)^c]$ are summable, this yields
since $\e_0<\k/4$,
\begin{equation}
    \label{eqlemmaprobaFncBis}
    \sum_{n=1}^{+\infty}
    \P\big[\fn(n)^c\big]
<
    \infty.
\end{equation}

To complete the proof of Lemma \ref{lemmaBorelCantelli}, let $0<\varepsilon<1/2$, and
define
\begin{eqnarray*}
    \mathcal{D}_n
 &:=& \left\{\exists t_n \in [(1-\e)r_{2n}, (1+\e)r_{2n}],
    \vphantom{\quad f[H_{X\circ \Theta_{H(Z_{2n-1})}} (t_n),
              (\lE\circ H)_{X\circ \Theta_{H(Z_{2n-1})}} (t_n)
              ]\in\D_n
    }\right.
    \\
 && \left.
    \vphantom{\exists t_n \in [(1-\e) r_{2n}, (1+\e) r_{2n}],
             }
    \qquad\qquad f\big[
                    H_{X\circ \Theta_{H(Z_{2n-1})}} (t_n),
                    (\lE\circ H)_{X \circ \Theta_{ H(Z_{2n-1})}} (t_n)
                  \big]
                  \in \D_n
    \right\},
\\
    \mathcal{E}_n
&:=&
    \Big\{
        \Big(1-5 r_{2n}^{-\d_0}/\k\Big)r_{2n}\leq
    F_{\wk\circ\Theta_{Z_{2n-1}}}
    (r_{2n})\leq\Big(1+5 r_{2n}^{-\d_0}/\k\Big)r_{2n}\Big\}.
\end{eqnarray*}
Let
$
    {\widetilde t_n}
:=
    F_{\wk\circ\Theta_{Z_{2n-1}}}
    (r_{2n})
$.
We have uniformly for large $n$,
\begin{equation}\label{InclusionDcapE}
    \mathcal{D}_n \cap \fn(n)
\supset
    \left\{
    f\big[H_{X \circ \Theta_{H(Z_{2n-1})}} \big({\widetilde t_n}\big), (\lE\circ H)_{X \circ
    \Theta_{H(Z_{2n-1})}} \big({\widetilde t_n}\big) \big] \in \D_n
    \right\}
    \cap
    \fn(n)
    \cap
    \mathcal{E}_n .
\end{equation}
Due to our assumption \eqref{e1p9}, $\sum_{n\geq 1} \P \{
f[H_{X\circ\Theta_{H(Z_{2n-1})}}({\widetilde t_n}),(\lE\circ
H)_{X\circ\Theta_{H(Z_{2n-1})}}({\widetilde t_n}) ]\in\D_n \}
=\infty$, since $X\circ\Theta_{H(Z_{2n-1})}$ is a diffusion
process in the $(-\k/2)$-drifted Brownian potential
$\wk\circ\Theta_{Z_{2n-1}}$,
which also gives $\P(\mathcal{E}_n)=\P(\EW(r_{2n}))$. In view of
\eqref{eqlemmaprobaFncBis}, \eqref{InclusionDcapE} and Lemma \ref{lemmaEncadrementFdeR},
this yields $\sum_{n\in\N}\P( \mathcal{D}_n\cap\fn(n))=+\infty$.

Define $x\wedge y:=\inf\{x,y\}$, $(x,y)\in\R^2$.
Since
$\e r_{2n}\leq r_{2n+1}/2$ for large $n$, the event $\mathcal{D}_n\cap \fn(n)$ is measurable
with respect to the $\sigma$-field generated by
$(\wk(x+Z_{2n-1})-\wk(Z_{2n-1}), -r_{2n-1}/2\leq x \leq Z_{2n}+r_{2n+1}/2-Z_{2n-1})$
and
$
\big(X\circ\Theta_{H(Z_{2n-1})}(t), \ 0\leq t \leq
    H_{X\circ\Theta_{H(Z_{2n-1})}}(-r_{2n-1}/2)
\wedge
    H_{X\circ\Theta_{H(Z_{2n-1})}}(Z_{2n}+r_{2n+1}/2-Z_{2n-1}\big)
$.
So, the events $\mathcal{D}_n\cap \fn(n)$, $n\geq 1$, are independent by the strong Markov Property,
because the intervals $\big[Z_{2n-1}-r_{2n-1}/2,Z_{2n}+r_{2n+1}/2\big)$, $n\geq 1$ are disjoint.
Hence,  Lemma \ref{lemmaBorelCantelli} follows by an application
of the \bc lemma. \hfill$\Box$





\mysection{Almost sure asymptotics  of  $\lm[H(r)]$ \label{SectEtudeLmdeHdeR}}


As a warm up, we first prove the following results, which are useful in Section \ref{SectJointStudy}.

\begin{theo}
 \label{lemmaLimsupL0}
Let $\k>0$. For any positive nondecreasing function
$a(\cdot)$, we have
$$
  \sum_{n=1}^{\infty}\frac{1}{n a(n)} \
  \left\{
  \begin{array}{l} <\infty
       \\
       =+\infty
  \end{array}
  \right.
  \Longleftrightarrow \limsup_{r\to\infty} \frac{\lm[H(r)]}{[r a(r)]^{1/\k}}=
  \left\{
  \begin{array}{l} 0
       \\
       +\infty
  \end{array}
  \right.
  \qquad\P\text{--a.s.}
 $$
\end{theo}

\medskip

\begin{theo}
 \label{lemmaLiminfL0}
 For $\k>0$,
$$
    \liminf\limits_{r\to+\infty}
    \frac{\lm[H(r)]}{ (r/ \log\log r)^{1/\k}}
=
    4\left(\frac{\k^2 }{2}\right)^{1/\k}\qquad
\P\text{--a.s.}
$$
\end{theo}


\subsection{Proof of Theorem \ref{lemmaLimsupL0}}
 \label{sectLimsupKsup1}

Let $r_n:=e^n$ and $Z_n:=\sum_{k=1}^n r_k$. Denote by $a(\cdot)$ be a
positive nondecreasing function. We begin with the upper bound
in Theorem \ref{lemmaLimsupL0}.

\noindent First, notice that  for $\Lb_\pm$ which is defined in
(\ref{e2p22}), and any positive $y$ and $r$, we have
\begin{equation}
    \P\left(\Lb_\pm(r)<\left(y r \right)^{1/\k}\right)
=
    \P\bigg[ \, \sup_{0\leq u\leq \tau_{\b}(\l t_\pm(r))}
    \b(u) <\frac{y r}{4^{\k}\k}\bigg]
=
    \exp\left(-\frac{\k^2 4^\k\p_\pm(r)}{2y}\right),
\label{eqProbaLiInf}
\end{equation}
by \eqref{borodin} and \eqref{e1p22}.
This together with Lemma \ref{lemmaApproxLxLi} gives,
for some $\alpha>0$, $\e>0$ and all large $r$,
\begin{equation}
    \P\left\{\lm[H(F(r))] >
    \left(r a(e^{-2}r)\right)^{1/\k}\right\}
\le
    1-
    \exp\bigg(-\frac{(1+\e)^\k \k^2
    4^\k\p_+(r)}{2a(e^{-2}r)} \bigg)
    +r^{-\alpha}
\le
    \frac{c_{7}}{a(e^{-2}r)}
    +r^{-\alpha},
\label{e1p24}
\end{equation}
since $1-e^{-x}\leq x$ for all $x\in\R$.
Assume $\sum_{n=1}^{+\infty} \frac{1}{n a(n)}<\infty$, which is
equivalent to $\sum_{n=1}^{+\infty} \frac{1}{a(r_n)}<\infty$. Then it follows from \eqref{e1p24} that
$
    \sum_{n=1}^{+\infty}\P\{\lm[H(F(r_n))]>[r_n
    a(r_{n-2})]^{1/\k}\}<\infty.
$

So by the Borel--Cantelli lemma, almost surely for all
large $n$, $\lm[H(F(r_n))]\leq [r_n a(r_{n-2})]^{1/\k}$. On the
other hand, $r_{n-1}\leq F(r_{n})$ almost surely for all large
$n$ (see \eqref{F}). As a consequence, almost surely for all
large $n$, $\lm[H(r_{n-1})]\leq [r_n a(r_{n-2})]^{1/\k}$.
Let $r\in[r_{n-2},r_{n-1}]$, for such large $n$. Then
$$
\lm[H(r)]\leq \lm[H(r_{n-1})] \leq [r_{n}a(r_{{n-2}})]^{1/\k}
\leq e^{2/\k} [r a(r)]^{1/\k}.
$$
Consequently,
\begin{equation}
    \label{eqMajorLimsupCprimeKsup1}
    \limsup_{r\to+\infty}\frac{\lm[H(r)]}{[r a(r)]^{1/\k}}\leq
    e^{2/\k} \qquad \P\text{--a.s.}
\end{equation}
Since $\sum_{n=1}^{+\infty} \frac{1}{n \e a(n)}$ is also finite, \eqref{eqMajorLimsupCprimeKsup1} holds
for $a(\cdot)$ replaced by $\e a(\cdot)$, $\e>0$. Letting $\e\to 0$
yields the ``zero'' part of Theorem \ref{lemmaLimsupL0}.


Now we turn to the proof of the "infinity" part. Assume
$\sum_{n=1}^{+\infty} \frac{1}{n a(n)}=+\infty$, that is,
$\sum_{n=1}^{+\infty} \frac{1}{a(r_n)} =+\infty$. Observe that we
may restrict ourselves to the case $a(x)\to+\infty$ when
$x\to+\infty$, since the result in this case yields the result
when $a$ is bounded.

By an argument similar to that leading to \eqref{e1p24}, we
have, for some $\alpha>0$ and all  large $r$,
$$
    \P\left\{\lm[H(F(r))] > \left(r a(e^2r)\right)^{1/\k}\right\}
\geq
    \frac{c_{8}}{a(e^2r)}
    -r^{-\alpha},
$$
which implies
$
\sum_{n=1}^{+\infty}\P\left\{(\lm\circ H\circ F )(r_{2n})>
[r_{2n}a(r_{2n+2})]^{1/\k}\right\}=+\infty.
$
Let $0<\e<1/2$ and recall that $Z_n=\sum_{k=1}^n r_k$; by Lemma \ref{lemmaBorelCantelli}, almost surely, there
exist infinitely many $n$ such that
$$
    \sup_{s\in[(1-\e)r_{2n},(1+\e)r_{2n}]}(\lE\circ H)_{X\circ{\Theta_{H(Z_{2n-1})}}}(s)
        >
    \left[r_{2n}a(r_{2n+2})\right]^{1/\k}.
$$
For such $n$, we have $\sup_{s \in [(1-\e)r_{2n}, (1+\e)r_{2n}]}\lm [H(Z_{2n-1}+s)]> [r_{2n}a(r_{2n+2})]^{1/\k}$.
Consequently,
$$
    \sup_{s\in[(1-\e)r_{2n},(1+\e)r_{2n}]}
    \frac{\lm(H(Z_{2n-1}+s))}{[(Z_{2n-1}+s) a(Z_{2n-1}+s)]^{1/\k}}
\geq
    c_{9},
$$
almost surely for infinitely many $n$. This gives
$$
    \limsup_{r\to+\infty} \frac{\lm[H(r)]}{[r a(r)]^{1/\k}}\geq c_{9}
\qquad
    \P\text{--a.s.}
$$
Replace $a(\cdot)$ by $a(\cdot)/\e$, and let $\e\to
0$. This yields the ``infinity'' part of Theorem
\ref{lemmaLimsupL0}.\hfill$\Box$


\subsection{Proof of Theorem \ref{lemmaLiminfL0}}
\label{SubSectPreuvelemmaliminfL0}
We fix $\e\in(0,1)$.
By Lemma \ref{lemmaApproxLxLi} and \eqref{eqProbaLiInf}, we get for some $\alpha>0$,
for every positive function $g$ and all large $r$,
\begin{equation}
    \label{eqProbaMajorL1}
    \P\left[ \lm[H(F(r))] < \left[r/g(r) \right]^{1/\k}\right]
\leq
    \exp \Big[ -\k^2 4^\k(1-\e)^\k\p_- (r)g(r)/2 \Big]
    + r^{-\alpha}.
\end{equation}
We choose $g(r):= \frac{2(1+\e)}{\k^24^\k(1-\e)^{\k+1}
\p_-(r)} \log\log r$. Let $s_n:=\exp(n^{1- \e})$.
It follows from \eqref{eqProbaMajorL1} that
$
    \sum_{n=1}^{\infty}
    \P\big\{\lm[H(F(s_n))]
            <
            \left[s_n/g(s_n) \right]^{1/\k}
    \big\}
<
    \infty.
$
Hence  by the \bc lemma, almost surely for all large $n$,
$$
\lm[H(F(s_n))]\geq [s_n/g(s_n)]^{1/\k}.
$$
On the other hand, by
\eqref{eqProbaEW} and the Borel-Cantelli lemma,
$s_n\geq F(s_{n-1})$ almost
surely for all large $n$, which implies that, for $r\in [s_n,
s_{n+1}]$,
$$
\lm[H(r)] \geq \lm[H(F(s_{n-1}))] \geq
[s_{n-1}/g(s_{n-1})]^{1/\k} \geq (1-\e) [r/g(r)]^{1/\k},
$$
since $s_{n-1}/s_{n+1}\to 1$ as $n\to+\infty$.
Consequently,
$$
    \liminf \limits_{r\to\infty}
    \frac{\lm[H(r)]}{(r/\log\log r)^{1/\k}}
\geq
    4\left( \frac{\k^2 }{2} \right)^{1/\k} \qquad
\P\text{--a.s.}
$$


Now we prove the inequality ``$\leq$". Let $\e\in(0,1/2)$, $r_n:=\exp(n^{1+\e})$,
$Z_n:=\sum_{k=1}^n r_k$, $n\geq 1$,
and ${\widetilde g}(r) :=
\frac{2(1-\e)}{\k^2 4^\k(1+\e)^{\k+1} \p_+(r)}\log\log r$. By Lemma
\ref{lemmaApproxLxLi} and \eqref{eqProbaLiInf}, for some $\alpha>0$ and all large
$r$,
$$
    \P \left[ \lm[H(F(r))] < \left[ r/{\widetilde g}(r)\right]^{1/\k} \right]
\geq
    \exp\Big[ - \k^2 4^\k(1+\e)^\k\p_+(r) {\widetilde g}(r)/2 \Big]
    -r^{-\alpha}.
$$
Therefore,
$$
    \sum_{n\geq 1}
        \P \left[ \lm[H(F(r_{2n}))]
                <
                \big[ r_{2n}/{\widetilde g} (r_{2n}) \big]^{1/\k}
           \right]
=
    +\infty.
$$
It follows from Lemma \ref{lemmaBorelCantelli}
that, almost surely, there are infinitely many $n$ such that
\begin{equation}\label{eqPourCiterAutreChose}
    \inf_{s\in[(1-\e)r_{2n},(1+\e)r_{2n}]}
    (\lE\circ H)_{X\circ\Theta_{H(Z_{2n-1})}}(s)
<
    \big[r_{2n}/{\widetilde g} (r_{2n})\big]^{1/\k}.
\end{equation}
On the other hand,
an application of Theorem \ref{lemmaLimsupL0} with $a(x)\sim_{x\to +\infty} (\log x)^2$
gives that almost surely for large $n$,
$
    \lm[H(Z_{2n-1})]
\leq
    [Z_{2n-1}\log^2 Z_{2n-1}]^{1/\k}
\leq
    \e
    \big[r_{2n}/{\widetilde g} (r_{2n})\big]^{1/\k},
$
 since $Z_p\leq pr_p\leq p\exp(-p^{\e})r_{p+1}$ for $p$
large enough. Therefore,
$
    \inf_{s\in[(1-\e)r_{2n},(1+\e)r_{2n}]}\lm[H(Z_{2n-1}+s)]
\leq
    (1+\e) \big[r_{2n}/{\widetilde g}(r_{2n})\big]^{1/\k}
$
almost surely, for infinitely many $n$,
where we used
$
    L_X^*[H(r+s)]
\leq
    L_X^*[H(r)]
    +
    (L^*\circ H)_{X\circ\theta_{H(r)}}(s)
$,
$r\geq 0$, $s\geq 0$.
Hence, for such $n$,
$$
    \inf_{s\in[(1-\e)r_{2n},(1+\e)r_{2n}]}
    \frac{\lm[H(Z_{2n-1}+s)]}{[ (Z_{2n-1}+s)/ \log\log (Z_{2n-1}+s)]^{1/\k}}
\leq
    (1+c_{10}\e)\left(\frac{\k^2 4^\k\p_+(r_{2n})}{2} \right)^{\! \! 1/\k}.
$$
This yields
$$
 \liminf_{r\to+\infty}\frac{\lm(H(r))}{(r/\log\log
 r)^{1/\k}}\leq 4 \left(\frac{\k^2}{2}\right)^{1/\k}\qquad
 \P\text{--a.s.,}
$$
proving Theorem \ref{lemmaLiminfL0}. \hfill$\Box$



\mysection{Proof of Theorems \ref{ThUpperLevyKinf1} and
\ref{thLowerLevyKinf1} and Corollary \ref{CorrolaryLILofX}}\label{SectEtudeHdeR}



Recall $\Ib_\pm$ from \eqref{eqI''CasKinf1} and
$c_4(\k)$
from \eqref{DefPsi}. By Fact \ref{FactBianeYor},
\begin{eqnarray}
    \Ib_\pm(r)
& \el &
    t_\pm(r)^{1/\k}\big\{ c_4(\k) \, S_\k^{ca} \pm c_{11} \,
    t_\pm(r)^{1-1/\k}\big\},
    \qquad 0<\k<1,
    \label{eqDefItierceKinf1}
\\
\label{eqDefItierceKeg1}
    \Ib_\pm(r)
&\el&
    4t_\pm(r)
    [8c_3 +(\pi/2)C_8^{ca} + 8\log t_\pm(r) ]
    \qquad \k=1,
\end{eqnarray}
where
 $c_{11}>0$ and $c_3>0$ are unimportant constants.
We have now all the ingredients to prove Theorems
\ref{ThUpperLevyKinf1} and \ref{thLowerLevyKinf1}.


\subsection{Proof of Theorem \ref{ThUpperLevyKinf1}}
 \label{SubSecUpperLevyH}


\subsubsection{Case $0<\k<1$}

We assume $0<\k<1$. Let $a(\cdot)$ be
a positive nondecreasing function.
Without loss of generality, we suppose that $a(r) \to \infty$
(as $r\to \infty$).

It is known  (see e.g. Samorodnitsky and Taqqu \cite{Papier_3_SamoTaqqu}, (1.2.8) p.~16) that
$$
    \P\big(S_\k^{ca}>x\big)
\underset{x\to +\infty}{\sim}
    c_{12} x^{-\k},
$$
where
$
f(x)
\underset{x\to +\infty}{\sim}
g(x)
$ means
$\lim\limits_{x\to+\infty} f(x)/g(x)=1$, and $c_{12}>0$ is a
constant depending on $\k$.

Recall $t_\pm(\cdot)$ from \eqref{e1p22}. By Lemma
\ref{lemmaApproxLxLi} and \eqref{eqDefItierceKinf1}, for some $\alpha>0$,
we have for large $r$,
\begin{eqnarray}
    \P\big[H(F(r))>\big(a(e^{-2}r) t_+(r)\big)^{1/\k}\big]
& \le &
    \frac{c_{13}}{a(e^{-2}r )} + r^{-\alpha}.
    \label{e1p42}
\end{eqnarray}
As in Section \ref{sectLimsupKsup1}, we
define $r_n:=e^n$ and $Z_n:=\sum_{k=1}^n r_k$.
Assume $\sum_{n\geq 1} \frac{1}{a(r_n)}<\infty$,
which is equivalent to $\sum_{n\geq 1}\frac{1}{n a(n)}<\infty$. By the
Borel--Cantelli lemma, almost surely for $n$ large enough,
\begin{equation}\label{eqHrn}
    H[F(r_n)]
\leq
    [a(r_{n-2})t_+(r_n)]^{1/\k}.
\end{equation}
On the other hand, by Lemma
\ref{lemmaEncadrementFdeR}, almost surely for all large $n$,
we have $r_{n+1}\le F(r_{n+2})$, which together with \eqref{eqHrn} implies that for
$r\in [r_n, r_{n+1}]$,
$$
H(r)
\leq H[F(r_{n+2})]
\leq [\p_+(r_{n+2})\k r_{n+2}a(r_{n})/\l]^{1/\k}
\leq c_{14} [ ra(r)]^{1/\k}.
$$
Therefore,
$\limsup_{r\to+\infty}\frac{H(r)}{[r a(r)]^{1/\k}}\leq c_{14}$
$\P\text{--a.s.},$
implying the ``zero'' part of Theorem
\ref{ThUpperLevyKinf1}, since we can replace $a(\cdot)$ by any
constant multiple of $a(\cdot)$.


To prove the ``infinity" part, we assume
$\sum_{n\geq 1} \frac{1}{ n a(n)} = +\infty$, and observe that, by an argument
similar to that leading to \eqref{e1p42}, we have, for some $\alpha>0$ and all $r$ large enough,
\begin{equation}
    \P\big[H(F(r))>\big(a(e^2r)t_-(r)\big)^{1/\k}\big]
\geq
    \frac{c_{15} }{a(e^2r)}-r^{-\alpha}.
\label{Maine1}
\end{equation}
It follows from Lemma \ref{lemmaBorelCantelli} that
$
    \sup_{s\in[(1-\e)r_{2n},(1+\e)r_{2n}]}
    H_{X\circ{\Theta_{H(Z_{2n-1})}}}(s)
>
    [a(r_{2n+2})t_-(r_{2n})]^{1/\k} ,
$
almost surely for infinitely many $n$. Since
$
    H(Z_{2n-1}+s)
\geq
    H_{X\circ{\Theta_{H(Z_{2n-1})}}}(s)
$
for all $s>0$, this implies, for these $n$,
\begin{equation}
    \sup_{s\in[(1-\e)r_{2n},(1+\e)r_{2n}]}
    H(Z_{2n-1}+s)/[a(Z_{2n-1}+s)(Z_{2n-1}+s)]^{1/\k}
\geq
    c_{16}.
\label{Maine2}
\end{equation}
This gives
$
\limsup_{r\to+\infty}\frac{H(r)}{[r a(r)]^{1/\k}}\geq c_{16}
$
$\P$-a.s.,
proving the ``infinity'' part in Theorem
\ref{ThUpperLevyKinf1}, in the case $0< \k<1$ by replacing $a(.)$ by any constant multiple of $a(.)$.\hfill$\Box$


\subsubsection{Case $\k=1$}

Let $r_n:=e^n$ and $Z_n:=\sum_{k=1}^n r_k$. We recall that there
exists a constant $c_{17}:=\frac{16}{\pi}$ such that
$\P(C_8^{ca}>x)\underset{x\to +\infty}{\sim} \frac{c_{17}}{x}$
(see e.g. Samorodnitsky et al. \cite{Papier_3_SamoTaqqu}, prop. 1.2.15
p.~16). Hence,
by Lemma \ref{lemmaApproxLxLi} and \eqref{eqDefItierceKeg1}, for some $\alpha>0$ and all large $r$,
\begin{equation}
 \P \left\{ H(F(r))>4t_+(r)(1+\e) \big[8c_3 +a(e^{-2}r) +
    8\log t_+ (r) \big] \right\}
\le
    c_{18}/a(e^{-2}r)
    + r^{-\alpha}.
    \label{littlehope}
\end{equation}
Assume $\sum_{n\geq 1} \frac{1}{na(n)} <\infty$. Then by
the \bc lemma, almost surely, for all large $n$,
$$
    H[F(r_n)]
\leq
    4(1+\e)t_+(r_n)
    [8c_3+a(r_{n-2})+ 8 \log (\p_+(r_n) \k r_n/8)].
$$
Under the additional assumption $\limsup_{r\to +\infty}(\log
r)/a(r)<\infty$, we have, almost surely, for all large $n$ and
$r\in [r_n, r_{n+1}]$ (thus $r \le F(r_{n+2})$ by Lemma
\ref{lemmaEncadrementFdeR}),
$$
    H(r)
\leq
    H[F(r_{n+2})]
\le
    c_{19} r_{n+2} [ a(r_{n})+ \log r_{n+2}]
\le
    c_{20}r a(r) .
$$
As in the case $0<\k<1$, this yields the ``zero'' part of Theorem
\ref{ThUpperLevyKinf1} in the case $\k=1$.


For the ``infinity" part, we assume
$\sum_{n\geq 1} \frac{1}{na(n)} =+\infty$. As in \eqref{littlehope}, we have, for some $\alpha>0$
$$
    \P \left\{ H(F(r))>4t_-(r)(1-\e) a(e^2r) \right\}
\geq
    c_{21}/a(e^2r) -r^{-\alpha},
$$
for large $r$.
%
As in the displays between \eqref{Maine1} and \eqref{Maine2}, this yields the "infinity part" of Theorem \ref{ThUpperLevyKinf1} in the case $\k=1$.
\hfill$\Box$


\subsection{Proof of Theorem \ref{thLowerLevyKinf1}}
\label{SubSecLowerLevyH}


\subsubsection{Case $0<\k<1$}
\label{SubSubSectPreuveThlowerLevyKinf1}

We have
$
    \E\big(e^{-t S_\k^{ca}}\big)
=
    \exp[-t^\k/\cos(\pi\k/2)]
$, $t\geq 0$,
e.g. by Samorodnitsky et al. (\cite{Papier_3_SamoTaqqu}, Proposition 1.2.12,
in the notation of \cite{Papier_3_SamoTaqqu}, $S_\k^{ca}$ is distributed as
$S_\k(1,1,0)$).
So by Bingham et al. (\cite{Bingham_Teugels} Example p. 349),
\begin{equation}\label{bertoin}
    \log\P\big( S_\k^{ca} <x\big)
\underset{x\to 0,\ x>0}{\sim} -
    c_{22} x^{-\k/(1-\k)},
\end{equation}
where $c_{22}:=(1-\k)\k^{\k/(1-\k)}[\cos(\pi\k/2)]^{-1/(1-\k)}$.
By Lemma \ref{lemmaApproxLxLi}, \eqref{eqDefItierceKinf1} and \eqref{bertoin}, for any
(strictly) positive
function $f$ such that $\lim_{x\to+\infty}f(x)=0$ and
$\e>0$ small enough, we have for large $r$,
\begin{equation}\label{eqH+inf1}
    \P\big[H(F(r)) < t_-(r)^{1/\k}f(r)\big]
\leq
    \exp \bigg[ -(c_{22}- \e)
        \left(\frac{ (1-\e) c_4(\k) }{f(r)+ (1-\e)c_{11} t_-(r)^{1-1/\k} }
    \right)^{\k/(1-\k)}
    \bigg]
    +r^{-\alpha}.
\end{equation}
We define for $\e>0$ and $r>1$,
$$
    f_\e^\pm(r)
:=
    (1\pm\e)
    c_4(\k)
    \bigg( \frac{(1\pm\e)(c_{22}\pm\e)}{(1\mp\e)\log\log r}\bigg)^{(1-\k)/\k}
    \pm c_{11}(1\pm\e) t_\pm(r)^{1-1/\k}.
$$
So, \eqref{eqH+inf1} gives
\begin{equation*}
    \P\big[H(F(r)) < t_-(r)^{1/\k}f_\e^-(r)\big]
\le
    (\log r)^{-(1+\e)/(1-\e)}
    +r^{-\alpha}.
\end{equation*}
With $s_n:=\exp(n^{1-\e})$,
this gives $\sum_n^{+\infty} \P\big[ H(F(s_n)) < t_-(s_n)^{1/\k} f_\e^-(s_n)\big] < \infty$,
which, by
the \bc lemma, implies that, almost surely, for all large $n$,
$
    H[F(s_n)]
\geq
    t_-(s_n)^{1/\k}f_\e^-(s_n)
$.

 Recall from Lemma \ref{lemmaEncadrementFdeR} that,
almost surely, for all large $n$, we have $F(s_n) \le
(1+\varepsilon) s_n$. Let $r$ be large. There exists $n$ (large)
such that $(1+\varepsilon) s_n \le r\le (1+2\varepsilon)s_n$.
Then if $r$ is large,
$$
    H(r)
\ge
    H[F(s_n)]
\ge
    t_-(s_n)^{1/\k} f_\e^-(s_n)
\ge
    t_-^{1/\k}
    \bigg(\frac{r}{1+2\e}\bigg)
    f_\e^-
    \left(\frac{r}{1+\e}\right) .
$$
Plugging the value of $t_-( \frac{r}{1+2\e})$
(defined in \eqref{e1p22}), this yields inequality "$\geq$" of
\eqref{Maine3} with
\begin{equation}
    c_2(\k)
:=
    8\psi(\k) c_{22}^{(1-\k)/\k}
=
    8 [\pi \k]^{1/\k}  (1-\k)^{\frac{1-\k}{\k}}\k /\big[2 \Gamma^2(\k)\sin(\pi\k)\big]^{1/\k}
   \label{definitionc2}
\end{equation}
where $c_{22}=c_{22}(\k)$ is defined after \eqref{bertoin} and $\psi$ and $c_4(\k)$ in
\eqref{DefPsi}.


To prove the upper bound, let $r_n:=\exp(n^{1+\e})$
and
$Z_n:=\sum_{k=1}^n r_k$.
By means of an argument similar to that leading to
\eqref{eqH+inf1},
we have
$
    \sum_{n\ge 1}
    \P\big[ H(F(r_{2n})) <t_+(r_{2n})^{1/\k} f_\e^+(r_{2n}) \big]
=
    +\infty
$.
So by Lemma
\ref{lemmaBorelCantelli}, for $0<\e<1/2$, there exist almost surely infinitely
many $n$ such that
$$
    \inf_{u\in[(1-\e)r_{2n},(1+\e)r_{2n}]}
    H_{X\circ{\Theta_{H(Z_{2n-1})}}}(u)
<
    [t_+(r_{2n})]^{1/\k}f_\e^+(r_{2n}).
$$
In addition, by Theorem \ref{ThUpperLevyKinf1},
$
    H(Z_{2n-1})
<
    \left[  Z_{2n-1}\log^2 Z_{2n-1}\right]^{1/\k}
\leq
    \e[t_+(r_{2n})]^{1/\k}f_\e^+(r_{2n})
$
almost surely for all large $n$,
since $\sum_{n\geq 1} 1/(n \log^2 n)<\infty$
and $Z_p\le p \exp(-p^\e)r_{p+1}$ for all large $p$ as before.
This yields almost surely for large $n$,
$$
    \inf_{v\in[Z_{2n-1}+(1-\e)r_{2n}, \, Z_{2n-1}+(1+\e)r_{2n}]} H(v)
<
    (1+\e)[t_+(r_{2n})]^{1/\k}f_\e^+(r_{2n}).
$$
Consequently,
$$
    \liminf_{r\to+\infty}\frac{H(r)}{r^{1/\k}(\log\log r )^{(\k-1)/\k}}
\leq
    8\psi(\k) c_{22}^{(1-\k)/\k}
=
    c_2(\k)
\qquad
    \P\text{--a.s.}
$$
This gives inequality "$\leq$" of \eqref{Maine3} and thus yields
Theorem \ref{thLowerLevyKinf1} in the case $0<\k<1$.




\subsubsection{Case $\k=1$}
\label{SubSubSectPreuveThlowerLevyKinf1K1}

Assume $\k=1$ (thus $\l=8$). By Samorodnitsky et al.
(\cite{Papier_3_SamoTaqqu}, Proposition 1.2.12),
$\E[\exp(-C_8^{ca})]=1$ (in the notation of
\cite{Papier_3_SamoTaqqu}, $C_8^{ca}$ is distributed as
$S_1(8,1,0)$). Hence,
\begin{equation}
    \label{ProbaSamoTaqqu}
    \P[ C_8^{ca} \le -\varepsilon\log r ]
\le
    r^{-\varepsilon}\,
    \E[ \exp( -C_8^{ca}) ] = r^{-\varepsilon},
\qquad
    r>0,
\end{equation}
for $\e>0$.
By Lemma \ref{lemmaApproxLxLi} and
\eqref{eqDefItierceKeg1}, we have if $\e>0$ is small enough, for all large $r$,
\begin{equation*}
    \P \big\{H[F(r)]\leq 32 t_-(r) (1-2\varepsilon) [c_3
    +
    \log t_-(r)] \big\}
\le
    \P\left( C_8^{ca}\le - \varepsilon \log r\right)
    +\P[\ELb(r)^c]
\le
    2r^{-\e}.
\end{equation*}
Let $s_n := \exp(n^{1-\varepsilon})$. Thus, by the \bc
lemma, almost surely, for all large $n$,
$$
    H[F(s_n)]
>
    32 t_-(s_n) (1-2\varepsilon) [c_3 + \log t_-(s_n)]
\ge
    4(1-3\e)s_n\log s_n.
$$
In view of the last part of Lemma
\ref{lemmaEncadrementFdeR}, this yields inequality "$\geq$" in
\eqref{Maine4} similarly as before \eqref{definitionc2}.
The inequality ``$\le$'', on the other hand,  follows immediately from Theorem
\ref{theoKawazuTanaka} (that $H(r)/(r\log r)\to 4$ in
probability).
Theorem \ref{thLowerLevyKinf1} is proved.
\hfill$\Box$

\subsection{Proof of Corollary \ref{CorrolaryLILofX}}
First, we need the following lemma, which says that $X$ does not go back too far on the left,
and so $X(t)$ is very close from $\sup_{0\leq s \leq t} X(s)$:

\begin{lemma}\label{LemmaNotTooFarOnTheLeft}
For every $\k>0$, there exists a constant $c_{23}(\k)$ such that $\P$ a.s. for large $t$,
\begin{equation}\label{eqLemmaNotTooFarOnTheLeft}
    0
\leq
    \sup_{0\leq s \leq t} X(s)
    -
    X(t)
\leq
    c_{23}(\k) \log t.
\end{equation}
\end{lemma}

Notice that this is not true in the recurrent case $\k=0$.
An heuristic explanation for $0\leq \k < 1$ would be that the valleys of height approximatively $\log t$
have a length of order $(\log t)^2$ in the case $\k=0$, whereas they have a height of order at most
$\log t$ in the case $0<\k<1$, see e.g. Andreoletti et al. (\cite{Andreoletti_Devulder}, Lem. 2.7).

\noindent{\bf Proof of Lemma \ref{LemmaNotTooFarOnTheLeft}:}
Let $\k>0$.
By Kawazu et al. (\cite{KawazuTanakaMaximum}, Theorem p. 79 applied with $c=\k/2$ to our $-X$),
there exists a constant $c_{24}(\k)>0$ such that
$
    \P\big[\inf_{u\geq 0} X(u) < - c_{24}(\k) \log n \big]
\leq
    1/n^2
$
for large $n$.
Since $\inf_{u\geq 0} X(H(n)+u)-n$ has the same law under $\P$ as
$\inf_{u\geq 0} X(u)$ due to the strong Markov property as explained before Lemma \ref{lemmaBorelCantelli},
this gives
$
    \sum_n \P\big[
                \inf_{u\geq 0} X(H(n)+u)-n
                < -c_{24}(\k) \log n
             \big]
<
    \infty
$.
So by the Borel-Cantelli lemma, almost surely for large $n$,
\begin{equation}\label{InegXetSonSup}
    \inf_{u\geq 0} X(H(n)+u)-n
\geq
    -c_{24}(\k) \log n.
\end{equation}
For $t>0$, there exists $n\in\N$ such that $H(n)\leq t <H(n+1)$.
We have by \eqref{InegXetSonSup}, almost surely if $t$ is large,
$$
    \sup_{0\leq s \leq t} X(s)
    -
    X(t)
\leq
    \sup_{0\leq s \leq H(n+1)} X(s)
    -
    X(t)
=
    n+1-X[H(n)+(t-H(n))]
\leq
    1+c_{24}(\k) \log n.
$$
Moreover, we have
$
    \log v
\leq
    2\log H(v)
$
$\P$ a.s. for large $v$, by Theorem \ref{theoKawazuTanaka} if $\k>1$
and by Theorem \ref{thLowerLevyKinf1} if $0<\k\leq 1$.
Hence almost surely for large $t$, with the same notation as before,
$$
    \sup_{0\leq s \leq t} X(s)
    -
    X(t)
\leq
    1+c_{24}(\k) \log n
\leq
    1+2 c_{24}(\k) \log H(n)
\leq
    1+2 c_{24}(\k) \log t.
$$
This proves the second inequality of \eqref{eqLemmaNotTooFarOnTheLeft}.
The first one is clear.
\hfill $\Box$

\noindent{\bf Proof of Corollary \ref{CorrolaryLILofX}:}
By Lemma \ref{LemmaNotTooFarOnTheLeft},
$$
    \limsup_{t\to\infty}\big[X(t)/(t/\log t)\big]
=
    \limsup_{t\to\infty}\Big[\Big(\sup_{0\leq s \leq t}X(s)\Big)/(t/\log t)\Big].
$$
So, \eqref{eqLimsupXcaskegal1} is equivalent to \eqref{eqLimsupXcaskegal1} with
$X(t)$ replaced by $\sup_{0\leq s \leq t}X(s)$.
The same remark also applies to \eqref{eqLimsupXcaskentre0et1} and \eqref{eqLiminfXcaskentre0et1}.

Now, we have
$
    \sup_{0\leq s \leq y}X(s)
\geq
    r
\Longleftrightarrow
    H(r)
\leq
    y
$,
$r>0$, $y>0$.
Consequently \eqref{eqLimsupXcaskentre0et1}, \eqref{eqLimsupXcaskegal1} and \eqref{eqLiminfXcaskentre0et1}
with $X(t)$ replaced by $\sup_{0\leq s \leq t} X(s)$
follow respectively from \eqref{Maine3}, \eqref{Maine4} and Theorem \ref{ThUpperLevyKinf1} applied to
$a(r)=(\log r)\dots (\log_{k-1} r)(\log_k r)^\alpha$.
Indeed for \eqref{eqLiminfXcaskentre0et1} when $\k=1$,
cases $k=1$, $\alpha\leq 1$ and $k=2$, $\alpha\leq 0$ follow from the case $k=3$, $\alpha=1$.
This proves Corollary \ref{CorrolaryLILofX}.
\hfill $\Box$



\mysection{Proof of Theorems
\ref{theoLimsupDifferentRWRE} to
\ref{theoLawLm}
}\label{SectJointStudy}

\smallskip
\noindent {\bf Proof of Theorem \ref{theoLiminfKsup1}: case
$\k>1$.} Follows from Theorems \ref{lemmaLiminfL0} and
\ref{theoKawazuTanaka}.\hfill$\Box$

\medskip
\noindent {\bf Proof of Theorem \ref{theoLimsupKsup1}.} Follows
from Theorems \ref{lemmaLimsupL0} and
\ref{theoKawazuTanaka}.\hfill$\Box$

\medskip
\noindent {\bf Proof of Theorem \ref{theoLawLm}.}
We first notice that for every $\k>0$, thanks to Lemma \ref{lemmaApproxLxLi} {\bf (i)},
\begin{equation}\label{eqConvergenceLawLHFr}
    \lm[H(F(r))]/ r^{1/\k}
\overset{\L}{\longrightarrow}
    4 \big[\k^2 /\l\big]^{1/\k} \big(\sup\nolimits_{0\leq u\leq \tau_\beta(\l)}\beta(u)\big)^{1/\k},
\end{equation}
where $\overset{\L}{\longrightarrow}$ denotes convergence in law under $\P$ as $r\to+\infty$.

We now assume $\k>1$. In this case,
$H(F(r))/r \to_{r\to+\infty} 4/(\k-1)$ $\P$--a.s.
by Lemma \ref{lemmaEncadrementFdeR} eq. \eqref{F} and Theorem \ref{theoKawazuTanaka} eq. \eqref{KawazuTanakakappaSuperieur1}.
This, combined with \eqref{eqConvergenceLawLHFr} leads to the convergence in law under $\P$
of $\lm(t)/t^{1/\k}$ to
$
    4[\k^2(\k-1)/(4\l)]^{1/\k}
    \big(\sup\nolimits_{0\leq u\leq \tau_\beta(\l)}\beta(u)\big)^{1/\k}
$.
Since $\sup\nolimits_{0\leq u\leq \tau_\beta(\l)}\beta(u)$ has by \eqref{borodin} the same law as $\l/(2\mathcal{E})$,
where $\mathcal{E}$ is an exponential variable with mean $1$, this proves Theorem \ref{theoLawLm} when $\k>1$.

We finally assume $\k=1 $. In this case,
$H(F[t/(4\log t)])/t\to_{t\to+\infty} 1$ in probability under $\P$
by Lemma \ref{lemmaEncadrementFdeR} and Theorem \ref{theoKawazuTanaka} eq. \eqref{KawazuTanakakappa1}.
This, combined with \eqref{eqConvergenceLawLHFr} leads to the convergence in law of
$\lm(t)/(t/\log t)$ to
$\l^{-1} \sup\nolimits_{0\leq u\leq \tau_\beta(\l)}\beta(u)$,
which proves Theorem \ref{theoLawLm} when $\k=1$.
\hfill$\Box$

\medskip
We now assume $0<\k\leq 1$, and need to prove Theorems
\ref{theoLimsupDifferentRWRE}, \ref{theoLiminfKsup1} and \ref{theoLiminfKinf1}. Unfortunately,
it follows immediately from Theorems \ref{ThUpperLevyKinf1} and \ref{thLowerLevyKinf1}
that there is no almost
sure convergence result for $H(r)$ in this case due to strong
fluctuations; hence a joint study of $\lm[H(r)]$ and
$H(r)$ is useful. In Section \ref{SubSecAlemma}, we prove a
lemma which will be needed later on. Section \ref{SubSec62} is
devoted to the proof of Theorems \ref{theoLimsupDifferentRWRE},
\ref{theoLiminfKsup1} and \ref{theoLiminfKinf1}   in the
case $0<\k<1$, whereas Section \ref{SubSec63} to the proof of
Theorems \ref{theoLiminfKsup1} and \ref{theoLiminfKinf1} in the
case $\k=1$.



\subsection{A lemma}\label{SubSecAlemma}

In this section we assume $0<\k\leq 1$. Let $\delta_1>0$
and recall the definitions of $t_\pm (r)$ from \eqref{e1p22}
and $\Lb_\pm(r)$ from \eqref{e2p22}.

\begin{lemma}
 \label{lemmaEncadrementQ}
Define $\EU(r):=\big\{ \Lb_- (r)=\Lb_+(r)\big\}$.
For all $\delta_2 \in (0, \delta_1)$ and all large $r$, we have
$$
    \P\big[\EU(r)^c\big]
\leq
    r^{-\d_2}.
$$
\end{lemma}

\medskip

\noindent{\bf Proof.}
Let $\delta_2 \in (0, \delta_1)$.
Observe that
\begin{equation}
    1
\le
    \bigg(\frac{\Lb_+(r)}{\Lb_-(r)}\bigg)^\k
\leq
    \max \bigg( 1,
                \frac{ \sup_{0\leq u\leq \tau_{\bti}
                \{[\p_+(r)-\p_-(r)] \k r\} } \bti(u)}
                {\sup_{0\leq u\leq \tau_{\beta} (\p_-(r)\k r)} \beta(u)}
         \bigg),
    \label{max}
\end{equation}
where $\bti(u):= \beta[u+ \tau_{\beta} (\p_-(r)\k r)]$, $u\ge 0$, is a Brownian motion independent of the random
variable $\sup_{0\le u\le \tau_{\beta} (\p_-(r)\k r)} \beta(u)$.
By \eqref{borodin} and the usual inequality $1-e^{-x}\le x$
(for $x\ge 0$),
\begin{eqnarray*}
    \P\bigg(\sup_{0\leq u\leq \tau_{\bti}\{[\p_+(r)-\p_-(r)]\k r\}} \bti(u)
                >
                [\p_+(r)-\p_-(r)]\k r^{1+\d_2}\bigg)
&\le&
    \frac{1}{2r^{\d_2}{}},
\\
    \P\bigg(\sup_{0\leq u\leq \tau_{\beta}(\p_-(r)\k r) } \beta(u)
                <
                \frac{\p_-(r)\k r}{4\d_2\log r}
    \bigg)
&=&
    \frac{1}{r^{2\d_2}}
\leq
    \frac{1}{2r^{\d_2}{}},
\end{eqnarray*}
for large $r$.
By definition, $\psi_\pm(r)= 1\pm c_{5}r^{-\d_1}$
(see (\ref{e1p22})). Therefore, we have for large $r$, with
probability greater than $1-r^{-\d_2}$,
$$
    \frac{ \sup_{0\leq u\leq \tau_{\bti} \{[\p_+(r)-\p_-(r)] \k r\}} \bti(u)}
    {\sup_{0\leq u\leq \tau_{\beta} (\p_-(r)\k r)}
    \beta(u)}
\le
    \frac{[\p_+(r)-\p_-(r)] \k r^{1+\d_2}}{\p_-(r)\k r/(4\d_2\log r)}
=
    \frac{8c_{5} \d_2 r^{-(\d_1-\d_2)} \log r }
    { 1-c_{5}r^{-\d_1}}
<
    1.
$$
This, combined with (\ref{max}), yields the
lemma.\hfill$\Box$


\subsection{Case $0<\k<1$}\label{SubSec62}

This section is devoted to the proof of Theorems
\ref{theoLimsupDifferentRWRE}, \ref{theoLiminfKsup1} and \ref{theoLiminfKinf1}
 in the case $0<\k<1$.


For any Brownian motion $(\beta(u),\ u\geq 0)$, let
$$
    N_\b
:=
    \frac{\int_0^{+\infty}x^{1/\k-2}L_{\b}(\tau_{\b}(\l),x)\text{d}x}
    {\big[ \sup_{0\leq u\leq \tau_{\b}(\l)}\b(u)\big]^{1/\k}} .
$$
So, in the notation of (\ref{e3p18}), (\ref{e1p22})
and (\ref{e2p22}),
$
    N_{\beta_{t_\pm(r)}}
=
    4[\k t_\pm(r) ]^{1/\k}K_{\b_{t_\pm(r)}}(\k)/\Lb_\pm(r)
$,
$r>0$.

On $\EL(r)\cap\ELb(r)\cap\EU(r)$ (the events $\EL(r)$ and
$\ELb(r)$ are defined in Lemma \ref{lemmaApproxLxLi}, whereas
$\EU(r)$ in Lemma \ref{lemmaEncadrementQ}), we have, for some
constant $c_{25}$, $\e>0$ small enough  and all large $r$,
\begin{eqnarray}
    \frac{H(F(r))}{\lm[H(F(r))]}
& \ge &
    \frac{4(1-\e)\k^{1/\k-2} t_-(r)^{1/\k}  \{
    K_{\b_{t_-(r)}} (\k) -c_{6} t_-(r)^{1-1/\k} \} }{(1+\e)
    \Lb_-(r) }
\nonumber\\
&\ge&
    (1-3\e) \k^{-2} N_{\beta_{t_-(r)}}
    -
    c_{25} t_-(r)/\Lb_-(r).
\label{eqCompQuotientHL>}
\end{eqnarray}
Similarly, on $\EL(r)\cap\ELb(r)\cap\EU(r)$, for some
constant $c_{26}$ and all large $r$,
\begin{equation}\label{eqCompQuotientHL<}
    \frac{H(F(r))}{\lm[H(F(r))]}
\le
    (1+3\e) \k^{-2} N_{\beta_{t_+(r)}}
    +
    c_{26} \frac{t_+(r)}{\Lb_+(r)}.
\end{equation}
Define
$
    \EV(r)
:=
    \big\{ c_{25} t_-(r)/\Lb_-(r) \leq \e, \;
    c_{26} t_+(r)/\Lb_+(r) \le \e \big\}
$.
By \eqref{eqProbaLiInf}, $\P[\EV(r)^c]\leq 1/r^2$ for large $r$.
Thus $\P[ \EL(r)\cap\ELb(r)\cap\EU(r) \cap \EV(r) ] \ge 1-
r^{-\alpha_1}$ for some $\alpha_1>0$ and all large $r$
by Lemmas \ref{lemmaApproxLxLi} and \ref{lemmaEncadrementQ}. In view
of (\ref{eqCompQuotientHL>}) and (\ref{eqCompQuotientHL<}), we
have, for some $\alpha_1>0$ and all large $r$,
\begin{equation}
\label{eqCompQuotientHL2}
    \P \bigg(
        (1-3\e) \k^{-2} N_{\beta_{t_-(r)}}  -\e \leq
        \frac{H(F(r))}{\lm[H(F(r))]}
    \leq
        (1+3\e) \k^{-2} N_{\beta_{t_+(r)}}  +\e
    \bigg)
\geq
    1-\frac{1}{r^{\alpha_1}}.
\end{equation}

We now proceed to the study of the law of $N_{\b}$. By the second
Ray--Knight theorem (Fact \ref{FactRN2}), there exists a
$0$--dimensional Bessel process $(U(x),\ x\geq 0)$, starting from $\sqrt{\l}$,
such that
\begin{eqnarray}
\label{e1p32}
   (L_\b(\tau_\b(\l),x),x\geq 0)
& = &
    \big(U^2(x),x\geq 0\big),
\\
\label{e3p31}
    \sup_{0\leq u \leq \tau_\b(\l)}\b(u)
& = &
    \inf\{x\geq 0,\ U(x)=0\}=:\z_U,
\\[-4mm]
\label{e4p31}
    N_\b & = & \z_U^{-1/\k}\int_0^{\z_U}x^{1/\k-2}U^2(x)\text{d}x.
\end{eqnarray}
By Williams' time reversal theorem (Fact
\ref{factWilliam}), there exists a $4$--dimensional Bessel
process $(R(s),\ s\geq 0)$, starting from $0$, such that
\begin{equation}
    \label{e1p31}
    (U(\z_U-s),\ s\leq \z_U)
\el
    (R(s),\ s\leq \g_a),
    \qquad a:=\sqrt{\l},
    \qquad
    \g_a:=\sup\{s\geq 0,\ R(s)=\sqrt{\l}\}.
\end{equation}
Therefore,
$$
N_\b\el\g_a^{-1/\k}\int_0^{\g_a}x^{1/\k-2}R^2(\g_a-x)\text{d}x
=\int_0^{1}(1-v)^{1/\k-2}\left(\frac{R(\g_a
v)}{\sqrt{\g_a}}\right)^2 \text{d}v.
$$
Recall (Yor \cite{Papier_3_Y11}, p.~52) that for any
bounded measurable functional $G$,
\begin{equation}\label{e2p31}
    \E\left[ G \left( \frac{R(\g_a u)}{\sqrt{\g_a}}, u\leq 1 \right) \right]
=
    \E \left( \frac{2}{R^2(1)} G\big(R(u),u\leq 1\big)\right).
\end{equation}
In particular, for $x>0$,
\begin{equation}\label{eqProbaDensite}
    \P\big(N_{\b}>x\big)
=
    \E\left( \frac{2}{R^2(1)} {\bf 1}_{ \{ \int_0^1
    (1-v)^{1/\k-2} R^2(v) \text{d} v>x \} } \right).
\end{equation}


\subsubsection{Proof of Theorem \ref{theoLiminfKinf1} (case $0<\k<1$)}

Fix $y>0$. By (\ref{eqProbaDensite}), for $r>1$,
\begin{eqnarray}
    \P(N_{\b}>y\log\log r)
& \le &
    \E \left( \frac{2}{R^2(1)}
            {\bf 1}_{ \{ \int_0^1 (1-v)^{1/\k-2} R^2(v) \text{d}v > y\log \log r ,\ R^2(1)
                     \le 1\} }
       \right)
\nonumber\\
&&
    + 2 \P\left( \int_0^1 (1-v)^{1/\k-2} R^2(v) \text{d}v >
    y\log\log r \right)
\nonumber\\
&:=&
    \Pi_1(r) + \Pi_2(r)
    \label{Pi1etPi2}
\end{eqnarray}
with obvious notation.

We first consider $\Pi_2(r)$.
Let
$
    {\mathcal{H}}
:=
    \big\{
        \big(t\in[0,1]\mapsto \int_0^t f(s)\text{d}s\big),
        \ f\in L^2\big([0,1],\R^4\big)
    \big\}
$.
As $R$ is the Euclidean norm
of a 4--dimensional Brownian motion $(\g(t),\ t\geq 0)$,
we have by Schilder's theorem (see e.g. Dembo and Zeitouni
\cite{Papier_3_DZ}, Thm. 5.2.3),
\begin{eqnarray}
 &&\lim_{r\to +\infty} \frac{1}{y\log\log r}
    \log \P \left( \int_0^1 (1-v)^{1/\k-2} R^2(v) \text{d}v >
    y\log\log r\right)
    \nonumber
    \\
 &=& -\inf\left\{\frac{1}{2}\int_0^1 \|\phi'(v)\|^2\text{d}v :
    \quad \phi\in {\mathcal {H}}, \; \int_0^1 (1-v)^{1/\k-2}
    \|\phi(v)\|^2 \text{d}v \geq 1 \right\}
    =:-c_1(\k),~~~~\hphantom{aa}
    \label{SchilderKinf1}
\end{eqnarray}
where $\| \cdot \|$ denotes the Euclidean norm.
For $\phi\in {\mathcal {H}}$,
$
    \|\phi(v)\|^2
=
    \big\|\int_0^v\phi'(u) \text{d}u\big\|^2
\leq
    v \int_0^1 \|\phi' (u)\|^2 \text{d}u
$,
where we applied Cauchy-Schwarz to each coordinate;
thus
$
    \int_0^1(1-v)^{1/\k-2} \|\phi(v)\|^2 \text{d} v
\le
    \big[\int_0^1(1-v)^{1/\k-2}v \text{d}v \big]
    \int_0^1 \|\phi' (v)\|^2 \text{d}v
$.
So, $c_1(\k)\in(0,\infty)$.

By \eqref{SchilderKinf1}, for $0<\e<1$ and large $r$,
\begin{equation}\label{PARTIEL}
    \Pi_2(r)
\leq
    (\log r)^{-(1-\e) y c_1(\k)}.
\end{equation}

Now, we consider $\Pi_1(r)$.
As $R$ is the Euclidean norm
of a 4--dimensional Brownian motion $(\g(t),\ t\geq 0)$, we have
$$
    \Pi_1(r)
=
    \E\left(
        \frac{2}{\|\g(1)\|^2}
        {\bf 1}_{ \{\|\g(1)\| \le 1\} }
        {\bf 1}_{ \{ \int_0^1(1-v)^{1/\k-2}
        \|\g(v) \|^2 \text{d}v > y \log \log r \}
    }\right) .
$$
By the triangular inequality, for any finite positive measure $\mu$ on
$[0,1]$,
$$
    \sqrt{\int_0^1 \|\g(v) \|^2 \text{d} \mu(v)}
\leq
    \sqrt{\int_0^1 \|\g(v) - v\g(1) \|^2 \text{d} \mu(v)}
    +
    \sqrt{\int_0^1 v^2 \text{d} \mu(v)} \, \|\g(1)\|.
$$
Therefore, applying this to
$\text{d} \mu(v)=(1-v)^{1/\k-2}\text{d}v$,
we have for large $r$,
\begin{equation*}
    \Pi_1(r)
 \le\E\left(\frac{2}{\|\g(1)\|^2} {\bf 1}_{ \{
    \int_0^1(1-v)^{1/\k-2} \|\g(v) -v\g(1) \|^2 \text{d}v >
    (\sqrt{y \log \log r}-c_{27})^2 \} }\right)
 :=\E\left(\frac{2}{\|\g(1)\|^2} {\bf 1}_{E}\right) ,
\end{equation*}
where $c_{27}:=\sqrt{\int_0^1 v^2 (1-v)^{1/\kappa-2}\text{d}v}$. By the independence of $\g(1)$ and $(\g(v)-v\g(1),\
v\in[0,1])$, the expectation on the right hand side is $=
\E \big(\frac{2}{\|\g(1)\|^2} \big) \P(E) = \P(E)$ (the last identity
being a consequence of (\ref{e2p31}) by taking $G=1$ there).
Therefore, $\Pi_1(r) \le \P(E)$.

\noindent Again, by the independence of $\g(1)$ and
$(\g(v)-v\g(1),\ v\in[0,1])$, we see that, by writing $c_{28}
:= 1/\P(\|\g(1)\|\le 1)$, $\Pi_1(r) \le c_{28}\, \P (E, \;
\|\g(1)\|\le 1)$. By another application of the triangular
inequality, this leads to, for large $r$:
$$
    \Pi_1(r)
\le
    c_{28} \,
    \P\left( \int_0^1(1-v)^{1/\k-2} \|\g(v) \|^2 \text{d}v
    >
    \Big(\sqrt{y \log \log r}-2c_{27}\Big)^2
    \right).
$$
In view of (\ref{SchilderKinf1}), we have, for all
large $r$,
$
    \Pi_1(r)
\leq
    (\log r)^{ -(1-\varepsilon) yc_1(\k)}
$.
Plugging this into \eqref{Pi1etPi2} and
\eqref{PARTIEL} yields that, for any $y>0$, $\varepsilon>0$ and
all large $r$,
\begin{equation}\label{InegQueueNbeta}
    \P(N_\b> y\log\log r)
\leq
    2(\log r)^{ -(1-\varepsilon)yc_1(\k)} .
\end{equation}



Let $0<\e<1/2$, and $s_n:=\exp(n^{1-\e})$.
We get
$
    \sum_{n=1}^{+\infty}
    \P \big\{\frac{H(F(s_n))}{\lm[H(F(s_n))]} > \frac{(1+4\e)\log\log s_n}{(1-\e)^3\k^2c_1(\k)} \big\}
<
    \infty
$
due to  \eqref{eqCompQuotientHL2} and \eqref{InegQueueNbeta}.
By the Borel--Cantelli lemma,
almost surely, for all large $n$,
\begin{equation}
    \label{conclusion1a}
    \frac{H(F(s_n))}{\lm[H(F(s_n))]}
\leq
    \frac{1+4\e}{(1-\e)^3\k^2 c_1(\k)}
    \log\log s_n.
\end{equation}

We now bound $\frac{H(F(s_{n+1}))}{H(F(s_n))}$. Observe that for
large $n$, $s_{n+1}-s_n\leq n^{-\e}s_n$. By Lemma
\ref{lemmaEncadrementFdeR}, almost surely for all large $n$,
\begin{eqnarray}
    H[F(s_{n+1})] - H[F(s_n)]
& \leq &
    H\big[\big(1+5  s_{n+1}^{-\d_0}/\k\big)s_{n+1}\big]
    -H\big[\big(1-5  s_n^{-\d_0}/\k\big)s_n\big]
\nonumber\\
& \leq &
    H\big[\big(1-5  s_n^{-\d_0}/\k\big)s_n+(2-\e)n^{-\e}s_n\big]
    -
    H\big[\big(1-5   s_n^{-\d_0}/\k\big)s_n\big]
\nonumber\\
 & = &  \inf \left\{ u\ge 0: \; {\widehat X}_n(u) > (2-\e)n^{-\e}s_n
 \right\},\label{Hn222}
\end{eqnarray}
where $\big({\widehat X}_n(u), \; u\ge 0\big)$ is, conditionally on $\wk$,  a diffusion
process in the random potential
$
    {\widehat W}_\k (x)
:=
    \wk\big[x+\big(1-\frac{5}{\k}s_n^{-\d_0}\big)s_n\big]
    -
    \wk\big[\big(1-\frac{5}{\k}s_n^{-\d_0}\big)s_n\big]
$,
$x\in \R$, starting from $0$.
We denote by ${\widehat H}_n(r)$ the hitting time of $r\geq 0$ by ${\widehat X}_n$, so that
\begin{equation}
    \label{Hn1}
    \inf \left\{ u\ge 0: \; {\widehat X}_n(u) > (2-\e)n^{-\e}s_n
    \right\}
 =
    {\widehat H}_n \big[ (2-\e)n^{-\e}s_n\big].
\end{equation}
Note that for any $r>0$, under $\P$, ${\widehat
H}_n(r)$ is distributed as $H(r)$. Therefore, applying
\eqref{e1p42} and Lemma \ref{lemmaEncadrementFdeR} to $r=
2n^{-\varepsilon}s_n$ yields that, for any
$0<\d_0< \frac{1}{2}$,
$$
    \sum_n
    \P\left[ {\widehat H}_n \left( \big[1- 5 (2n^{-\varepsilon}s_n)^{-\d_0}/\k\big] 2n^{-\varepsilon}s_n \right)
            >
            \big[n (\log n)^{1+\varepsilon} t_+(2n^{-\varepsilon}s_n)\big]^{1/\k}
       \right]
<
    \infty.
$$
Since $
    \big[1- \frac{5}{\k} (2n^{-\varepsilon} s_n)^{-\d_0}\big]
    2n^{-\varepsilon}s_n
\ge
    (2-\varepsilon) n^{-\varepsilon}s_n
$
(for large $n$), it follows from the
Borel--Cantelli lemma that, almost surely for all large $n$,
$
    {\widehat H}_n \left[ (2-\varepsilon) n^{-\varepsilon}s_n \right]
\le
    \big[n (\log n)^{1+\varepsilon} t_+(2
    n^{-\varepsilon}s_n) \big]^{1/\k} .
$
This, together with (\ref{Hn222}) and (\ref{Hn1}),
yields that, almost surely for all large $n$,
$$
    H[F(s_{n+1})] - H[F(s_n)]
\le
    \big[n (\log n)^{1+\varepsilon} t_+(2 n^{-\varepsilon}s_n) \big]^{1/\k}
\le
    c_{29} \big[n^{1-\varepsilon}(\log n)^{1+\varepsilon} s_n \big]^{1/\k} .
$$
Recall from Lemma \ref{lemmaEncadrementFdeR} and
Theorem \ref{thLowerLevyKinf1} that, almost surely, for all
large $n$, $H[F(s_n)]\ge H[(1-\e)s_n]
\ge \frac{c_{30} s_n^{1/\k}}{(\log\log s_n)^{1/\k-1}}$, which
yields
$$
    \frac{H[F(s_{n+1})]}{H[F(s_n)]}
\leq
    1+ \frac{c_{29} [n^{1-\varepsilon}(\log n)^{1+\varepsilon} s_n ]^{1/\k}}
            {c_{30} s_n^{1/\k}/(\log\log s_n)^{1/\k-1}}
\leq
    c_{31} (\log s_n)^{1/\k}
    (\log\log s_n)^ {(2+\e)/\k-1}.
$$
In view of \eqref{conclusion1a}, this yields that,
almost surely, for large $n$ and $t\in [H(F(s_n)), \,
H(F(s_{n+1}))]$,
$$
    \frac{t}{\lm(t)}
\le
    \frac{H[F(s_n)]}{\lm[H(F(s_n))]}
    \frac{H[F(s_{n+1})]}{H[F(s_n)]}
<
    c_{32}(\log s_n)^{1/\k}
    (\log\log s_n)^{(2+\e)/\k}.
$$
Since, almost surely for all large $n$,
$
    \log H[F(s_n)]
\geq
    \log H[(1-\e)s_n]
\geq
    \frac{1-\varepsilon}{\k} \log s_n
$
(this is seen first by Lemma \ref{lemmaEncadrementFdeR},
and then by Theorem \ref{thLowerLevyKinf1}), we have proved that
$$
\liminf_{t\to+\infty}\frac{\lm(t)}{t(\log t)^{-1/\k}(\log\log
t)^{-(2+\varepsilon)/\k}}
\ge
    c_{33}
\qquad \P\text{--a.s.}
$$
Since $\varepsilon\in (0, \frac{1}{2})$ is arbitrary,
this proves Theorem \ref{theoLiminfKinf1} in the case
$0<\k<1$.\hfill$\Box$


\subsubsection{Proof of Theorem \ref{theoLiminfKsup1} (case $0<\k<1$)}

By \eqref{eqProbaDensite}, for any $s>0$ and $u>0$,
\begin{eqnarray*}
    \P(N_{\b}>s)
&\ge&
    \frac{2}{u} \, \P \left( \, \int_0^1(1-v)^{1/\k-2}
    R^2(v)\text{d}v > s, \; R^2(1)\le u\right)
\\
&\ge&
    \frac{2}{u} \, \P \left( \, \int_0^1(1-v)^{1/\k-2}
    R^2(v)\text{d}v >  s\right) - \frac{2}{ u} \, \P
    \left( R^2(1)> u\right) .
\end{eqnarray*}
The first probability term on the right hand side is
taken care of by (\ref{SchilderKinf1}), whereas for the second,
we have
$
\frac{1}{u} \log \P( R^2(1)> u ) \to -\frac{1}{2}
$,
for $u\to \infty$,
since $R^2(1)$ has a chi-squared distribution with $4$ degrees of freedom.
Taking $u:= \exp( \sqrt{\log \log r}\,)$ leads
to: for any $y>0$,
$$
    \liminf_{r\to \infty}
    \frac{\log\P ( N_{\b}> y \log \log r)}{\log\log r}
\geq
    - y c_1(\k).
$$
Plugging this into \eqref{eqCompQuotientHL2} yields that, for
$r_n := \exp(n^{1+\varepsilon})$,
$$
    \sum_{n\geq 1}
    \P\left(
        \frac{(H\circ F )(r_{2n})} {(\lm\circ H\circ F )(r_{2n}) }
    >
        \frac{(1-3\e) \log\log r_{2n}}{\k^2c_1(\k)(1+\e)^3} -\e
    \right)
=
    +\infty.
$$
Let $Z_n := \sum_{k=1}^n r_k$. By Lemma
\ref{lemmaBorelCantelli} (in its notation), almost surely,
for infinitely many $n$,
\begin{equation}\label{eqZ2}
    \sup_{u\in[(1-\e)r_{2n},(1+\e)r_{2n}]}
    \frac{H_{X\circ\Theta_{H(Z_{2n-1})}}(u)}{(\lE\circ H)_{X\circ\Theta_{H(Z_{2n-1})}}(u) }
>
    \frac{(1-8\e) \log\log r_{2n}}
    {\k^2 c_1(\k)},
\end{equation}
if $\e>0$ is small enough.
Observe that
\begin{equation}\label{juin2526}
    (\lE\circ    H)_{X\circ\Theta_{H(Z_{2n-1})}}(u)
=
    \sup_{x\in\R}L_{\widetilde{X}_n}\big(\widetilde{H}_n(u),x\big)
=:
    L_{\widetilde{X}_n}^*\big(\widetilde{H}_n(u)\big),
\end{equation}
where $\big(\widetilde{X}_n(v),\ v\geq 0\big)$ is a diffusion
process in the random potential $\wk(x+Z_{2n-1})-\wk(Z_{2n-1})$, $x\in\R$,
$\big(L_{\widetilde{X}_n}(t,x),\ t\geq 0,\ x\in\R\big)$ is its
local time and
$\widetilde{H}_n(r):=\inf\big\{t>0,\ \widetilde{X}_n(t)>r\big\}$, $r> 0$. Hence, for any $u>0$, under $\P$,
the left hand side of \eqref{juin2526} is distributed as
$\lm(H(u))$. Applying  \eqref{eqProbaMajorL1} and Lemma
\ref{lemmaEncadrementFdeR} to $\widetilde{r}_{2n} :=
(1-\e)^2r_{2n}$, there exists $c_{34}>0$ such that
$$
    \sum_n
    \P\left[L_{\widetilde{X}_n}^*\left(\widetilde{H}_n\big[\big(1+5(\widetilde{r}_{2n})^{-\d_0}/\k\big)\widetilde{r}_{2n}\big]\right)<c_{34}
[r_{2n}/\log\log r_{2n}]^{1/\k}\right]
 <\infty.
$$
Since
$
    \big(1+\frac{5}{\k}(\widetilde{r}_{2n})^{-\d_0}\big)\widetilde{r}_{2n}
\leq
    (1-\e)r_{2n}
$
for large $n$, the \bc lemma gives that, almost surely, for all large $n$,
\begin{equation}\label{juin17}
    c_{34} \big[r_{2n}/\log\log r_{2n}\big]^{1/\k}
\leq
    L_{\widetilde{X}_n}^*\left(\widetilde{H}_n[(1-\e)r_{2n}]\right)
\leq
    L_{\widetilde{X}_n}^*\left( \widetilde{H}_n(u)\right)
\end{equation}
for any $u\in[(1-\e)r_{2n},(1+\e)r_{2n}]$. Applying Theorem
\ref{lemmaLimsupL0}, we have almost surely for large $n$,
\begin{equation*}
    \lm[H(Z_{2n-1})]
  \leq  [Z_{2n-1}\log^2 Z_{2n-1}]^{1/\k}
  \leq  \e [r_{2n}/\log\log r_{2n}]^{1/\k}
  \leq  (\e/c_{34}) L_{\widetilde{X}_n}^*\left(
    \widetilde{H}_n(u)\right)
\end{equation*}
for $u\in[(1-\e)r_{2n},(1+\e)r_{2n}]$, since $Z_k\le k
\exp(-k^\e)r_{k+1}$ for large $k$. Hence,
\begin{equation}
    \label{e3p36}
    \lm[H(Z_{2n-1}+u)]\leq (1+\e/c_{34})L_{\widetilde{X}_n}^*\left( \widetilde{H}_n(u)\right).
\end{equation}
On the other hand,
we have by Theorem \ref{ThUpperLevyKinf1}, almost surely, for all large
$n$,
$$
\log\log r_{2n}\geq (1-\e) \log\log H(Z_{2n-1}+u),\qquad
u\in[(1-\e)r_{2n},(1+\e)r_{2n}].
$$
Consequently, almost surely for infinitely many $n$, by
\eqref{e3p36} and \eqref{eqZ2},
\begin{eqnarray*}
 & & \inf_{v\in[Z_{2n-1} + (1-\e)r_{2n},
    Z_{2n-1} +(1+\e)r_{2n}]} \frac{\lm[H(v)]}{H(v)/\log\log H(v)}
    \\
 & \leq & (1+c_{35}\e) \inf_{u\in[(1-\e)r_{2n},(1+\e)r_{2n}]}
    \frac{(\lE\circ H )_{X\circ\Theta_{H(Z_{2n-1})}}(u)}
    {H_{X\circ\Theta_{H(Z_{2n-1})}}(u)/\log\log r_{2n}}
  \leq  (1+c_{36}\e)\k^2c_1(\k),
\end{eqnarray*}
proving Theorem \ref{theoLiminfKsup1} in the case
$0<\k<1$.\hfill$\Box$


\subsubsection{Proof of Theorem \ref{theoLimsupDifferentRWRE}}

Assume  $0<\k<1$. Fix $x>0$, and let $r_n:=\exp(n^{1+\e})$.
Since $\P(N_\b<x)>0$, \eqref{eqCompQuotientHL2} implies
$
    \sum_{n\in\N}
    \P\left(
        \frac{(H\circ F )(r_{2n})} {(\lm\circ H\circ F)(r_{2n}) }
    <
        \frac{(1+3\e)x}{\k^2} +\e
    \right)
=
    +\infty
$.
By Lemma \ref{lemmaBorelCantelli}, for small $\e>0$, almost surely for
infinitely many $n$,
\begin{equation}
    \label{Thm1.7}
    \inf_{u\in[(1-\e)r_{2n}, \ (1+\e)r_{2n}]}
    \frac{H_{X\circ\Theta_{H(Z_{2n-1})}}(u)}
    {(\lE\circ H)_{X\circ\Theta_{H(Z_{2n-1})}}(u) }
<
    \frac{(1+3\e)x}{ \k^2}
    +\e .
\end{equation}
 With the same notation as in \eqref{juin2526},
$H_{X\circ\Theta_{H(Z_{2n-1})}}(u)=H(Z_{2n-1}+u)-H(Z_{2n-1})$ is
the hitting time $\widetilde{H}_n(u)$ of $u$ by the diffusion
$\widetilde{X}_n$. For any $u$, under $\P$, it has the same
distribution as $H(u)$. Hence, applying \eqref{eqH+inf1} and Lemma
\ref{lemmaEncadrementFdeR} to $\widetilde{r}_{2n}=(1-\e)^2 r_{2n}$
leads to (for $0<\d_0<1/2$)
$$
    \sum_n
    \P\left[
        \widetilde{H}_n\left(\Big(1+\frac{5}{\k}(\widetilde{r}_{2n})^{-\d_0}\Big) \widetilde{r}_{2n}\right)
    <
        r_{2n}^{1/\k}/\log r_{2n}
    \right]
<
    \infty.
$$
Since $\big(1+\frac{5}{\k}(\widetilde{r}_{2n})^{-\d_0}\big)
\widetilde{r}_{2n}<(1-\e)r_{2n}$ for large $n$, it follows from
the \bc lemma that, almost surely, for all large $n$,
\begin{equation}\label{InegHuIssudeR}
    \frac{r_{2n}^{1/\k}}{\log r_{2n}}
\leq
    \widetilde{H}_n[(1-\e)r_{2n}]
\leq
    \inf_{u\in[(1-\e)r_{2n},(1+\e)r_{2n}]}H_{X\circ\Theta_{H(Z_{2n-1})}}(u).
\end{equation}
On the other hand, by Theorem \ref{ThUpperLevyKinf1},
$
    H(Z_{2n-1})
\leq
    [Z_{2n-1}\log^2 Z_{2n-1}]^{1/\k}
\leq
    \e \frac{r_{2n}^{1/\k}}{\log r_{2n}}
$
almost surely, for all large $n$.
This and \eqref{InegHuIssudeR} give, for $u \in [(1-\e)r_{2n}, \ (1+\e)r_{2n}]$,
$H(Z_{2n-1}+u)\le (1+\e)H_{X\circ\Theta_{H(Z_{2n-1})}}(u)$.
Plugging this into (\ref{Thm1.7}) yields that, almost surely,
for infinitely many $n$,
$$
    \inf_{u\in[(1-\e)r_{2n}, \ (1+\e)r_{2n}]}
    \frac{H(Z_{2n-1}+u)}{\lm(H(Z_{2n-1}+u))}
<
    \frac{(1+\e)(1+3\e)x}{\k^2} +\e(1+\e) .
$$
Hence
$
    \limsup_{t\to+\infty}\frac{\lm(t)}{t}
\geq
    \frac{\k^2}{x}
$,
a.s. Sending $x\to 0$
completes the proof of Theorem
\ref{theoLimsupDifferentRWRE}.\hfill$\Box$


\subsection{Case $\k=1$}\label{SubSec63}

This section is devoted to the proofs of Theorems
\ref{theoLiminfKsup1} and \ref{theoLiminfKinf1} in the case
$\k=1$ (thus $\l=8$; since $\l=4(1+\k)$).
Let
$$
    N_\b(t)
:=
    \frac{1}{\sup_{0\leq u\leq \tau_{\b}(8)}\b(u)}
    \left[ \,
        \int_0^1
        \frac{L_{\b}(\tau_{\b}(8),x) -8}{x} \text{d}x + \int_1^{+\infty}
        \frac{L_\b(\tau_{\b}(8), x)}{x} \text{d}x + 8\log t
    \right].
$$
Exactly as in (\ref{eqCompQuotientHL2}), we have, for
some $\alpha_1>0$, any $\varepsilon\in (0, 1/3)$, and all
large $r$,
\begin{equation}
    \label{eqCompQuotientHLCasKeg1}
    \P \left( (1-3\e)N_{\beta_{t_-(r)}}[t_-(r)] \le
    \frac{H(F(r))}{\lm[H(F(r))]} \le
    (1+3\e)N_{\beta_{t_+(r)}}[t_+(r)] \right)
\geq
    1-
    \frac{1}{r^{\alpha_1}} ,
\end{equation}
where $t_\pm (\cdot)$ are defined in (\ref{e1p22}), and $C_\beta$ in \eqref{e2p18}.
(Compared to (\ref{eqCompQuotientHL2}), we no longer have the
extra ``$\pm \e$'' terms, since they are already taken care of by
the presence of $8\log t$ in the definition of $N_\b(t)$).

With the same notation as in \eqref{e1p32} and \eqref{e3p31}, the
second Ray--Knight theorem (Fact \ref{FactRN2}) gives
\begin{eqnarray}
    N_\b(t)
&=&
    \frac{1}{\z_U} \left[ \, \int_0^1 \frac{U^2(x)-8}{x}
    \text{d}x + \int_1^{+\infty} \frac{U^2(x)}{x} \text{d}x +
    8\log t \right]
\label{EqNbetatkappa1Bis}
\\
&=&
    \frac{1}{\z_U} \left[ \, \int_0^{\z_U} \frac{U^2(x)-8}{x}
    \text{d}x + 8\log\z_U + 8\log t \right],
\label{EqNbetatkappa1}
\end{eqnarray}
since $U(x)=0$ for every $x\geq \z_U$.


\subsubsection{Proof of Theorem \ref{theoLiminfKinf1} (case
$\k=1$)}


We have $\l=8$ in the case $\k=1$.
Since $\sup_{x>0} \frac{\log x}{x}<\infty$, we have
\begin{equation}\label{InegNbeta1}
    N_\b(t)
\leq
    c_{37}
    +
    \frac{1}{\z_U} \, \int_0^{\z_U} \frac{|U^2(x)-8|}{x} \text{d}x
    + \frac{8\log t}{\z_U}.
\end{equation}
We claim that for some constant $c_{38}>0$,
\begin{equation}
    \label{U-LG}
    \limsup_{y\to +\infty} \frac{1}{y} \log\P\left(
    \frac{1}{\z_U} \, \int_0^{\z_U} \frac{|U^2(x)-8|}{x}
    \text{d}x > y \right) \le -c_{38} .
\end{equation}
Indeed, $\z_U = \sup_{0\le u\le \tau_\beta(8)} \beta(u)$ by
definition (see \eqref{e3p31}), which, in view of
\eqref{borodin}, implies that
$
\P( \z_U > z) = 1-e^{-4/z} \
\leq
4/ z
$
for $z>0$. Therefore, if we write $p(y)$ for the
probability expression at (\ref{U-LG}), we have, for any $z>0$,
$$
    p(y)
\leq
    \frac{4}{z}
    + \P\left( \frac{1}{\z_U} \, \int_0^{\z_U}
    \frac{|U^2(x)-8|}{x} \text{d}x >  y, \; \z_U \le z\right) .
$$
In the notation of \eqref{e1p31}--\eqref{e2p31}, this
yields
\begin{eqnarray}
    p(y)
 &\le&
    \frac{4}{z} + \P\left( \frac{1}{\g_a} \int_0^{1} \frac{|R^2
    (\g_a v)-8|}{1-v} \text{d}v >y, \; \gamma_a \le z \right)
    \nonumber
    \\
 &=&
    \frac{4}{ z}
    + \E\left( \frac{2}{R^2(1)} {\bf 1}_{ \{ \int_0^1
    \frac{|R^2(v)-R^2(1)|}{1-v} \text{d}v > y, \; R^2(1) \ge 8/z\}
    }\right)
    \nonumber
    \\
 &\le&
    \frac{4}{ z}
    + \frac{z}{ 4} \, \P\left( \, \int_0^1
    \frac{|R^2(v)-R^2(1)|}{1-v} \text{d}v > y \right).
    \label{Thm1.5}
\end{eqnarray}

In order to apply Schilder's theorem as in (\ref{SchilderKinf1}),
let $\phi\in {\mathcal{H}}$. As before between \eqref{SchilderKinf1} and \eqref{PARTIEL},
we have
$
    \|\phi(t)\|
\leq
    \sqrt{t} \, \big[\int_0^1 \|\phi' (s)\|^2\text{d}s \big]^{1/2}
$.
Similarly,
$
    \big|
        \|\phi(u)\|-\|\phi(1)\|
    \big|
\leq
    \|\phi(u)-\phi(1)\|
\leq
    \sqrt{1-u}
    \, \big[\int_0^1 \|\phi' (s)\|^2\text{d}s \big]^{1/2}
$.
Hence,
\begin{eqnarray*}
    \int_0^1 \frac{\big|\|\phi(u)\|^2-\|\phi(1)\|^2\big|}{1-u} \text{d}u
& = &
    \int_0^1 \frac{\big|\|\phi(u)\|-\|\phi(1)\|\big|}{1-u}
             \big[\|\phi(u)\| +\|\phi(1)\|\big] \text{d}u
\\
& \leq &
    2\left(\, \int_0^1 \frac{\text{d}u}{\sqrt{1-u}} \right)
    \int_0^1 \|\phi' (s)\|^2 \text{d}s.
\end{eqnarray*}
Consequently,
$$
    c_{39}
:=
    \inf \bigg\{ \frac{1}{2} \int_0^1 \big\|\phi'(u)\big\|^2 \text{d}u :
\quad
    \phi\in {\mathcal{H}}
    ,\
    \int_0^1 \frac{\big|\|\phi(u)\|^2-\|\phi(1)\|^2\big|}{1-u} \text{d}u > 1 \bigg\}
\in
    (0,\infty).
$$
Applying Schilder's theorem gives that
$
\limsup_{y\to +\infty}
\frac{1}{ y} \log\P\left( {\int_0^1
\frac{|R^2(v)-R^2(1)|}{1-v} \text{d}v> y}\right) \le -c_{39}
$.
 Plugging this into (\ref{Thm1.5}), and taking $z=
\exp( \frac{c_{39}}{ 2} y)$ there, we obtain the claimed
inequality in (\ref{U-LG}), with $c_{38} := c_{39}/ 2$.

On the other hand, by \eqref{borodin} and \eqref{e3p31},
$$
\P\left( \frac{8\log t}{\z_U} > 2(1+2\e) (\log t)\log \log t
\right) =\frac{1}{(\log t)^{1+2\e}} .
$$
This, combined with \eqref{InegNbeta1} and \eqref{U-LG} gives, for all large $t$,
$$
\P\left\{ N_\b(t) > 2(1+3\e) (\log t)\log \log
t\right\} \le \frac{2}{(\log t)^{1+2\e}}.
$$

Let $s_n:=\exp(n^{1-\e})$. By \eqref{eqCompQuotientHLCasKeg1},
$
\sum_{n=1}^{+\infty} \P \left( \frac{H(F(s_n))}{\lm[H(F(s_n))]}
> 2(1+3\e)^2 (\log s_n) \log\log s_n \right) < \infty
$,
 which, by means of the Borel--Cantelli lemma,
implies that, almost surely, for all large $n$,
\begin{equation}
\frac{H(F(s_n))}{\lm[H(F(s_n))]} \le 2(1+3\e)^2 (\log s_n)
\log\log s_n .
\label{exp31}
\end{equation}
Now we give an upper bound for
$\frac{H(F(s_{n+1}))}{H(F(s_n))}$. By Lemma
\ref{lemmaEncadrementFdeR}, almost surely for $n$ large enough,
$F(s_n)\geq (1-\e)s_n$. An application of Theorem
\ref{thLowerLevyKinf1} yields that, almost surely, for large $n$,
\begin{equation}
    \label{pape2K1}
    H[F(s_n)]\geq H[(1-\e)s_n]\geq 4(1-2\e)s_n\log s_n.
\end{equation}
With the same notation and the same arguments as in \eqref{Hn222}
and \eqref{Hn1}, almost surely for all large $n$,
$
    H[F(s_{n+1})]-H[F(s_n)]
\leq
    \widehat{H}_n[(2-\e)n^{-\e}s_n]
$.
Moreover, $\widehat{H}_n(r)$ is distributed as $H(r)$ under $\P$
for any $r>0$. Hence, applying Lemma \ref{lemmaEncadrementFdeR}
and \eqref{littlehope} to $r=\widetilde{s}_n:=2n^{-\e}s_n$ and
$a(e^{-2}\widetilde{s}_n)=8n(\log n)^{1+\e}$  for
$0<\d_0<\frac{1}{2}$, we get
$$
    \sum_n\P\left[\widehat{H}_n\left((1-5(\widetilde{s}_n)^{-\d_0}/\k)\widetilde{s_n}\right)
            >
                32(1+\e)t_+(\widetilde{s}_n)
                \big[c_3+n(\log n)^{1+\e}+\log t_+(\widetilde{s}_n)\big]
            \right]
<
    \infty.
$$
Since
$[1-\frac{5}{\k}(\widetilde{s}_n)^{-\d_0}]\widetilde{s}_n\geq
(2-\e)n^{-\e}s_n$ (for large $n$), the \bc lemma yields that
$$
    \widehat{H}_n\big((2-\e)n^{-\e}s_n\big)
\leq
    32(1+\e)t_+(2n^{-\e}s_n)[c_3+n(\log n)^{1+\e}+\log t_+(2n^{-\e}s_n)],
$$
almost surely for large $n$.
Hence,
$
    H[F(s_{n+1})]-H[F(s_n)]
\leq
    c_{39}s_n(\log s_n)(\log n)^{1+\e}
$.
Hence, by \eqref{pape2K1}, we have, almost surely, for all
large $n$,
\begin{equation*}
    H[F(s_{n+1})]/H[F(s_n)]
    \leq c_{40}(\log\log s_n)^{1+\e}.
\end{equation*}
Let $t\in [H(F(s_n)), \, H(F(s_{n+1}))]$. By \eqref{exp31}, if $t$ is large enough,
$$
    \frac{t}{\lm(t)}
\leq
    \frac{H[F(s_n)]}{\lm[H(F(s_n))]}\frac{H[F(s_{n+1})]}{H[F(s_n)]}
<
    3c_{40}(\log s_n)(\log\log s_n)^{2+\e}.
$$
Since almost surely for large $n$,
$
    \log H[F(s_n)]
\ge
    \log H[(1-\e)s_n]
\ge
    \log s_n
$ (by Lemma
\ref{lemmaEncadrementFdeR} and Theorem \ref{thLowerLevyKinf1}),
this yields
$$
\liminf_{t\to+\infty}\frac{\lm(t)}{t/[(\log t)(\log\log
t)^{2+\e}]}\geq \frac{1}{3c_{40}}\qquad\P\text{--a.s.}
$$
Theorem \ref{theoLiminfKinf1} is proved in the case
$\k=1$.\hfill$\Box$


\subsubsection{Proof of Theorem \ref{theoLiminfKsup1} (case $\k=1$)}

Again, $\l=8$. Let $0<\e<1/2$.
Recall that $\z_U=\sup_{0\leq u\leq \tau_\b(8)}\b(u)$,
and that
$
N_\beta(t) = \frac{1}{\z_U} \left[ \, \int_0^{\z_U}
\frac{U^2(x)-8}{x} \text{d}x + 8\log\z_U + 8\log t \right]
$ (see \eqref{e3p31} and \eqref{EqNbetatkappa1}).
This time, we need to bound $N_\beta(t)$ from below.
By \eqref{borodin} for large $z$,
$$
    \P \left( 8\frac{\log\z_U}{\z_U} <- z \right)
\leq
    \P\left( -\frac{1}{\z_U^2} < -z \right)
=
    \P\left( \z_U < \frac{1}{\sqrt{z}} \right)
=
    \exp\left( - 4\sqrt{z}\right).
$$
By \eqref{borodin} again,
$$
\P\left(\frac{8\log t}{\z_U}>2(1-\e)(\log t)\log \log
t\right) = \frac{1}{(\log t)^{1-\e}}.
$$
On the other hand, for all large $y$,
$
    \P \big(\frac{1}{\z_U} \, \int_0^{\z_U} \frac{|U^2(x)-8|}{x} \text{d}x >y \big)
\le
    e^{-c_{41} y}$ (see (\ref{U-LG})). Assembling these
pieces yields that, for all large $t$,
$$
\P[N_\b(t)>2(1-2\e)(\log t)\log \log t]\geq\frac{1}{2(\log
t)^{1-\e}}.
$$


Let $r_n:=\exp(n^{1+\e})$. In view
of \eqref{eqCompQuotientHLCasKeg1} and Lemma \ref{lemmaBorelCantelli}, we get almost surely for infinitely many $n$,
%
\begin{equation}
    \label{eqZ2K1}
    \sup_{u \in[(1-\e)r_{2n},\, (1+\e)r_{2n}]} \frac{H_{X \circ
    \Theta_{H(Z_{2n-1})}} (u)} {(\lE\circ H)_{X \circ
    \Theta_{H(Z_{2n-1})}}(u) }
>
    2(1-2\e)(1-3\e)(\log r_{2n})
    \log\log r_{2n} .
\end{equation}
The expression on the left hand side of \eqref{eqZ2K1}
is ``close to'' $H(r_{2n})/\lm[H(r_{2n})]$, but we need to prove
this rigorously. With the same argument as in the displays between
\eqref{juin2526} and \eqref{juin17}, we get that there exists $c_{42}>0$ such that, almost surely
for large $n$,
\begin{equation*}
    \inf_{u\in[(1-\e)r_{2n},(1+\e)r_{2n}]}(\lE\circ H)_{X\circ\Theta_{H(Z_{2n-1})}}(u)
\geq
     c_{42} r_{2n}/\log\log r_{2n}.
\end{equation*}
Observe that $Z_k\le k \exp(-k^\e)r_{k+1}$ for large $k$, as in the paragraph after \eqref{eqPourCiterAutreChose}.
Exactly as in the case $0<\k<1$, we apply Theorem \ref{lemmaLimsupL0}, to
see that almost surely for large $n$,
\begin{equation*}
\lm[H(Z_{2n-1})]
  \leq
\e r_{2n}/\log\log r_{2n}
 \leq  (\e/c_{42}) \inf_{u\in[(1-\e)r_{2n},(1+\e)r_{2n}]}(\lE\circ
H
)_{X\circ\Theta_{H(Z_{2n-1})}}(u) ,
\end{equation*}
which implies, for all $u\in[(1-\e)r_{2n},(1+\e)r_{2n}]$,
\begin{equation} \label{e3p36K1}
    \lm[H(Z_{2n-1}+u)]\leq (1+\e/c_{42})(\lE\circ H
    )_{X\circ\Theta_{H(Z_{2n-1})}}(u).
\end{equation}
By Theorem \ref{ThUpperLevyKinf1}, almost surely for all
large $n$, $\sup_{u\in[(1-\e)r_{2n},(1+\e)r_{2n}]} \log
H(Z_{2n-1}+u) \le (1+\e) \log r_{2n}$. In view of \eqref{e3p36K1} and then \eqref{eqZ2K1},
there are almost surely infinitely many $n$ such that
\begin{eqnarray*}
& &
    \inf_{v\in[Z_{2n-1}+(1-\e)r_{2n},Z_{2n-1}+(1+\e)r_{2n}]}
    \frac{\lm[H(v)]} {H(v)/[(\log H(v))\log\log H(v)]}\\
& \leq &
    (1+c_{43}\e) \inf_{u\in[(1-\e)r_{2n},(1+\e)r_{2n}]}\frac{(L^*\circ H)_{X\circ\Theta_{H(Z_{2n-1})}}(u)}
    {H_{X\circ\Theta_{H(Z_{2n-1})}}(u) [(\log r_{2n})\log\log r_{2n}]^{-1}}
\leq
    (1+ c_{44}\e)/2.
\end{eqnarray*}
This proves Theorem \ref{theoLiminfKsup1} in the case
$\k=1$.\hfill$\Box$

\noindent{\bf Remark:} Assume $\k=1$.
We also prove that in this case, $\P$ almost surely,
\begin{equation}\label{InegLimsupCasK1}
    \limsup_{t\to+\infty}\lm(t)/t\
\geq
    8/[ c_{17}\pi]
=
    1/2.
\end{equation}
This is in agreement with Theorem 1.1 of Gantert and Shi \cite{Papier_3_GS1} for RWRE.
However, we could not prove whether this $\limsup$ is finite or not, contrarily to the cases $\k\in(0,1)$
and $\k>1$, and to the case of RWRE, for
which the maximum local time at time $t$ is clearly less than $t/2$.

We now prove \eqref{InegLimsupCasK1}.
With the same notation as in \eqref{e1p32} and \eqref{e3p31},
let
$
    {\widehat C}_U
:=
    \int_0^1 \frac{U^2(x)-8}{x} \text{d}x
    +
    \int_1^{+\infty} \frac{U^2(x)}{x} \text{d}x
$,
$\e\in(0,1/3)$ and $y:=(1+\e)^2 c_{17}\pi/[8(1-\e)]$.
We have for $z>0$, by \eqref{EqNbetatkappa1Bis},
$$
    \P[N_\b(t)<y]
=
    \P\big[{\widehat C}_U+8\log t<y \z_U\big]
\geq
    \P\big[(z+8)\log t<y \z_U, {\widehat C}_U\leq z \log t\big].
$$
Notice that
$
    {\widehat C}_U
\el
    C_\beta
 \el
    8c_3 +(\pi/2)C_8^{ca}
$
first by \eqref{e1p32} and
\eqref{e2p18}, then by Fact \ref{FactBianeYor}.
So,
$
    \P\big[{\widehat C}_U > z \log t\big]
=
    \P\big[ C_8^{ca} > (2z/\pi) \log t-16c_3/\pi\big]
\sim_{t\to+\infty}
    \pi c_{17}/(2z\log t)
$
(see before \eqref{littlehope}).
Moreover,
$
    \P[(z+8)\log t<y \z_U]
\sim_{t\to+\infty}
    4y/[(z+8)\log t]
$
by \eqref{e3p31} and \eqref{borodin}.
Thus,
\begin{eqnarray*}
    \P[N_\b(t)<y]
& \geq &
    \P\big[(z+8)\log t<y \z_U\big]
    -\P\big[ {\widehat C}_U > z \log t\big]
\\
& \geq &
    [4(1-\e)y/(z+8)-(1+\e)\pi c_{17}/(2z)]/\log t
\\
& = &
    (1+\e)c_{17}\pi[(1+\e)/(z+8)-1/z]/(2\log t)
\end{eqnarray*}
So we can choose $z$ so that
$
    \P[N_\b(t)<y]
\geq
    c_{45}/\log t
$
for some constant $c_{45}>0$.
We now set $r_k:=k^k$, $k\in\N^*$.
This and \eqref{eqCompQuotientHLCasKeg1} give for some $\a_1>0$,
\begin{eqnarray*}
    \P\bigg[\frac{\lm[H(F(r_{2n}))]}{H(F(r_{2n}))}
            >
            [(1+3\e)y]^{-1}
      \bigg]
& \geq &
    \P\bigg[N_{\beta}[t_+(r_{2n})]
            <
            y
      \bigg]
    -
    \frac{1}{r_{2n}^{\a_1}}
\geq
    \frac{c_{45}}{2\log r_{2n}}
\geq
    \frac{c_{45}}{5n\log n}
\end{eqnarray*}
for large $n$.
Hence by Lemma \ref{lemmaBorelCantelli}  in its notation,
$\P$ almost surely, there exist infinitely many $n$ such that
for some $t_n \in [(1-\e)r_{2n}, (1+\e)r_{2n}]$,
\begin{equation}\label{InegLHXK1}
    \frac
    {(\lE\circ H)_{X \circ\Theta_{H(Z_{2n-1})}}(t_n) }
    {H_{X \circ \Theta_{H(Z_{2n-1})}} (t_n)}
>
            [(1+3\e)y]^{-1}.
\end{equation}
Notice that
$
    Z_{2n-1}
=
    \sum_{k=1}^{2n-1} k^k
\leq
    (2n-1)^{2n-1}
    +
    \sum_{k=1}^{2n-2} (2n-2)^{2n-2}
\leq
    2(2n)^{2n-1}
=
    r_{2n}/n
$.
We have by Theorem \ref{ThUpperLevyKinf1}, almost surely for all large $n$,
\begin{equation}\label{InegHZK1}
    H(Z_{2n-1})
\leq
    Z_{2n-1} (\log Z_{2n-1}) (\log \log Z_{2n-1})^2
\leq
    \e r_{2n}\log r_{2n}
\end{equation}
On the other hand, first by Lemma \ref{lemmaEncadrementFdeR} and Lemma \ref{lemmaApproxLxLi}, then by \eqref{eqDefItierceKeg1}
and since $\E[\exp(-C_8^{ca})]=1$ as before \eqref{ProbaSamoTaqqu}, for every $\e>0$ small enough,
\begin{eqnarray*}
&&
    \P\big[H_{X \circ \Theta_{H(Z_{2n-1})}} ((1-\e)r_{2n})
           <
           (1-10\e) 32 t_-(r_{2n}) [\log t_-(r_{2n}) +c_3]
      \big]
\\
& \leq &
    \P\big[(1-\e) \Ib_-((1-2\e)r_{2n})
           <
           (1-10\e) 32 t_-(r_{2n}) [\log t_-(r_{2n}) +c_3]
      \big]
    +2r_{2n}^{-\a_1}
\\
& \leq &
    \P\big[ C_8^{ca}
           <
           -\e (16/\pi)  \log r_{2n}
      \big]
    +2 r_{2n}^{-\a_1}
\leq
    2r_{2n}^{-16\e/\pi}.
\end{eqnarray*}
Hence, thanks to the Borel Cantelli lemma, almost surely for large $n$,
\begin{eqnarray*}
    H_{X \circ \Theta_{H(Z_{2n-1})}} (t_n)
& \geq &
    H_{X \circ \Theta_{H(Z_{2n-1})}} ((1-\e)r_{2n})
\\
& \geq &
           (1-10\e) 32 t_-(r_{2n}) [\log t_-(r_{2n}) +c_3]
\geq
    (1-11\e) 4 r_{2n} \log r_{2n}.
\end{eqnarray*}
This together with \eqref{InegLHXK1} and \eqref{InegHZK1} gives $\P$ almost surely for infinitely many $n$,
\begin{eqnarray*}
    \frac{\lm[H(Z_{2n-1}+t_n)]}{H(Z_{2n-1}+t_n)}
& \geq &
    \frac{(\lE\circ H)_{X \circ\Theta_{H(Z_{2n-1})}}(t_n) }
         {H_{X \circ \Theta_{H(Z_{2n-1})}} (t_n)}
    \frac{H_{X \circ \Theta_{H(Z_{2n-1})}} (t_n)}
         {H(Z_{2n-1})+H_{X \circ \Theta_{H(Z_{2n-1})}} (t_n)}
\\
& \geq &
    [(1+3\e)y]^{-1} (1+\e)^{-1}
\geq
    (1-10\e)
    8/[ c_{17}\pi],
\end{eqnarray*}
for small $\e$.
As before, let $t\to +\infty$, and then $\e\to 0$. This proves \eqref{InegLimsupCasK1} since $c_{17}=\frac{16}{\pi}$
as before \eqref{littlehope}.


\mysection{Proof of Lemma \ref{lemmaApproxLxLi}}
 \label{sectannexe}

This section is devoted to the proof of Lemma
\ref{lemmaApproxLxLi}. The basic idea goes back to Hu et
al.~\cite{Papier_3_HSY1}, but requires considerable refinements
due to the complicated nature of the process
$x\mapsto L_X(t,x)$ and to the fact that we are interested in the joint law of
$\big(\lm[H(.)], H(.)\big)$.
Throughout the proof we consider the annealed probability $\P$.

Let $\k>0$ and $\e\in(0,1)$. We fix $r>1$.
Recall that $A(x)=\int_0^x e^{\wk(u)} \text{d}u$, and
$A_{\infty}=\lim_{x\to+\infty} A(x)<\infty$, a.s. As in Brox
(\cite{Papier_3_B3}, eq. (1.1)), the general diffusion theory leads to
\begin{equation}\label{eqChangeOfTimeAndScale}
    X(t)
=
    A^{-1}[B(T^{-1}(t))],\qquad t\ge 0,
\end{equation}
where $(B(s),\ s\geq 0)$ is a Brownian motion independent of $W$,
and for  $0\leq u <\sigma_B(A_{\infty})$,
$
T(u):=\int_0^u\exp\{-2\wk[A^{-1}(B(s))]\}\text{d}s$
 ($A^{-1}$ and $T^{-1}$ denote respectively the
inverses of $A$ and $T$). The local time of $X$ can be written as
(see Shi \cite{Papier_3_S4}, eq. (2.5))
\begin{equation}
    \label{tempsLocal}
    \lx(t,x)
=
    e^{-\wk(x)}\lB\big(T^{-1}(t),A(x)\big),
\qquad
    t\geq 0,\
    x\in\R.
\end{equation}

As in (\ref{H}), $H(\cdot)$ denotes the first hitting time of $X$.
Then as in Shi (\cite{Papier_3_S2}, eq. (4.3) to (4.6)),
\begin{equation}\label{eqHder}
    H(u)
=
    T[\s_B(A(u))]
=
    \int_{-\infty}^0+\int_0^{A(u)} e^{-2\wk[A^{-1}(x)]} L_B(\s_B[A(u)],x)
    \text{d}x
=:
    H_-(u)+H_+(u)
\end{equation}
for $u\geq 0$.
Recall $F$ from \eqref{e1p5} and notice that $F(r)>0$ on $\EW(r)$ if $r$ is large enough.
By scaling since $\wk$ and then $A(F(r))$ are independent of $B$, and then by
the first Ray--Knight theorem (Fact \ref{FactRN1}), there
exists a squared Bessel process of dimension $2$, starting from
$0$ and denoted by $\big(R_2^2(s),\ s\geq 0\big)$, independent of $\wk$,
such that
$$
    \left(\frac{\lB\{\s_B[A(F(r))],A(F(r))-s A(F(r))\}}{A(F(r))},\ s\in[0,1] \right)
=
    \left( R_2^2(s),\ s\in [0,1] \right).
$$
Hence, it is more convenient to study $\lm[H(.)]$
instead of $\lm(t)$. We consider
$$
    \lp[H(u)]
:=
    \sup\limits_{x\geq 0} L_X(H(u),x)
=
    \sup\limits_{0\leq x \leq u} \big\{ e^{-\wk(x)} \lB[\s_B(A(u)),A(x)] \big\},
\qquad u\geq 0.
$$
In particular,
$$
\lp[H(F(r))]=\sup_{x\in [0,F(r)]} \left\{e^{-\wk(x)}A(F(r))
R_2^2\left[\frac{A(F(r))-A(x)}{A(F(r))}\right]\right\}.
$$
Moreover, by Lamperti's representation theorem (Fact
\ref{FactLamperti}), there exists a Bessel process $\r=(\r(t),\
t\geq 0)$, of dimension $(2-2\k)$, starting from $\r(0)=2$, such
that
for all $t\geq 0$, $e^{\wk(t)/2}=\frac{1}{2}\r(A(t))$.
Now, let
$$
\rt(t):=\r(A_{\infty}-t), \qquad  0\leq t\leq A_{\infty}.
$$
By Williams' time reversal theorem (Fact
\ref{factWilliam}), $\rt$ is a Bessel process of dimension
$(2+2\k)$, starting from $0$.
Since $\wk$ and $A(F(r))$ are independent of
$R_2$, $u\mapsto \sqrt{A(F(r))}R_2(u/A(F(r)))$ is a
$2$--dimensional Bessel process, starting from $0$ and independent
of $\rt$. We still denote by $R_2$ this new Bessel process.
We obtain
\begin{equation*}
    \lp[H(F(r))]
 =  4\sup\limits_{x\in[0,F(r)]} \frac{R_2^2
    [A(F(r))-A(x)]}{\rt^2 [A_{\infty}-A(x)]}
    \\
  =  4\sup\limits_{v\in[0,A(F(r))]} \frac{R_2^2(v)}{\rt^2
    [A_{\infty}-A(F(r))+v]}.
\end{equation*}
Doing the same transformations on $H_+(F(r))$ (see \eqref{eqHder}) and recalling
that $A_{\infty}-A(F(r)) =\d(r)=\exp(-\k r/2)$ and so is deterministic thanks to the random function $F$, we obtain
\begin{eqnarray*}
& &
    \big(\lp[H(F(r))],H_+[F(r)]\big)
\\
& = &
    \left(4\sup\limits_{v\in[0,A(F(r))]}
    \frac{R_2^2(v)}{\rt^2[\d(r)+v]}, 16\int_0^{A(F(r))}
    \frac{R_2^2(s)}{\rt^4[\d(r)+s]}\text{d}s\right)\\
& = &
    \left(4\sup\limits_{u\in[0,\d(r)^{-1}A(F(r))]}
    \frac{R_2^2[\d(r)u]}{\rt^2[\d(r)(1+u)]},
    16\int_0^{\d(r)^{-1}A(F(r))}
    \frac{R_2^2[\d(r)u]\d(r)\text{d}u}{\rt^4[\d(r)(1+u)]}\right).
\end{eqnarray*}
We still denote by $R_2$ the $2$-dimensional Bessel
process $u\mapsto \frac{1}{\sqrt{\d(r)}}R_2(\d(r)u)$. We
define
\begin{equation}\label{eqImportanceFdeR}
    \rb(u)
=
    \frac{1}{\sqrt{\d(r)}}\rt[\d(r)(1+u)],\qquad u\geq 0.
\end{equation}
Notice that
$\rb(u)$
is a
$(2+2\k)$--dimensional Bessel process, starting from
$\rt(\d(r))/\sqrt{\d(r)}$ and independent of $R_2$.

Recall (see e.g. Karlin and Taylor \cite{Papier_3_KT2} p.~335) that a
Jacobi process $(Y(t),\ t\geq 0)$  of dimensions $(d_1,d_2)$ is a solution of the
stochastic differential equation
\begin{equation}
    \label{eqDefJacobi}
    \text{d} Y(t) = 2\sqrt{Y(t) (1-Y(t))} \, \text{d} \bhat(t) +
    \big[d_1-(d_1+d_2)Y(t)\big] \text{d} t,
\end{equation}
where $\big(\bhat(t),\ t\geq 0\big)$  is a standard Brownian motion.

According to Warren and Yor (\cite{Papier_3_WY1} p.~337), there exists a
Jacobi process $(Y(t),\ t\geq 0)$ of dimensions $(2,2+2\k)$,
starting from $0$, independent of $\big(R_2^2(t)+\rb^2(t), \ t\geq 0\big)$, such that
\begin{equation}
\forall u\geq 0, \qquad
\frac{R_2^2(u)}{R_2^2(u)+\rb^2(u)} = Y\circ \ty(u), \qquad
    \ty(u):=\int_0^u\frac{\text{d}s}{R_2^2(s)+\rb^2(s)}.
    \label{defThetaY}
\end{equation}
In particular, $(\ty(t),\ t\geq 0)$ is independent of
$Y$. As a consequence, for all $r\geq 0$,
\begin{eqnarray*}
 && (\lp[H(F(r))],H_+[F(r)])
    \\
 & = & \left(4\sup\limits_{u\in[0,\d(r)^{-1}A(F(r))]}
    \frac{Y\circ\ty(u)}{1- Y\circ\ty(u)},
    16\int_0^{\d(r)^{-1}A(F(r))} \frac{[Y\circ\ty(u)]
    \ty'(u) \text{d}u}{[1-Y\circ \ty(u)]^2} \right)
    \\
 &=& \bigg( 4 \sup\limits_{u\in[0,\varphi(r)]} \frac{Y(u)}{1-Y(u)}, 16
    \int_0^{\varphi(r)} \frac{Y(u)}{(1-Y(u))^2} \text{d}u \bigg),
\end{eqnarray*}
where
\begin{equation}
    \label{e10p7}
    \varphi(r)
:=
    \ty\big[\d(r)^{-1}A(F(r))\big].
\end{equation}

Define $\a_{\k}:= 1/(4+2\k)$ and let
$T_Y(\a_\k):=\inf\{t>0, Y(t)=\a_{\k}\}$
be the hitting time of $\a_\k$ by $Y$. We introduce
\begin{equation}
    \label{e1p7}
    \begin{array}{ll}
    \lbar(r)
:=
    4\sup\limits_{u\in[0,T_Y(\a_\k)]}
    \cfrac{Y(u)}{1-Y(u)},
&
\qquad
    \hbar(r)
:=
    16\int\limits_0^{T_Y(\a_\k)}
    \cfrac{Y(u)}{(1-Y(u))^2}\text{d}u,
\\
    L_0(r)
:=
    4\sup\limits_{u\in[T_Y(\a_\k),\varphi(r)]}
    \cfrac{Y(u)}{1-Y(u)},
& \qquad
    H_0(r)
:=
    16\int\limits_{T_Y(\a_\k)}^{\varphi(r)}
    \cfrac{Y(u)}{(1-Y(u))^2}\text{d}u.
    \end{array}
\end{equation}
We have on
$
E_{9}:=\{
    T_Y(\a_\k)
\leq
    64\log r
\leq
    \k  r/(2\l)
\leq
    \varphi(r)
\}
$,
\begin{equation}
    \left( \lp[H(F(r))] , \; H_+(F(r)) \right) =
    \left( \max\{\lbar(r),L_0(r)\}, \; \hbar(r)+H_0(r)\right),
    \label{e3p7}
\end{equation}
\begin{equation}
   \label{eqMajorAvantAlphaK}
    \lbar(r)
\leq
    \frac{4\a_\k}{1-\a_\k}
\qquad \text{and}\qquad
    \hbar(r)
\leq
    \frac{16\a_\k}{(1-\a_\k)^2}T_Y(\a_\k)
\leq
    \frac{2^{10}\a_\k}{(1-\a_\k)^2}\log r.
\end{equation}
Observe that
$
    S(y)
:=
    \int_{\a_\k}^y \frac{\text{d}x}{x(1-x)^{1+\k}}
$, $y\in(0,1)$
is a scale function of $Y$, as in Hu et al. (\cite{Papier_3_HSY1}, eq. (2.27)).
Hence $t\mapsto S[Y(t+T_Y(\alpha_\k))]$ is a continuous local martingale, so by Dubins-Schwarz theorem,  there exists a
Brownian motion $(\b(t),\ t\geq 0)$ such that for all $t\geq 0$,
\begin{equation}\label{eqDefBetaEtU}
    Y[t+T_Y(\a_\k)]=S^{-1}\{\b[U_Y(t)]\},
\qquad
    U_Y(t)
:=
    4 \int_0^t \frac{\text{d}s}{Y[s+T_Y(\a_\k)]
    \{1-Y[s+T_Y(\a_\k)] \}^{1+2\k}}.
\end{equation}

The rest of the proof of Lemma \ref{lemmaApproxLxLi} requires some
more estimates, stated as Lemmas
\ref{lemmaApproxTy}--\ref{remarquenegatif}
below. Lemmas \ref{lemmaApproxTy}--\ref{lemmaApproxJCasKinf1} deal only with Bessel processes, Jacobi processes and
Brownian motion, and may be of independent interest, whereas Lemma \ref{remarquenegatif}
gives an upper bound for the total time spent by $X$ on $(-\infty,0)$, and for the maximum local time of $X$ in $(-\infty, 0)$.
We defer the proofs of Lemmas \ref{lemmaApproxTy}--\ref{lemmaApproxJCasKinf1} to Section \ref{SectProofLemmasFin},
and we complete the proof of Lemma
\ref{lemmaApproxLxLi}.

\bigskip
\begin{lemma}\label{lemmaApproxTy}
Let $(R(t),\ t\geq 0)$ be a Bessel process of dimension $d>4$,
starting from $R_0\el\widetilde{R}_{d-2}(1)$, where
$\big(\widetilde{R}_{d-2}(t),\ t\in[0,1]\big)$ is $(d-2)$--dimensional
Bessel process. For any $\d_3\in (0, \, \frac{1}{2})$ and all
large $t$,
$$
    \P \left\{
        \left|
            \frac{1}{\log t} \int_0^{t}
            \frac{\textnormal{d}s}{R^2(s)} - \frac{1}{d-2}
        \right|
        >
        \frac{1}{(\log t)^{(1/2)-\d_3}}
    \right\}
\leq
\exp \left( - c_{46}\, (\log
t)^{2\d_3} \right).
$$
\end{lemma}

\bigskip
\begin{lemma}\label{lemmaApproximationU}
Let  $\d_1>0$ and define
\begin{equation}
    \EG
:=
    \left\{
        \tau_{\beta}\big[\big(1-v^{-\d_1}\big)\l v \big]
        \leq
        U_Y(v)
        \leq
        \tau_{\beta}\big[\big(1+v^{-\d_1}\big)\l v\big]
    \right\}.
\label{eqLemmaApproxUbis}
\end{equation}
If $\d_1$ is small enough, then for all large $v$,
$\P(\EG^c)\leq \frac{1}{v^{1/4-5\d_1}}$.
\end{lemma}

In the two previous lemmas, taking respectively $\frac{1}{(\log t)^{(1/2)-\d_3}}$ and
$v^{-\d_1}$ instead of simply some fixed $\tilde \e>0$ is necessary to obtain Lemma \ref{lemmaApproxLxLi}
with $\psi_\pm(r)$ instead of simply $1\pm \tilde \e$ in the definition of $\Lb_\pm(r)$ and $\Ib_\pm(r)$,
which itself is necessary for example to prove Lemma \ref{lemmaEncadrementQ}.

\bigskip
\begin{lemma}\label{lemmaApproxJCasKinf1}
Let $(\b(s),\ s\geq 0)$ be a Brownian motion, and
$\l=4(1+\k)$ as before. We define
\begin{equation}
    J_{\b}(\k, t)
:=
    \int_0^1 y(1-y)^{\k-2}L_\b\big[\tau_\b(\l), S(y)/t \big]\text{d}y,
\qquad
    0<\k\leq 1,\ t\geq 0.
\label{e1p18}
\end{equation}
Let $0<d<1$ and recall that $0<\e<1$.

\noindent {\bf (i)} Case $0<\k<1$: recall $K_\beta(\k)$ from \eqref{e3p18}.
There exist $c_{47}>0$ and
$c_{48} >0$ such that for all large $t$, on an
event $\ER$ of probability greater than $1-c_{47}/t^d$, we
have
\begin{equation}
\label{eqLemmaJcasKinf1} (1-\e)K_\b(\k)- c_{48} t^{1-1/\k}\leq
{\k}^{2-1/\k}t^{1-1/\k}J_\b(\k,t)\leq
(1+\e)K_\b(\k)+ c_{48} t^{1-1/\k}.
\end{equation}
\noindent {\bf (ii)} Case $\k=1$: recall $C_\beta$ from \eqref{e2p18}.
 There exists $c_{47}>0$ such
that for $t$ large enough, on an event $\ER$ of probability
greater than $1-c_{47}/t^d$,
\begin{equation}
\label{eqLemmaJcasKeg1} (1-\e)[C_\b + 8\log t]\leq
J_\b(1,t)\leq (1+\e)[C_\b + 8\log t].
\end{equation}
\end{lemma}

\bigskip

\begin{lemma}
 \label{remarquenegatif}
 Let $\k>0$ and define
 \begin{equation*}
     \lneg
   :=   \sup\limits_{r\geq 0} \ \sup\limits_{x<0} \lx(H(r),x) =
     \sup\limits_{t\geq 0} \ \sup\limits_{x<0} \lx(t,x),
    \qquad
     H_-(+\infty)
   :=  \lim_{r\to+\infty}H_-(r).
\end{equation*}
 There exist $c_{49}>0$ and $c_{50}>0$ such that for all large
 $z$,
 \begin{eqnarray}
     \label{eqProbaLneg}
     \P\big[ \lneg>z \big]
& \leq &
    c_{49}z^{-\k/(\k+2)},
\\
     \P\big[H_-(+\infty)>z\big]
& \leq &
    c_{50}[(\log z)/z]^{\k/(\k+2)}.
    \label{eqHmoins}
 \end{eqnarray}
\end{lemma}

\noindent{\bf Proof of Lemma \ref{remarquenegatif}:} This lemma is proved in Andreoletti et al. (\cite{Andreoletti_Devulder}, Lemma 3.5,
which is true for every $\k>0$).
More precisely, \eqref{eqHmoins} is proved in (\cite{Andreoletti_Devulder}, eq. (3.29)),
whereas \eqref{eqProbaLneg} is proved in (\cite{Andreoletti_Devulder}, eq. (3.31)).
\hfill$\Box$

\medskip

By admitting Lemmas
\ref{lemmaApproxTy}--\ref{lemmaApproxJCasKinf1}, we can now
complete the proof of Lemma \ref{lemmaApproxLxLi}.

\medskip

\noindent {\bf Proof of Lemma \ref{lemmaApproxLxLi}: part (i).}
Notice that
\begin{equation}
    S(y) \underset{y\to 1}{\sim} \int_{\a_\k}^y
    \frac{\text{d}s}{(1-s)^{1+\k}} \underset{y\to 1} {\sim}
    \frac{1}{\k} \frac{1}{(1-y)^\k}.
    \qquad
    \label{equivalent}
    \frac{y}{1-y} \underset {y\to 1} {\sim} [\k S(y)]^{1/\k}.
\end{equation}
Define
$
    \LA(r)
:=
    4\big[\sup\nolimits_{u\in[T_Y(\a_\k), \varphi(r)]} \k S(Y(u))\big]^{1/\k}
$.
We have,
\begin{equation}
    \LA(r)
 =
    4
    \Big[\sup\limits_{u\in[0,\varphi(r)-T_Y(\a_\k)]}
    \k\beta\big(U_Y(u)\big)\Big]^{1/\k}
=
    4\Big[\sup\limits_{s\in[0,U_Y(\varphi(r)-T_Y(\a_\k))]}
            \k\beta(s)
      \Big]^{1/\k} .
    \label{eqEgaliteLi}
\end{equation}
Recall $L_0$ from \eqref{e1p7}. By
\eqref{equivalent}, there exists a constant $c_{51}>0$ depending
on $\e$ such that
\begin{equation}
    \left\{ \LA(r)>c_{51} \right\} \subset \left\{ (1-\e)\LA(r)
    \le L_0(r) \le (1+\e)\LA(r) \right\} .
    \label{eqInegLxLi}
\end{equation}

We look for an estimate of $U_Y[\varphi(r)-T_Y(\a_\k)]$ appearing
in the expression of $\LA(r)$ in the right hand side of \eqref{eqEgaliteLi}. Recall
(see Dufresne~\cite{Papier_3_D1}, or Borodin et al. \cite{BorodinSalminem} IV.48 p. 78) that $A_{\infty}\el 2/\g_\k$,
where $\g_\k$ is a gamma variable of parameter $(\k,1)$, with density $e^{-x}x^{\k-1}\un_{(0,\infty)}(x)/\Gamma(\k)$. Since
$A(F(r))\le A_{\infty}$, we have
$$
    \P\big[A(F(r))>r^{2/\k}\big]
\leq
    \P\big[\g_\k<2r^{-2/\k}\big]
\leq
    2^\k r^{-2}/[\k\Gamma(\k)].
$$
On the other hand, by definition,
$
    A(F(r))
=
    A_{\infty}-\delta(r)
=
    A_{\infty}- e^{-\k r/2}
$
(see \ref{e1p5}),
which implies
$$
    \P \left[ A(F(r)) <\frac{1}{ 2\log r} \right]
=
    \P \left[\frac{2}{\g_\k} < \frac{1}{ 2\log r} + \delta(r) \right]
\leq
    \frac{1}{ r^2}
$$
for large $r$.
Consequently,
$$
    \P \left\{
        \k r/2 - 2 \log\log r
        \leq
        \log\big[\d(r)^{-1}A(F(r))\big]
        \leq
        \k r/2 +(2/\k)\log r
    \right\}
\geq
    1- c_{52}  r^{-2}.
$$
Recall that $\varphi(r)=\ty[\d(r)^{-1}A(F(r))]$ by \eqref{e10p7}.
Thus, for large $r$,
$$
    \P \left\{
        \ty\big[\exp(  \k r/2 - 2\log\log r)\big]
        \leq
        \varphi(r)
        \leq
        \ty\big[ \exp(  \k r/2 +(2/\k) \log r)\big]
    \right\}
\geq
    1- c_{52} r^{-2}.
$$
By definition,
$
    \ty (u)
=
    \int_0^u \frac{\text{d}s}{R_2^2(s) + \rb^2(s)}
$.
Notice that
$
    \big(R_2^2(t)+\rb^2(t),\ t\ge 0\big)
$
is a $(4+2\k)$--dimensional squared
Bessel process starting from $\Rt_{2+2\k}^2[\d(r)]/\d(r)$
by the additivity property of squared Bessel processes (see e.g. Revuz et al. \cite{Papier_3_RY3}, XI th. 1.2).
So,  it
follows from Lemma \ref{lemmaApproxTy}
applied with $d=4+2\k$ and $\delta_3= 1/4$,
that there exist constants
$c_{53}>0$ and $c_{54}>0$, such
that
\begin{equation}\label{gamma}
        \P \big\{
            \k  r/\l - c_{53}r^{1/2+\d_3}
            \leq
            \varphi(r)
            \leq
            \k  r/\l + c_{53}r^{1/2+\d_3}
    \big\}
\geq
    1- c_{54} r^{-2} ,
\end{equation}
for large $r$,
where $\l = 4(1+\k)$, as before.

In order to study
$T_Y(\a_\k)$, we go back to the
stochastic differential equation in \eqref{eqDefJacobi}
satisfied by the Jacobi process $Y(\cdot)$, with $d_1=2$ and
$d_2=2+2\k$. Note that $Y(s)\in(0,1)$ for any $s>0$.
By the Dubins--Schwarz theorem, there exists a
Brownian motion $\big(\Bhat(s),\ s\geq 0\big)$ such that
$$
    Y(t)
=
    \Bhat\left(4\int_0^t Y(s)(1-Y(s))\text{d} s\right)+\int_0^t
    [2-(4+2\k)Y(s)]\text{d} s,
\qquad
    t\geq 0.
$$
Recall that $\a_\k=1/(4+2\k)$, and let $t\ge 2\a_\k$.
We have, on the event $\{T_Y(\a_\k)\ge t \}$,
$$
    \inf_{0\le s\le 4 t}\Bhat(s)
\leq
    \Bhat\left(4\int_0^{t} Y(s)(1-Y(s))\text{d} s\right)
\leq
    \a_\k-t
\leq
    -\frac{t}{2},
$$
since $Y(s)\leq \a_\k\leq 1$ if
$
0\leq s \leq t
\leq T_Y(\a_\k)
$.
As a consequence, for $t\ge 2\a_\k$,
\begin{equation}
    \label{eqLemmaHittingTimeYConclusion}
    \P[T_Y(\a_\k)>t]
\leq
    \P\left(\inf_{0\le s\le 4t} \Bhat(s)\le -\frac{t}{2}\right)
=
    \P\left(\big|\Bhat(4t)\big|\geq \frac{t}{2}\right)
\leq
    2 \exp \left(-\frac{t}{32} \right).
\end{equation}
In particular,
$
\P[T_Y(\a_\k)> 64\log r ] \le 2 r^{-2}
$
for large $r$. Plug this into (\ref{gamma}),
let $c_{55}>c_{53}$
and define
$
    {\underline \varphi}
=
    {\underline \varphi} (r)
:=
    \k r/\l- c_{55} r^{1/2+\d_3}
$
and
$
    {\overline \varphi}
=
    {\overline \varphi} (r)
:=
    \k r/\l
    + c_{53}r^{1/2+\d_3}
$.
This gives,
$$
    \P \left\{ U_Y({\underline \varphi})\leq U_Y[\varphi(r)-T_Y(\a_\k)] \leq
        U_Y({\overline \varphi})
    \right\}
\geq
    1- c_{56} r^{-2}
$$
for large $r$. By Lemma \ref{lemmaApproximationU}, for small $\d_1>0$
and all large $r$,
$$
    \P\left\{
            \tau_\beta \left[ \big(1- ({\underline \varphi})^{-\d_1} \big) \l {\underline \varphi} \right]
        \leq
            U_Y[\varphi(r)-T_Y(\a_\k)]
        \leq
            \tau_\beta\left[ \big(1+  ({\overline \varphi})^{-\d_1} \big) \l {\overline \varphi}\right]
    \right\}
\geq
    1- r^{-c_{57}}.
$$
We choose $\d_1$ such that $\d_1<1/2-\d_3$. Then for
large $r$, we have
$
    \big(1- ({\underline \varphi})^{-\d_1} \big) \l {\underline \varphi}
\geq
    \big[1-2(\frac{\l}{ \k})^{\d_1}r^{-\d_1}\big] \k r
=
    \l t_-(r)
$,
and
$
    \big(1 + ({\overline \varphi})^{-\d_1} \big) \l {\overline \varphi}
\leq
    \big[1+2(\frac{\l}{ \k})^{\d_1}r^{-\d_1}\big] \k r=\l t_+(r)
$
(see \eqref{e1p22}). Thus,
\begin{equation}
    \P\left\{
        \tau_\beta[\l t_-(r)] \le U_Y[\varphi(r)-T_Y(\a_\k)]
        \le \tau_\beta[\l t_+(r)]
    \right\}
\geq
    1- r^{-c_{57}}.
    \label{eqEncadrementUti}
\end{equation}
With
$
    \Lb_\pm(r)
=
    4\big[\sup_{s\in[0,\tau_{\b}(\l t_\pm(r))]} \k\b(s)\big]^{1/\k}$ (see (\ref{e2p22})),
\eqref{eqEncadrementUti} and \eqref{eqEgaliteLi} give for large $r$,
\begin{equation}\label{eqProbaInegaliteL0}
    \P\left\{ \Lb_-(r)\leq \LA(r) \leq \Lb_+(r) \right\}
\geq
    1-r^{-c_{57}}.
\end{equation}

By \eqref{eqProbaLiInf},
$
    \P\big\{ \Lb_-(r) > r^{(1-\d_1)/\k} \big\}
\geq
    1-r^{-1}
$
for large $r$. Applying \eqref{eqInegLxLi} and \eqref{eqProbaInegaliteL0}, this yields
$$
\P\left\{ (1-\e)r^{(1-\d_1)/\k} < (1-\e)\Lb_-(r)\leq L_0(r)\leq
(1+\e)\Lb_+(r) \right\} \ge 1-  r^{-c_{58}} .
$$
Recall that
$
    \P[T_Y(\a_\k)> 64\log r ]
\leq
    2r^{-2}
$
for large $r$, which together with \eqref{gamma} gives $\P(E_{9}^c)\leq (c_{54}+2)r^{-2}$.
In view of \eqref{e3p7} and \eqref{eqMajorAvantAlphaK}, for large $r$,
\begin{equation}\label{eqAvantCalculProbaE30}
\P\left\{ (1-\e)r^{(1-\d_1)/\k} < (1-\e)\Lb_-(r)\leq
\lp[H(F(r))]\leq (1+\e) \Lb_+(r) \right\} \ge 1-
r^{-c_{59}} .
\end{equation}

On the other hand, applying Lemma \ref{remarquenegatif} to
$z=r^{(1-2\d_1)/\k}$ gives
$
    \P[ \sup_{x<0} L_X [H(F(r)),x] > r^{(1-2\d_1)/\k}]
\leq
    \P[\lneg >r^{(1-2\d_1)/\k}]
\leq
    c_{49}/r^{(1-2\d_1)/(\k+2)}
$
for large $r$. This implies
$$
    \P\left\{
        (1-\e)\Lb_-(r)\leq \lm[H(F(r))]\leq
        (1+\e)\Lb_+(r)
    \right\}
\geq
    1- \frac{1}{ r^{c_{59}}}
    -
    \frac{c_{49}}{ r^{(1-2\d_1)/(\k+2)} },
$$
proving the first part of Lemma
\ref{lemmaApproxLxLi}.\hfill$\Box$


\bigskip

\noindent{\bf Proof of Lemma \ref{lemmaApproxLxLi}: part (ii).}
In this part, we assume $0<\k\leq 1$.

\noindent Recall that
$
    H_0(r)
=
    16\int_0^{\varphi(r)-T_Y(\a_\k)}
    \frac{Y[u+T_Y(\a_\k)]}{(1- Y[u+T_Y(\a_\k)])^2} \text{d}u
$
and
$
Y[t+T_Y(\a_\k)]=S^{-1}\{\b[U_Y(t)]\}
$,
see \eqref{e1p7} and \eqref{eqDefBetaEtU}. As in Hu et al.~(\cite{Papier_3_HSY1} p.~3923, calculation of $\Upsilon_\mu$),
this and again \eqref{eqDefBetaEtU} lead to:
\begin{eqnarray}
    H_0(r)
& = &
    4\int_0^{\varphi(r)- T_Y(\a_\k)}
    (Y[u+T_Y(\a_\k)])^2
    (1-Y[u+T_Y(\a_\k)])^{2\k-1}\text{d} U_Y(u)
\nonumber\\
& = &
    4\int_0^1 x (1-x)^{\k-2} L_{\beta} [U_Y(\varphi(r)- T_Y(\a_\k)), S(x)] \text{d} x.
\label{eqFormuleH0}
\end{eqnarray}
Recall that
$
    t_\pm(r)
=
    \big[1\pm 2 (\frac{\l}{  \k})^{\delta_1}r^{-\d_1} \big] \frac{\k}{ \l} r
$,
$\b_v(s)=\b(v^2s)/v$
and let $J_{\b}$ be as in \eqref{e1p18}. We have,
\begin{eqnarray*}
\int_0^1  x (1-x)^{\k-2} L_{\beta}\{ \tau_{\beta} [\l
    t_\pm(r)], S(x) \} \text{d}x
 & = & t_\pm(r) \int_0^1  x (1-x)^{\k-2} L_{\beta_{t_\pm(r)}}\{
    \tau_{\beta_{t_\pm(r)}} (\l), \frac{S(x)}{t_\pm(r)} \} \text{d}x\\
&  = & t_\pm(r) J_{\beta_{t_\pm(r)}} [\k,t_\pm(r)].
\end{eqnarray*}
By \eqref{eqEncadrementUti} and \eqref{eqFormuleH0}, we have for large $r$,
\begin{equation*}
    \P \left[4t_-(r) J_{\beta_{t_-(r)}} [\k,t_-(r)]
      \le H_0 (r) \le 4 t_+(r) J_{\beta_{t_+(r)}} [\k,t_+(r)] \right] \ge
    1- r^{-c_{57}}.
\end{equation*}
Now, apply Lemma \ref{lemmaApproxJCasKinf1} to
$d=1/2$. So there exist $c_{6}>0$ and $c_{60}>0$ such that for large $r$,
\begin{equation}
    \P \left\{ (1-\e)\Ib_- (r)\leq H_0(r)\leq (1+\e)\Ib_+ (r) \right\} \ge
    1-
    r^{-c_{60}} ,
    \label{eqI''CasGeneral}
\end{equation}
where $\Ib_\pm (r)$ is defined in (\ref{eqI''CasKinf1}).

In the case $0<\k<1$, we know that $\P(E_{9}^c)\leq (c_{54}+2)r^{-2}$ for large $r$ as proved before \eqref{eqAvantCalculProbaE30},
so by \eqref{eqMajorAvantAlphaK},
$
    \P\big[\hbar(r)\le c_{61}\log r\big]
\geq
    \P(E_{9})
\geq
    1-(c_{54}+2)r^{-2}
$
for some $c_{61}$ and all large $r$.
On the other hand, by Lemma
\ref{remarquenegatif}, $\P[H_-(F(r))\le \e r] \ge
\P[H_-(+\infty)\le \e r] \ge 1- {c_{62}
r^{-(1-\d_1)\k/(\k+2)}}$, for all large $r$. Consequently, by
\eqref{eqI''CasGeneral} and \eqref{e3p7}, for large $r$,
$$
    \P\big\{
        (1-\e)\Ib_-(r)\leq H(F(r))\leq (1+\e)\Ib_+(r)+
        (4\e\l/\k) t_+(r)
    \big\}
\geq
    1-r^{-c_{63}}.
$$
Changing the value of $c_{6}$, this proves Lemma \ref{lemmaApproxLxLi} (ii) in the
case $0<\k<1$.

Now we consider the case $\k=1$. As before,
$
    \P[H_-(F(r))+\hbar(r)\leq 2\e r]
\geq
    1-r^{-c_{64}}
$
(for large $r$).
Moreover,
$
    \P\big[C_{\b_{t_\pm(r)}}>-\pi\log r\big]
\geq
    1-r^{-2}
$
by Fact \ref{FactBianeYor} and \eqref{ProbaSamoTaqqu}.
Therefore, \eqref{eqI''CasKinf1} gives
$
    \P\big\{\Ib_+(r)\geq 16t_+(r)\log r\big\}
\geq
    1-r^{-2}
$.
Consequently, for large $r$,
$$
    \P\big[0\le H_-(F(r))+\hbar(r) \le \e \Ib_+(r)\big]
\ge
    1- r^{-c_{65}} ,
$$
which, in view of \eqref{eqI''CasGeneral}, yields
that, for large $r$,
$$
    \P\big[ (1-\e)\Ib_-(r)\le H(F(r))\le (1+2\e)\Ib_+(r)\big]
\geq
    1-r^{-c_{66}} .
$$
This proves Lemma \ref{lemmaApproxLxLi} (ii) in the
case $\k=1$. \hfill$\Box$


\mysection{Proof of lemmas \ref{lemmaApproxTy}--\ref{lemmaApproxJCasKinf1}}\label{SectProofLemmasFin}

This section is devoted to the proof of Lemmas \ref{lemmaApproxTy}--\ref{lemmaApproxJCasKinf1}. For the sake of
clarity, the proofs of these lemmas are presented in separated
subsections.


\subsection{Proof of Lemma \ref{lemmaApproxTy}}

First, notice that we can not apply Talet (\cite{Papier_3_Ta2}, Lem. 3.2 eq. (3.4))
since her constant $c_3$ depends on her (fixed) $\delta$, whereas we would like to take
her $\delta=(\log t)^{\delta_3-1/2}\to_{t\to+\infty} 0$, which is necessary for example for our Lemma \ref{lemmaEncadrementQ}.
A similar remark applies for Talet (\cite{Papier_3_Ta2},  Prop. 5.1) and our Lemma \ref{lemmaApproximationU}.
So we need different estimates than in her paper.

Let $d>4$ and $R_0\el\widetilde{R}_{d-2}(1)$, where
$\widetilde{R}_{d-2}$ is a $(d-2)$--dimensional Bessel process. We
consider a $d$--dimensional Bessel process $R$, starting from
$R_0$.
We introduce
$\theta(t):=\int_0^t R^{-2}(s) \text{d}s$. It\^o's
formula gives $\log R(t)=\log R_0
+M(t)+\frac{d-2}{2}\theta(t)$, where $M(t):=\int_0^t
R(s)^{-1}\text{d}\bhat(s)$ and $\big(\bhat(t),\ t\geq 0\big)$ is a
Brownian motion. By the Dubins--Schwarz theorem, there exists
a Brownian motion $\big(\Bt(t),\ t\geq 0\big)$ such that
$M(t)=\Bt(\theta(t))$ for all $t\geq 0$. Accordingly,
\begin{equation}
    \label{ito}
    (d-2)\theta(t)/2 = \log R(t) - \log R_0
    - \Bt(\theta(t)), \qquad t\geq 0.
\end{equation}

Let $\delta_3 \in (0, \frac{1}{ 2})$, $0<\e<1$, and
$
    x
=
    x(t)
:=
    {\frac{d-2}{ 6}} \frac{1}{ (\log t)^{(1/2)-\delta_3}}
$.
We have (see e.g. G{\"o}ing-Jaeschke et al. \cite{Going-Yor}, eq. (50)),
$
    \P\big(R_0^2\in\text{d}u\big)
=
    u^{d/2-2}e^{-u/2}\un_{(0,\infty)}(u)/\big[\Gamma(d/2-1)2^{d/2-1}\big]
$.
So for large $t$,
%
%
%
\begin{equation}
    \label{eqTpuissanceEpsilon1}
    \P\left(\left|\frac{\log R_0}{\log t}\right|> x\right)
=
    \P\left( \frac{\log R_0}{\log t}> x \right)
    +\P\left( \frac{\log R_0}{\log t} < -x \right)
\leq
    \exp\left( - (1-\varepsilon) \frac{t^{2x}}{ 2}\right)
    + \frac{c_{67}}{t^{x(d/2-1)}} .
\end{equation}

Denote by  $n:= \lceil d\rceil$ the smallest integer such that $n\ge
d$. Since an $n$-dimensional Bessel process can be realized as the
Euclidean modulus of an $\R^n$-valued Brownian motion, it follows
from the triangular inequality that
$
R(t)\leq_{\mathcal{L}}
 R_0 +\widehat{R}_n(t)
$,
where $(\widehat{R}_n(t),\ t\geq 0)$ is an
$n$-dimensional Bessel process starting from $0$.
Consequently, for large $t$,
$
    \P\big(R(t)>t^{(1/2)+x}\big)
\leq
    \P\big(\widehat{R}_n(t)> t^{(1/2)+x}/2 \big)
+
    \P\big(R_0> t^{x} \big)
\le
    2\exp( - (1- \e){t^{2x}/8}) ,
$
and
$
    \P\big(R(t)<t^{(1/2)-x} \big)
\leq
    2 t^{-x}
$, e.g. since $\big|\tilde \beta(t)\big|\leq_{\mathcal{L}} R(t)$.
Therefore, for large $t$,
\begin{equation}
    \label{eqTpuissanceEpsilon11}
    \P\left(\left|\frac{\log R(t)}{\log t}-\frac{1}{2}\right|> x \right)
\leq
    2\exp\left( - (1- \e) \frac{t^{2x}}{ 8}\right) +
    2 t^{-x} .
\end{equation}

Define
$
    \EK
:=
    \left\{ \left|\frac{\log R(t)}{\log t} - \frac{1}{2}
    \right| \le x \right\} \cap \left\{ \left|
    \frac{\log R_0}{\log t} \right| \le x \right\}
$
and
\begin{equation*}
    \ET
:=
    \left\{ \frac{d-2}{2} \theta(t) < \log t\right\},
\qquad
    \EY
:=
    \bigg\{ \sup_{0\le s \leq 2(\log t)/(d-2)} \big|\Bt(s)\big| \leq
    x\log t\bigg\}.
\end{equation*}
By \eqref{eqTpuissanceEpsilon1} and
\eqref{eqTpuissanceEpsilon11}, we have for large $t$,
\begin{equation}
\label{Maine10}
    \P(\EK^c)
\le
    3\exp\big[ - (1- \e) t^{2x}/8\big]
    +
    3 t^{-x} .
\end{equation}
We now estimate $\P(\EK\cap\ET^c)$. We first observe that on
$\EK$, we have, by \eqref{ito},
$$
    \big| \Bt(\theta(t)) + (d-2) \theta(t)/2- (\log t)/2 \big|
\le
    2x\log t.
$$
 We claim that
$
    \EK\cap\ET^c
\subset
    \big\{\big|\Bt(\theta(t))\big|>\frac{d-2}{6} \theta(t)\big\}
$
for large $t$.
Indeed, on the event $\EK\cap\ET^c \cap \big\{\big|\Bt(\theta(t))\big| \le \frac{d-2}{6} \theta(t)\big\}$,
$$
(d-2)\theta(t)/2 \leq (2x+ 1/2 ) \log
t- \Bt(\theta(t)) \le (2x+1/2) \log t +
(d-2) \theta(t)/6 ,
$$
which implies $\frac{d-2}{2}\theta(t) \le
(\frac{3}{4}+3x)\log t$. This, for large $t$,
contradicts $\frac{d-2}{2}\theta(t)\geq \log t$ on $\ET^c$.
Therefore,
$
    \EK\cap\ET^c
\subset
    \big\{\big|\Bt(\theta(t))\big|>\frac{d-2}{6} \theta(t)\big\}
$ holds for all
large $t$, from which it follows that
\begin{eqnarray*}
    \P\big(\EK\cap\ET^c\big)
& \le &
    \P\bigg( \, \sup\limits_{s\ge 2(\log t)/(d-2)}
    \frac{\big|\Bt(s)\big|}{s} > \frac{d-2}{6} \bigg)
=
    \P\bigg( \sup_{u \ge 1} \frac{\big|\Bt(u)\big|}{u}
    >
    \sqrt{\frac{(d-2)\log t}{18} }
    \bigg)
    \\
 &=&\P\left( \sup_{0\le v\le 1} \big|\Bt(v)\big| >
    \sqrt{(d-2)(\log t)/18}
    \right)
\le
    4\exp\big[ -  (d-2)(\log t)/36\big],
\end{eqnarray*}
because $u\mapsto  u \Bt(1/u)$ is a Brownian motion
and
$
    \sup_{0\le v\le 1} \Bt(v)
\el
    \big|\Bt(1)\big|
$.
Since
$
    \P\big(\EY^c\big)
\leq
    4\exp\big[- \frac{d-2}{4} x^2 \log t\big]
$
for large $t$, this and \eqref{Maine10} give for large $t$,
\begin{equation*}
    \P\big(\EK^c \cup \ET^c \cup \EY^c\big)
\le
    \P\big(\EK^c\big)+\P\big(\EK\cap\ET^c\big)+P\big(\EK\cap\ET\cap\EY^c\big)
\leq
    \exp\big(- c_{68} \, x^2 \log t\big).
\end{equation*}
Since
$
    \EK\cap \ET \cap \EY
\subset
    \big\{ \big| \frac{\theta (t)}{ \log t} - \frac{1}{ d-2}\big| \le \frac{6x}{ d-2}\big\}
$
by \eqref{ito},
this completes the proof of
Lemma \ref{lemmaApproxTy}.\hfill$\Box$


\subsection{Proof of Lemma \ref{lemmaApproximationU}}

Let $v>0$. Recall that for every $x\geq 0$, $\b_v(x) = (1/v)\b(v^2x)$,
and notice that $v^2
\tau_{\beta_v}(x)=\tau_{\beta}(x v)$ almost surely. Then,
\begin{equation}
    \EG
=
    \big\{\tau_{\beta_v}\big[\big(1-v^{-\d_1}\big)\l\big]\leq U_Y(v)/v^2\leq
    \tau_{\beta_v}\big[\big(1+v^{-\d_1}\big)\l\big]\big\}.
\label{eqLemmaApproxU}
\end{equation}
For $\d_1>0$, define
$
    \EGb
:=
    \big\{\sup\nolimits_{0\leq s\leq
    \tau_{\beta_v}(2\l)}|\e_1(v,s)|<v^{-\d_1}\big\}
$, where
\begin{equation*}
    \e_1
=
    \e_1(v,s)
:=
    \frac{1}{4}
    \int_0^1 (1-x)^{\k}
    \left[
        L_{\beta_v}\bigg(s, \frac{S(x)}{v}\bigg)-L_{\beta_v}(s,0)
    \right]\text{d}x,
\qquad
    s\geq 0.
\end{equation*}
By Hu et al. (\cite{Papier_3_HSY1} eq. (2.34) p.~3924),
$\EGb\subset \EG$. Thus it remains to prove that for $\d_1$ small
enough, $\P(\EGb^c)\leq 1/v^{1/4-5\d_1}$ for large $v$. Notice
that for $s\geq 0$,
\begin{eqnarray}
    |\e_1|
& \leq &
    \bigg(
    \int_{\{S(x)>\sqrt{v}\}}+\int_{\{S(x)<-\sqrt{v}\}}+
    \int_{\{|S(x)|\leq\sqrt{v}\}} \bigg)
    \frac{(1-x)^{\k}}{4}\left|L_{\beta_v}
    \bigg(s, \frac{S(x)}{ v}\bigg)-L_{\beta_v}(s,0)\right|\text{d}x
\nonumber\\
&=:&
    \e_2(v,s)+\e_3(v,s)+\e_4(v,s).
\label{voila2}
\end{eqnarray}


Since $S(y)=\int_{\a_{\k}}^y\frac{\text{d}x}{x(1-x)^{1+\k}}$, we have
$
1-S^{-1}(u)\underset{u\to+\infty}{\sim}(\k u)^{-1/\k}
$.
So, 
we have
\begin{eqnarray*}
    \sup_{0\leq s \leq \tau_{\beta_v}(2\l)}\e_2(v,s)
& \leq &
    \frac{1}{4}\int_{1-\left(\frac{2}{\k\sqrt{v}}\right)^{1/\k}}^1 \
    (1-x)^{\k}\sup_{0\leq s \leq \tau_{\beta_v}(2\l)}
    \sup\limits_{u\geq
    0}\left[L_{\beta_v}(s,u)+L_{\beta_v}(s,0)\right]\text{d}x
    \\
& \leq &
    [2/(\k\sqrt{v})]^{\frac{1}{\k}+1}
    \sup\nolimits_{u\geq 0}
    \left[L_{\beta_v}(\tau_{\beta_v}(2\l),u)+2\l\right],
\end{eqnarray*}
%
for all large $v$.
By the second Ray--Knight theorem (Fact
\ref{FactRN2}),
$Q:=(L_{\beta_v}(\tau_{\beta_v}(2\l),u), \ u\geq 0)$ is a
$0$--dimensional squared Bessel process starting from $2\l$.
Moreover, $x\mapsto x$ is a scale function of $Q$ (see e.g. Revuz et al. \cite{Papier_3_RY3} p. 442).
Hence, for large $v$,
\begin{equation}
    \P\bigg(\sup_{0\leq s \leq \tau_{\beta_v}(2\l)}\e_2(v,s)
            \geq
            \left[2/(\k\sqrt{v})\right]^{1/\k+1}\sqrt{v}
    \bigg)
\leq
    \P\left(\sup_{u\geq 0}Q(u)\geq \frac{\sqrt{v}}{2}\right)
=
    \frac{4\l}{\sqrt{v}}.\label{eqInegEpsilon2}
\end{equation}
Similarly (this time, using $S(x)\sim\log x$, $x\to 0$), we
have, for large $v$,
\begin{equation}
\label{eqInegEpsilon3}
\P\left[\sup\nolimits_{0\leq s \leq
\tau_{\beta_v}(2\l)}\e_3(v,s)\geq
\exp(-\sqrt{v}/2) \sqrt{v}\right] \leq
4\l/\sqrt{v}.
\end{equation}

To estimate $\e_4(v,s)$, we note that
\begin{equation}
    \e_4(v,s)
\leq
    \sup\limits_{|u| \le 1/\sqrt{v}}
    \big|L_{\beta_v}(s,u)-L_{\beta_v}(s,0)\big|.
\label{eqInegEpsilon4}
\end{equation}
Let $\e\in(0,1/2)$, $t_v>0$, $\g\geq 1$ and define
$(M)_t^*:=\sup_{0\leq s\leq t}|M(s)|$ for $t>0$ and any Brownian motion $(M(s),\ s\geq 0)$.
Applying Barlow and Yor (\cite{Papier_3_BY1}, (ii) p. 199) to the
continuous martingale $\beta_v(.\wedge t_v)$
and its jointly continuous local time $(L_{\beta_v}(u\wedge t_v, a), \ u\geq 0,\ a\in\R)$,
we see that for
some constant $C_{\g,\e}>0$,
$$
    \left\|\sup\limits_{0\leq s\leq t_v,a\neq b}
           \frac{\big|L_{\beta_v}(s,b)-L_{\beta_v}(s,a)\big|}{|b-a|^{1/2-\e}}
    \right\|_{\gamma}
\leq
    C_{\g,\e}\left\|[(\beta_v)_{t_v}^*]^{1/2+\e}\right\|_{\g}
=
    C_{\g,\e}\left\|[(\beta)_{t_v}^*]^{1/2+\e}\right\|_{\g},
$$
where $\|.\|_{\g}=\E(|.|^\g)^{1/\g}$. Then, by Chebyshev's inequality and a change of scale, for $\a>0$,
\begin{equation}
    \P\left(
        \sup\limits_{0\leq s\leq t_v, \, a\neq b}
        \frac{\big|L_{\beta_v}(s,b)-L_{\beta_v}(s,a)\big|}{|b-a|^{1/2-\e}}
        \geq
        \a
    \right)
\le
    \frac{(\sqrt{t_v})^{(1/2+\e)\g}}{\a^{\g}}
    \left[C_{\g,\e}\left\| [(\beta)_{1}^*]^{1/2+\e}
    \right\|_{\g}\right]^\g.
\label{eqApplicationBarlowYor}
\end{equation}
On
$
    \EH
:=
    \Big\{ \sup\nolimits_{0\leq s\leq \tau_{\beta_v}(2\l),a\neq b}
            \frac{\left|L_{\beta_v}(s,b)-L_{\beta_v}(s,a)\right|}{|b-a|^{1/2-\e}}
            \leq
            v^{\frac{1}{2}\left(\frac{1}{2}-2\e\right)}
    \Big\}
$,
we have by \eqref{eqInegEpsilon4},
\begin{eqnarray}
    \sup\limits_{0\leq s \leq \tau_{\beta_v}(2\l)}\e_4(v,s)
& \leq &
    v^{-\frac{1}{2}(\frac{1}{2}-\e)}v^{\frac{1}{2}
    \left(\frac{1}{2}-2\e\right)}=v^{-\e/2}.
\label{eqInegEpsilon4surEH}
\end{eqnarray}
We now choose $\g:=2$ and
$t_v:=v^{\frac{1/4-\e}{1/2+\e}}$.
Since
$
    \P[\tau_{\beta_v}(2\l)>t_v]
=
    \P[L_{\beta_v}(t_v,0)<2\l]
=
    \P\big[\sup_{0\leq s \leq t_v}\beta(s)<2\l\big]
=
    \P[|\beta(t_v)|<2\l]
\leq
    4\l/\sqrt{t_v}
$
by L\'evy's theorem (see e.g. Revuz et al. \cite{Papier_3_RY3} VI th. 2.3),
we get for all large $v$ (if $\e$ is small enough).
\begin{eqnarray*}
    \P\big[\EH(v)^c\big]
& \leq &
    \P\big[\tau_{\beta_v}(2\l)>t_v\big]
    +
    \P\left(\sup\limits_{0\leq s\leq t_v,a\neq b}
    \cfrac{|L_{\beta_v}(s,b)-L_{\beta_v}(s,a)|}{|b-a|^{1/2-\e}} \geq
    v^{\frac{1}{2}\left(\frac{1}{2}-2\e\right)}  \right)
\\
& \leq &
    4\l\ v^{\frac{\e-1/4}{1+2\e}}
    +c_{69}v^{-1/4+\e}
\leq
   v^{-1/4+2\e}/2.
\end{eqnarray*}
Combining  this with  \eqref{voila2}, \eqref{eqInegEpsilon2},
\eqref{eqInegEpsilon3} and \eqref{eqInegEpsilon4surEH}, we obtain
that, for $\e>0$ small enough,
$$
\P\left(\sup\nolimits_{0\leq s\leq
\tau_{\beta_v}(2\l)}|\e_1(v,s)|\geq 2v^{-\e/2} \right)\leq
v^{-1/4+2\e}.
$$
This gives, with the choice of $\d_1:=2\e/5$,
$
    \P(\EG^c)
\leq
    \P(\EGb^c)
\leq
    v^{-1/4+5\d_1}
$
(for large $v$).\hfill$\Box$

\subsection{Proof of Lemma \ref{lemmaApproxJCasKinf1}}

Assume $0<\k\leq 1$. Consider $0<d<1$,
 $\e\in(0,1/2)$ such that $d(1/2+\e)+(\e-1)(1/2-\e)<0$, $M_\e>0$,  and a Brownian motion $(\beta(t),\ t\geq 0)$.
We can write for $t>0$,
\begin{eqnarray*}
    J_{\b}(\k,t)
& = &
    \bigg(\int_{0}^{S^{-1}(-t^{\e})}
          +\int_{S^{-1}(-t^{\e})}^{\a_\k}
          +\int_{\a_\k}^{S^{-1}(M_\e)}
          +\int_{S^{-1}(M_\e)}^{1}
    \bigg)
          y(1-y)^{\k-2} L_\b\bigg(\tau_\b(\l),\frac{S(y)}{t}\bigg)\text{d}y
\\
& := &
    J_1+J_2+J_3+J_4.
\end{eqnarray*}

We begin by estimating $J_1$.
Since $S(x)\sim_{x\to 0}\log x$,  we have
$
J_1 \leq
\exp(-t^{\e}/2)\sup_{s\geq 0}
\tilde Q(s)
$
for large $t$,
where $\tilde Q$ is a $0$--dimensional squared Bessel process starting
from $\l$ (by the second Ray--Knight theorem stated in Fact
\ref{FactRN2}, applied to $-\beta$). Hence,
we get
$
\P\left[J_1 \geq \exp(-t^\e/2)t^d \right]
\leq \l/t^d.
$

\noindent Fix a large constant $\g>0$ such that $d(1/2+\e+1/\g)+(\e-1)(1/2-\e)<0$, and define
\begin{equation*}
    \EP
:=
    \big\{\tau_\b(\l)\leq t^{2d}\big\},
\qquad
    \EQ
:=
    \left\{\sup_{0\leq s\leq t^{2d}, \, a\neq b}
        \frac{|L_\b(s,b)-L_\b(s,a)|}{|b-a|^{1/2-\e}}\leq t^{d(1/2+\e+1/\g)}
    \right\}.
\end{equation*}
Recall that $S(\a_\k)=0$.
To estimate $J_2$, we note that, on $\EP\cap \EQ$, uniformly  for all
large $t$,
\begin{equation*}
    J_2
\leq
    \left[\int_0^{\a_k}\hspace{-3mm}
          \frac{y\text{d}y}{(1-y)^{2-\k}}
    \right]
    \sup_{-t^{\e-1}\leq b\leq 0} L_\b(\tau_\b(\l),b)
\leq
    \frac{\a_\k\big[\l+ t^{d(1/2+\e+1/\g)}
    (t^{\e-1})^{\frac{1}{2}-\e}\big]}{(1-\a_k)^{2-\k}}
\leq
    \frac{2\a_\k\l}{(1-\a_k)^{2-\k}} .
\end{equation*}
Notice that $\P(\EP^c)\leq 2\l t^{-d}$ as proved after \eqref{eqInegEpsilon4surEH},
and that $\P(\EQ^c)\leq c_{70} t^{-d}$ (by \eqref{eqApplicationBarlowYor} with $t^{2d}$ instead of $t_v$).
Therefore, there exists $c_{71}>0$ such that for large $t$,
\begin{equation}
\label{eqProbaJ2}
    \P\left( J_2 \le c_{71} \right)
\geq
    \P(\EP\cap \EQ)
\geq
    1 - c_{72} t^{-d}.
\end{equation}

We now turn to $J_3$.
As already noticed after \eqref{voila2}, we have
$
1-S^{-1}(u)\underset{u\to+\infty}{\sim}(\k u)^{-1/\k}
$.
Therefore, we can choose $M_\e>0$ such that
\begin{equation}\label{eqUtiliteMepsilon}
    \forall u\geq M_\e,
\qquad
    \frac{[1-S^{-1}(u)]^{2\k-1}}{(\k u)^{1/\k-2}}
\in
    (1-\e,1+\e)
\quad \text{and}\quad
    S^{-1}(u)
\geq
    1-\e.
\end{equation}
On the event $\EP\cap\EQ$, uniformly for all large $t$,
\begin{eqnarray*}
    J_3
& \leq &
    \sup\limits_{0\leq x\leq M_\e/t} L_\b(\tau_\b(\l),x)
    \int_{\a_\k}^{S^{-1}(M_\e)}y(1-y)^{\k-2}\text{d}y
\\
& \leq &
    c_{73}
    \left[\l+t^{d(1/2+\e+1/\g)}(M_\e/t)^{\frac{1}{2}-\e}\right]
\leq
    2\l
    c_{73}
.
\end{eqnarray*}
Consequently,
$
    \P[J_3\leq 2\l c_{73}]
\geq
    \P(\EP \cap \EQ)
\geq
    1-c_{72}t^{-d}
$ for large $t$.

Now we write
$$
    J_4
=
    \k^{1/\k-2}\ t^{1/\k-1}
    \int_{M_\e/t}^{+\infty}\left[S^{-1}(t x)\right]^2
                           \frac{\left[1-S^{-1}(t x)\right]^{2\k-1}}{(\k t)^{1/\k-2}}
                           L_\b(\tau_\b(\l),x)\text{d}x.
$$
Therefore, \eqref{eqUtiliteMepsilon} leads to
\begin{equation}
(1-\e)^3\int_{M_\e/t}^{+\infty}x^{1/\k-2}L_\b(\tau_\b(\l),x)\text{d}x
 \leq   \k^{2-1/\k}t^{1-1/\k}J_4
   \leq
(1+\e)\int_{M_\e/t}^{+\infty}x^{1/\k-2}L_\b(\tau_\b(\l),x)\text{d}x.
\label{eqApproxJ3a}
\end{equation}

\bigskip\noindent
{\bf Proof of Lemma \ref{lemmaApproxJCasKinf1}: part (i).} We
first assume $0<\k<1$.

On $\EP\cap\EQ$, for large $t$, we have
$\int_0^{M_\e/t}x^{1/\k-2}L_\b(\tau_\b(\l),x)\text{d}x\leq
c_{74} t^{1-1/\k}$. Recall $K_\beta$ from \eqref{e3p18}.
It follows from \eqref{eqApproxJ3a} and \eqref{eqProbaJ2} that,
for large $t$,
$$
    \P\left[ (1-\e)^3K_\b(\k)-(1-\e)^3c_{74} t^{1-1/\k}
             \le
             \k^{2-1/\k}t^{1-1/\k} J_4
             \le
             (1+\e)K_\b(\k)
             \right]
\ge
    1- c_{72} t^{-d}.
$$
Since $J_{\b}(\k,t) = J_1+J_2+J_3+J_4$,
we get for large $t$,
$$
    \P\left\{
        (1-\e)^3 K_\b(\k)-c_{48}t^{1-1/\k}
        \leq
        {\k}^{2-1/\k}t^{1-1/\k}J_\b(\k,t)
        \leq
        (1+\e)K_\b(\k)+c_{48}t^{1-1/\k}
    \right\}
\geq
    1-c_{75} t^{-d},
$$
for some $c_{48}>0$,
proving the lemma in the case $0<\k<1$.\hfill$\Box$

\bigskip\noindent
{\bf Proof of Lemma \ref{lemmaApproxJCasKinf1}: part (ii).} We
assume $\k=1$, thus $\l=8$.

By the definition of $C_\b$ (see \eqref{e2p18}), we have
\begin{equation*}
\int_{M_\e/t}^{\infty}\frac{L_\b(\tau_{\beta}(8),x)}{x}\text{d}x
  =  C_\b -\int_0^{M_\e/t}
\frac{L_\b(\tau_{\beta}(8),x)-8}{x}\text{d}x+8\log t-8\log
    M_\e.
\end{equation*}
On $\EP\cap\EQ$, for large $t$,
$$
\int_0^{M_e/t} \frac{|L_\b(\tau_\b(8),x)-8|}{x}\text{d}x  \leq
 \int_0^{M_e/t} \frac{t^{d(1/2+\e+1/\g)}x^{1/2-\e}}{x}\text{d}x
\leq\e.
$$
As in \eqref{ProbaSamoTaqqu}, $\P(C_\b+8\log t <\log t )\leq t^{-7}$. Therefore, by  (\ref{eqApproxJ3a}) and \eqref{eqProbaJ2}, we
have for large $t$,
$$
    \P\left\{ (1-\e)^4 [C_\b +8\log t]\le J_4 \le (1+\e)^2 [C_\b
              +8\log t]
      \right\}
\ge 1- c_{76} t^{-d} .
$$
Since $J_{\b}(1,t) = J_1+J_2+J_3+J_4$,
we get for large $t$,
$$
\P\left\{ (1-\e)^4 [C_\b +8\log t] \le J_\b(1,t)\le
(1+\e)^3 [C_\b +8\log t] \right\} \ge 1- c_{77} t^{-d}.
$$
This proves the lemma in the case $\k=1$.
\hfill$\Box$


\medskip
\noindent {\bf  Acknowledgements}

\smallskip
I would like to thank Zhan Shi for many helpful discussions.
I am also grateful to an anonymous referee for a very careful reading of the first version of the paper and for valuable
comments which helped improve the presentation of the paper, including a recommendation to merge the first version of this paper
with \cite{DevulderPreprint}.


\end{document}